\documentclass[12pt,reqno]{amsart}
\usepackage[utf8x]{inputenc}
\usepackage{amssymb,amsmath,latexsym}
\usepackage{mathrsfs}
\usepackage{epsfig}
\usepackage{color}
\usepackage{hyperref}

\usepackage{epic,eepic,wrapfig,color, ifthen}
\usepackage{url}
\newcommand{\eps}{\varepsilon}
\newcommand{\vp}{\varphi}
 \newcommand{\Rd}{{\R^{d}}}
\newtheorem{thm}{Theorem}[section]
\newtheorem{cor}[thm]{Corollary}
\newtheorem{lemma}[thm]{Lemma}
\newtheorem{prop}[thm]{Proposition}
\newtheorem{defn}[thm]{Definition}
\newtheorem{definition}[thm]{Definition}

\newtheorem{remark}[thm]{Remark}
\newtheorem{example}[thm]{Example}
\numberwithin{equation}{section}
\renewcommand{\H}{\mathbb{H}}
\newcommand{\formula}[2][nolabel]
{\ifthenelse{\equal{#1}{nolabel}}
 {\begin{align*} #2 \end{align*}}
 {\ifthenelse{\equal{#1}{}}
  {\begin{align} #2 \end{align}}
  {\begin{align} \label{#1} #2 \end{align}}
 }
 
}
\def\({\left(} \def\){\right)} 
\def\pf{{\medskip\noindent {\bf Proof. }}}
\def\qed{{\hfill $\Box$ \bigskip}}
\def\R{{\mathbb R}}
\def\P{{\mathbb P}}
\def\E{{\mathbb E}}
\def\1{{\bf 1}}

\def\WUS{\text{WUS}}
\def\WLS{\text{WLS}}
\def\WS{\text{WS}}
\def\diag{\text{diag}}
\def\lo{\text{Log}}
\def\scL{{\mathscr L}}

\newcommand{\gener}{{\cal L}}

\newcommand{\ip}[1]{\langle #1 \rangle}
\DeclareMathOperator{\dist}{dist}

\newcommand{\cal}[1]{\mathcal{#1}}
\def\sA {{\cal A}} \def\sB {{\cal B}} \def\sC {{\cal C}}

  \def\sL {{\cal L}}
  \def\sO {{\cal O}}
  
  \def\sU {{\cal U}}
 \def\sW {{\cal W}}

 \def\bE {{\mathbb E}}

\def\bP {{\mathbb P}}  \def\bR {{\mathbb R}}

\def\R {{\mathbb R}}
\def\nn{\nonumber}
\def\ll {{\ell^*}}
\def\tll {{\wh{\ell}}}
\def\tPhi {{\wh{\Phi}}}

\def\wt{\widetilde}
\def\wh{\widehat}

\def\E{{\mathbb E}}
\def\P{{\mathbb P}}
\def\diam{\text{diam}}
\def\dist{\text{dist}}

\def\bea{\begin{align*}}
\def\eea{\end{align*}}
\def\bee{\begin{equation}}
\def\eee{\end{equation}}

\makeatletter
\@addtoreset{equation}{section}

\makeatother

\begin{document}
\allowdisplaybreaks

\title[Estimates of Dirichlet heat kernels for unimodal L\'evy processes]
{ \bf Estimates of Dirichlet heat kernels for unimodal L\'evy processes with low intensity of small jumps}

\author
{Soobin Cho}
\address[Cho]{Department of Mathematical Sciences,
Seoul National University,
Building 27, 1 Gwanak-ro, Gwanak-gu
Seoul 08826, Republic of Korea}
\thanks{The research of Soobin Cho is supported by the POSCO Science Fellowship of POSCO TJ Park Foundation.}
\curraddr{}
\email{soobin15@snu.ac.kr}

\author{Jaehoon Kang
}
\address[Kang]{Department of Mathematical Sciences, KAIST,  291 Daehak-ro,
	Yuseong-gu, Daejeon 34141, Republic of Korea}
\thanks{The research of Jaehoon Kang is supported by the National Research Foundation of Korea(NRF) grant funded by the Korea government(MSIT) (No. NRF-2019R1A5A1028324).}
\curraddr{}
\email{jaehoon.kang@kaist.ac.kr}

\author{Panki Kim}
\address[Kim]{Department of Mathematical Sciences and Research Institute of Mathematics,
Seoul National University,
Building 27, 1 Gwanak-ro, Gwanak-gu
Seoul 08826, Republic of Korea}
\thanks{The research of Panki Kim is supported by
 the National Research Foundation of Korea (NRF) grant funded by the Korea government (MSIP) : NRF-2016K2A9A2A13003815.}
\curraddr{}
\email{pkim@snu.ac.kr}

 \date{}

\maketitle

\begin{abstract}
In this paper, we study transition density functions for pure jump unimodal L\'evy processes killed upon leaving an open set $D$. Under some mild assumptions on the L\'evy density, we establish two-sided Dirichlet heat kernel estimates when the open set $D$ is $C^{1, 1}$.
Our result covers the case that the L\'evy densities of unimodal L\'evy processes are   regularly varying functions whose indices are equal to the Euclidean dimension. 
This is the first results on two-sided Dirichlet heat kernel estimates for L\'evy processes such that the lower scaling
 index of the L\'evy densities  is not necessarily strictly bigger than the Euclidean  dimension. 
\end{abstract}



\bigskip

\section{Introduction and Main results}

\subsection{Introduction}

The transition density of a L\'evy process killed  upon leaving  an open set $D$ (called the Dirichlet heat kernel of the process in $D$) is the  fundamental solution of the  equation $\partial_t u = \sL u$ with zero exterior condition on $D^c$, where $\sL$ is the infinitesimal generator of the L\'evy process.
 When the sample paths of the L\'evy process are discontinuous, such generator is a non-local operator. Hence,  transition densities of killed L\'evy processes with jumps play significant roles in the study of non-local operators with killings. However,  except for a few special cases, it is impossible to find an explicit expression
of the Dirichlet heat kernel.
Thus, obtaining  sharp two-sided estimates on the Dirichlet heat kernels for discontinuous L\'evy processes is a fundamental problem both in probability theory and analysis. 

The first result on this topic was done in  \cite{CKS10}. In \cite{CKS10},  
the third named author, jointly with Chen and Song, established  
sharp two-sided small time estimates on the Dirichlet heat kernel of  an isotropic  $\alpha$-stable process ($0<\alpha<2$) killed upon leaving a $C^{1,1}$ open set $D$ in $\R^d$. They also obtained large time estimates on the Dirichlet heat kernel when $D$ is bounded.
 After \cite{CKS10}, much has been developed on the  Dirichlet heat kernel estimates for discontinuous Markov processes. In \cite{CT11}, the authors  obtained global Dirichlet heat kernel estimates for an isotropic $\alpha$-stable process $(0<\alpha<2)$ in a half-space like or an exterior $C^{1,1}$ open set $D$ in $\R^d$. Then in \cite{BGR10}, the authors 
 succeeded in extending that result for  a $\kappa$-fat open subset $D$ in $\R^d$, and suggested a  factorization formula for the Dirichlet heat kernel.   Very recently, in \cite{CKSV20}, the first and third named authors, jointly with Song and Vondra\v{c}ek, obtained a general factorization formula for  the Dirichlet heat kernel  in a metric measure space. We refer to \cite{CK16, CKS11, CKS16} for the Dirichlet heat kernel estimates for isotropic L\'evy processes with non-vanishing Gaussian component,  \cite{CKS12a, CKS12b} for relativistic stable processes (see also \cite{GKK20,KK14}), and  mixed $\alpha$-stable processes  in \cite{BGR14b,CKS14}. Note that, every aforementioned result on the Dirichlet heat kernel estimates cover neither  processes with high intensity of small jumps nor  processes with low intensity of small jumps, namely, whose L\'evy measure has the form $\nu(dx)=|x|^{-d-2}\ell(|x|^{-1})dx$ or $\nu(dx)=|x|^{-d}\ell(|x|^{-1})dx$ for a function $\ell$ slowly varying at infinity (in Karamata's sense). The case of high intensity of small jumps  was treated for subordinate Brownian motions in \cite{KM17} by the third named author jointly with Mimica using the result in \cite{M16}. After that, the third named author and Bae extended  that result to a subordinate Brownian motion with non-vanishing Gaussian component in \cite{BK17}.

The purpose of this paper is treating the other case, low intensity of small jumps,  in the study of the Dirichlet heat kernel estimates. We consider the Dirichlet heat kernel estimates  on  isotropic unimodal L\'evy processes  without Gaussian components whose L\'evy measure has a radially non-increasing density which is comparable to $|x|^{-d}\ell(|x|^{-1})$ for a function $\ell$  satisfying weak scaling conditions at infinity with possibly non-positive  lower scaling index (see Definition \ref{d:sv} for the notion of the weak scaling condition). 
Typical examples of such processes are geometric stable processes and iterated geometric stable processes. (See, e.g., \cite[Page 112]{BBKRSV} for the definitions of these processes.) We refer to \cite{GR, KM12, KM14} for the scale invariant version of Harnack inequality and the Green function estimates for these processes.
To the authors best knowledge, our result (Theorems \ref{t:main} and \ref{t:main2}) is the first result on the Dirichlet heat kernel estimates for L\'evy  processes with low intensity of small jumps both in small time and large time. 
Our paper is motivated by the recent paper  \cite{GRT19} where sharp two-sided estimates on the heat kernel in the whole space $\R^d$ for pure jump isotropic unimodal L\'evy processes (without killing)  having the L\'evy measure in the form $\nu(dx)=|x|^{-d}\ell(|x|^{-1})dx$ for a bounded function $\ell$ slowly varying at infinity. Unlike \cite{GRT19}, in this paper,  we allow the function $\ell$ to be unbounded and  not slowly varying at infinity. Hence, our result is even new for the whole space case.

 In this paper, we first derive heat kernel estimates for small time and the whole space by using the results  and methods from \cite{GRT19}. Our heat kernel estimates in $\R^d$  have two forms depending on whether $\ell$ is bounded or unbounded. If the lower scaling index of $\ell$ is positive, then our results can be written in the form of $c_1(p(t, 0)\wedge t\nu(x)) \le p(t,x)\le c_2(p(t, 0)\wedge t\nu(x))$, which coincides with the main result in \cite{BGR14a}.  Hereinafter, $p(t,x)$ denotes the transition density function which is also called as the heat kernel.

Next, we study  behaviours of the process near the boundary of $D$. To do this, we use the boundary Harnack principle and gradient estimates on harmonic functions for pure jump isotropic L\'evy processes. These  results were obtained in \cite{GK18} and \cite{KuR16}, respectively, under some mild assumptions. Under a set of  conditions that give the boundary Harnack principle and gradient estimates (see the condition {\bf (B)} below), we obtain two-sided estimates on the mean exit time from the intersection of an open ball and $D$,  and survival probability in $D$.

Using heat kernel estimates in the whole space and boundary behaviours of the process, we establish small time  two-sided Dirichlet heat kernel estimates for isotropic unimodal L\'evy processes in $C^{1, 1}$ open sets  (see Theorem \ref{t:main} below). 
Even with
 the heat kernel estimates and the precise boundary behaviour of the mean exit time and the survival probability on hand, it is highly non-trivial to obtain Dirichlet heat kernel estimates  because  the lower scaling index of $\ell$ can be $0$ and the heat kernel can be unbounded.  
 
  For bounded $C^{1,1}$ open sets in $\R^d$, we  also  obtain large time estimates (see Theorem \ref{t:main2} below). Since the killed semigroup $\{P_t^D, t>0\}$ may not be compact operators for all $t>0$ even for a bounded open set $D\subset \R^d$, our method is different from ones for obtaining large time Dirichlet heat kernel estimates of stable processes or mixed stable processes.

  Then we obtain two-sided estimates on the Green function in a bounded $C^{1,1}$ open subset in $\R^d$ (see Theorem \ref{t:green} below). This result partially extends the result in \cite{KM14} where the Green function estimtates on subordinate Brownian motions,   whose L\'evy-Khintchine exponent possibly has  the lower scaling index $0$, are treated (see Remark \ref{r:green_re} below).

  The paper ends with an explicit example on the Dirichlet heat kernel estimates and the Green function estimates for some L\'evy processes including geometric stable processes and iterated geometric stable processes. 

\vspace{2mm}

\noindent{\bf Notations}:  We will use the symbol ``$:=$'' to denote a definition, which is read as ``is defined to be.'' For $a,b\in \R$, we denote $a\wedge b:=\min\{a,b\}$ and $a\vee b:=\max\{a,b\}$. For two  functions $f, g$ and a constant $c>0$, the notation $f\asymp g$ for $x>c$ means that there are strictly positive constants $c_1$ and $c_2$ such that $c_1g(x)\leq f (x)\leq c_2 g(x)$ for all $x>c$.  We denote an open ball by $B(x,r):=\{y\in\bR^d: |x-y|<r\}$ and the diagonal set by $\diag:=\{(x,x) : x \in \R^d \}$. For an open set $D$ in $\R^d$, we denote $\delta_D(x):=\dist(x,\partial D)$.

Upper case letters with subscripts $C_i, i=0, 1, 2, \dots,$ and the constants $\kappa_1, \kappa_2, \alpha_1$ and $\alpha_2$ will remain the same throughout this paper.
Lower case letters $c$'s without subscripts denote strictly positive constants  whose values are unimportant and which  may change even within a line, while values of lower case letters with subscripts
$c_i, i=0,1,2,  \dots$, are fixed in each proof, and the labeling of these constants starts anew in each proof.
The notation $c=c(a,b,\ldots)$ denotes a constant depending on $a, b,  \ldots$.

\subsection{Setup}
To describe our results, we introduce the notions of the weak scaling conditions, almost monotonicity and some geometric properties of subsets of $\R^d$.

\begin{defn}\label{d:sv}
{\rm 
Let $f:(0,\infty) \to (0,\infty)$ be a given (Lebesgue) measurable function.\\
(i) For $\alpha_1 \in \R$ and $c_1 > 0$, we say that $f$ satisfies $\WLS^{\infty}(\alpha_1, c_1)$ (resp. $\WLS^0(\alpha_1, c_1)$) if there exists $c>0$ such that
\begin{equation}\label{e:lsc}
\qquad \frac{f(R)}{f(r)}\ge c \left(\frac{R}{r}\right)^{\alpha_1}  \quad \text{for all} \;\; c_1<r \le R \;\; (\text{resp.} \;\; 0<r \le R \le c_1).
\end{equation}
Similarly, for $\alpha_2 \in \R$ and $c_2 > 0$, we say that $f$ satisfies $\WUS^{\infty}(\alpha_2, c_2)$ (resp. $\WUS^0(\alpha_2, c_2)$) if there exists $c>0$ such that
\begin{equation}\label{e:usc}
\qquad \frac{f(R)}{f(r)}\le c \left(\frac{R}{r}\right)^{\alpha_2}  \quad \text{for all} \;\; c_2<r \le R \;\; (\text{resp.} \;\; 0<r \le R \le c_2).
\end{equation}
If $f$ satisfies \eqref{e:lsc} (resp. \eqref{e:usc}), then we call $\alpha_1$ (resp. $\alpha_2$) the lower scaling index (resp. the upper scaling index) of the function $f$. 
If $f$ satisfies both $\WLS^{\infty}(\alpha_1, c_3)$ and $\WUS^{\infty}(\alpha_2, c_3)$ for some $\alpha_1, \alpha_2\in\R$ and $c_3 > 0$, we say that $f$ satisfies $\WS^{\infty}(\alpha_1,\alpha_2,c_3)$.\\
(ii) We say that $f$ is almost increasing if there exists $c_0>0$ such that 
$$f(x) \asymp \sup_{y \in [c_0, x]}f(y) \quad \text{for all} \;\; x>c_0.$$
Similarly, we say that $f$ is almost decreasing if there exists $c_0>0$ such that 
$$f(x) \asymp \inf_{y \in [c_0, x]}f(y) \quad \text{for all} \;\; x>c_0.$$
}
\end{defn}

\begin{definition}\label{c11}
{\rm (i) Let $d \ge 2$. An open set $D$ 
in $\bR^d$ is said to be a (uniform) $C^{1,1}$ open set if there exist a localization radius $R_0>0$ and a constant $\Lambda>0$ such that for 
every $Q\in\partial D$, there exist a $C^{1,1}$ function $\Gamma=\Gamma_{Q}: \bR^{d-1}\to \bR$
satisfying $\Gamma(0)=0$, $ \nabla\Gamma (0)=(0, \dots, 0)$, $\| \nabla
\Gamma \|_\infty \leq \Lambda$, $| \nabla \Gamma (\wt x)-\nabla \Gamma (\wt w)|
\leq \Lambda |\wt x- \wt w|$ for $\wt x, \wt w\in \bR^{d-1}$ and an orthonormal coordinate system $CS_{Q}:y=( \wt y, \, y_d)$ with its origin at $Q$
such that 
$$B(Q, R_0 )\cap D= \big\{y=(\wt y, y_d) \in B(0, R_0) \mbox{
in } CS_{Q}: y_d > \Gamma ( \wt y) \big\}.$$ 
The pair $( R_0, \Lambda)$ is called the characteristics of the $C^{1,1}$ open set $D$. \\
(ii) An open set $D$ in $\R$ is said to be a $C^{1,1}$ open set if there exists a localization radius $R_0>0$ such that $D$ is an union of open intervals of length at least $R_0$ and distanced one from another at least $R_0$. \\
(iii) A bounded set $D$ in $\R^d$ is said to be of scale $(r_1, r_2)$ if there exist $x_1, x_2 \in \R^d$ such that $B(x_1, r_1) \subset D \subset B(x_2, r_2)$.

}
\end{definition}

Let $Y=(Y_t, t\ge0)$ be a  L\'evy process in $\bR^d$ with the L\'evy-Khintchine exponent $\psi$. Then,
\begin{align*}
\bE\Big[\exp\big(i\ip{\xi,Y_t}\big)\Big]=\int_{\R^d}e^{i\ip{\xi, x}}p(t,dx)=e^{-t\psi(\xi)}, \quad \xi \in \R^d,
\end{align*}
where $p(t, dx)$ is the transition probability of $Y$. If $Y$ is a pure jump symmetric L\'evy process with L\'evy measure $\nu$, then $\psi$ is of the form
\begin{align*}
\psi(\xi)=\int_{\R^d}(1-\cos\ip{\xi, x})\nu(dx), \quad \xi \in \R^d,
\end{align*}
where $\int_{\R^d} (1\wedge |x|^2)\nu(dx)<\infty$. 

A measure $\mu(dx)$ is isotropic unimodal if it is absolutely continuous on $\R^d\setminus\{0\}$ with a radial and radially non-increasing density. A L\'evy process $Y$ is isotropic unimodal if $ p(t, dx)$ is isotropic unimodal for all $t > 0$. This is equivalent to the condition that the L\'evy measure $\nu(dx)$ of $Y$ is isotropic unimodal if $Y$ is pure jump L\'evy process. (See, \cite{W83}.) 

 Throughout this paper, we always assume that $Y$ is a pure jump isotropic unimodal L\'evy process with the L\'evy-Khintchine exponent $\psi$. 
With a slight abuse of notation, we will use the notations  $\psi(|x|) = \psi(x)$  and $\nu(dx)=\nu(x)dx=\nu(|x|)dx$ for $x \in \R^d$. Then, throughout this paper, we also assume the following condition {\bf (A)} holds.  A smooth function $\vp:[0,\infty) \to [0,\infty)$ is called a \textit{Bernstein function} if  $(-1)^{n-1}\vp^{(n)}(\lambda)\ge 0$ for all $n\ge 1$ and $\lambda>0$.

\vspace{2mm}

\setlength{\leftskip}{4mm}

\noindent {\bf (A)} The L\'evy measure $\nu$ on $\R^d$ is infinite and there exist constants $\kappa_1, \kappa_2>0$ and a continuous function $\ell:(0, \infty) \to (0, \infty)$ satisfying $\WS^{\infty}(\alpha_1, \alpha_2, 1)$ for some 
 $-d < \alpha_1 \le \alpha_2 <2 $ such that
\begin{equation}\label{A}
\kappa_1 r^{-d} \ell(r^{-1}) \le \nu(r) \le \kappa_2 r^{-d} \ell(r^{-1})\quad\text{for all}\;\;r>0.
\end{equation}
 If $d>1$, then we assume further that either  $\alpha_1>-1$  or $\psi(\xi)= \vp(|\xi|^2)$ for a Bernstein function $\vp$.

\setlength{\leftskip}{0mm}

\smallskip

Note that, since we allow the constant $\alpha_1$ to be negative, the map $r \mapsto \ell(r^{-1})$ can be  increasing near zero.

Here, we enumerate other main conditions which we will assume later.

\setlength{\leftskip}{4mm}

\vspace{2mm}

\noindent {\bf (B)}
$\nu(r)$ is absolutely continuous such that $r \mapsto -\nu'(r)/r$ is non-increasing on $(0, \infty)$ and there exists a constant $c_0>1$ such that $\nu(r) \le c_0 \nu(r+1)$ for all $r \ge 1$;

\smallskip

\noindent{\bf (C)} $\ell(r)$ satisfies $\WUS^{0}(\gamma, 1)$  for some $\gamma<2$;

\smallskip

\noindent{\bf (S-1)}
$\limsup_{r \to \infty} \ell(r) < \infty$;

\smallskip

\noindent{\bf (S-2)}
$\limsup_{r \to \infty} \ell(r) = \infty$ and $\ell(r)$ is almost increasing;

\smallskip

\noindent {\bf (L-1)}
$\liminf_{r \to \infty} \ell(r) = 0$ and $\ell(r)$ is almost decreasing;

\smallskip

\noindent {\bf (L-2)}
$0<\liminf_{r \to \infty} \ell(r) \le \limsup_{r \to \infty} \ell(r) < \infty$;

\smallskip

\noindent{\bf (D)}  If $d=1$, then 
$\alpha_2 < 1$  where $\alpha_2$ is the constant in {\bf (A)}.

\setlength{\leftskip}{0mm}

\vspace{2mm}

\begin{remark}\label{r:remark0}
{\rm

Let $B=(B_t,\, t\ge 0)$ be a Brownian motion in
$\R^d$  and $S=(S_t,\, t\ge 0)$ be a driftless subordinator
 independent of $B$.
The process $X=(X_t:\, t\ge 0)$ defined by  $X_t=B_{S_t}$ is called  a  \textit{subordinate Brownian motion} (SBM). Every SBM is an
isotropic unimodal L\'evy process. 
Let $\vp$ be \textit{the Laplace  exponent} of the subordinator $S$, namely,
$$\E\big[\exp(-\lambda S_t)\big]=\exp\big(-t\vp(\lambda)\big), \qquad \lambda \ge 0.$$
It is known that the Laplace exponent $\vp$ is a Bernstein function with $\vp(0)=0$. Since 
$S$ has no drift, $\vp$ has the representation $\vp(\lambda)= \int_0^\infty (1-e^{\lambda t})\mu(dt)$ where $\mu$ is a measure on $(0,\infty)$ satisfying $\int_0^\infty (1 \land t)\mu(dt)<\infty$, called the L\'evy measure of $\vp$.
Note that  the characteristic exponent of $X$ is   $\vp(|\xi|^2)$.

A function $f:(0,\infty)\to [0,\infty)$ is said to be \textit{completely monotone},  if  $(-1)^nf^{(n)}\ge 0$ on $(0,\infty)$ for every  $n\geq 0$.
A Bernstein function is said to be a \textit{complete Bernstein function}, if its
L\'evy measure has a completely monotone density.
\smallskip

\noindent (i) Suppose that $\vp$ is a complete Bernstein function such that $\lim_{\lambda \to \infty}\vp(\lambda)=\infty$ and $\vp'$ satisfies $\WUS^\infty(-\delta, 1)$ for some $\delta \in (\frac{1}{2},1]$. Suppose further that, if $d=2$, then $\vp'$ satisfies $\WLS^\infty(-\delta_0,1)$ for some $\delta_0 \in (0,2)$, and if $d=1$, then  $\vp'$ satisfies $\WLS^\infty(-\delta_0,1)$ for some $\delta_0 \in (\frac{1}{2},2\delta-\frac{1}{2})$. Then according to \cite[Proposition 2.6]{KM14}, a SBM with the characteristic exponent $\vp(|\xi|^2)$ satisfies {\bf (A)} with $\ell$ such that $\ell(r) \asymp r^2\vp'(r^2)$ for $r \ge 1$.

\noindent (ii) Let $X$ be a SBM with the characteristic exponent $\vp(|\xi|^2)$ for a complete Bernstein function $\vp$. Then by  \cite[Remark 1.4]{GKK20} and \cite[Lemma 7.4]{BGPR}, it satisfies {\bf (B)}. (See also  the proof of \cite[Proposition 3.5(b)]{KSV12}.)
}
\end{remark}

\begin{remark}\label{r:remark1}
{\rm

\noindent (i)  {\bf (A)} implies that $\nu(r)$ satisfies $\WLS^{0}(-d-\alpha_2, 1)$. Therefore, under {\bf (A)}, for every $R>0$, there exists $c>0$ such that
\begin{align}\label{e:doubling}
\nu(2r) \ge c \nu(r) \quad \text{for all} \;\; r \in (0, R].
\end{align}
On the other hand, {\bf (C)} implies that $\nu(r)$ satisfies $\WLS^{\infty}(-d-\gamma,1)$ for some $\gamma<2$. Thus, {\bf (C)} implies that for every $R>0$, there exists $c>0$ such that
\begin{align}\label{e:doubling2}
\nu(2r) \ge c \nu(r) \quad \text{for all} \;\; r \in [R, \infty).
\end{align}

\smallskip

\noindent (ii) If {\bf (A)} holds with $\alpha_1>0$, then {\bf (S-2)} holds. (See, \cite[Section 1.5]{BGT}.) 

}
\end{remark}

\subsection{Main results.}

We define for $r>0$,
\begin{align*}
 K(r)&:=r^{-2}\int_0^r s\ell(s^{-1})ds, \qquad L(r):=\int_r^{\infty} s^{-1} \ell(s^{-1})ds, \\
h(r)&:=K(r)+L(r).
\end{align*}
Since  {\bf (A)} holds, we see that 
\begin{align*}
K(r) \asymp r^{-2}\int_{|y|\le r}|y|^2\,\nu(y)dy \quad \text{and} \quad L(r) \asymp \int_{|y|>r} \nu(y) dy,
\end{align*}
which are the functions introduced in \cite{Pr81}.  We also define
\begin{align*}
\ell^*(r):= \sup_{u \in [1,r]} \ell(u) \quad \text{for} \;\; r \ge 1
\end{align*}
and denote by $\ell^{-1}$ the right continuous inverse of $\ell^*$, that is,
\begin{align}\label{ell_inv}
\ell^{-1}(t):= \inf\{r \ge 1: \ell^*(r)>t\} \quad \text{for} \;\; t>0.
\end{align}

 Now, we are ready to state our main results. Recall that $\delta_D(x)=\dist(x,\partial D)$.

\begin{thm}\label{t:main}
Suppose that $Y$ is a pure jump isotropic unimodal L\'evy process satisfying  {\bf(A)} and {\bf (B)}. Let $D$ be a $C^{1,1}$ open set in $\R^d$ with characteristics $(R_0,\Lambda)$. If $D$ is unbounded, we further assume that {\bf (C)} holds. 
Then, the following estimates hold:

\smallskip

\noindent (i) If {\bf (S-1)} holds, then for every $T>0$, there exist positive constants $c_1=c_1(d, \psi, T, R_0, \Lambda)$, $c_2=c_2(d,\psi, T)$ and $c_3=c_3(d, \psi, T, R_0, \Lambda)>1$ such that
\begin{align*}
& c_3^{-1} \left(1\wedge \frac{1}{tL(\delta_D(x))} \right)^{1/2}\left(1\wedge \frac{1}{tL(\delta_D(y))} \right)^{1/2} t \nu(|x-y|)\exp\big(-c_1th(|x-y|)\big) \\[3pt]
&\le p_D(t,x,y) \le c_3 \left(1\wedge \frac{1}{tL(\delta_D(x))} \right)^{1/2}\left(1\wedge \frac{1}{tL(\delta_D(y))} \right)^{1/2} t \nu(|x-y|)\exp\big(-c_2th(|x-y|)\big),
\end{align*}
for all $(t,x,y) \in (0,T] \times (D \times D \setminus \diag)$. 

\smallskip

\noindent (ii) If {\bf (S-2)} holds, then for every $T>0$ and $\eta > 0$,  there exist positive constants $a_0=a_0(d,\psi)$, $c_4=c_4(d, \psi, T, R_0, \Lambda),$ $c_5=c_5(d,\psi)$ and $c_6=c_6(d, \psi, T, \eta, R_0,\Lambda)>1$  such that
\begin{align*}
& c_6^{-1} \left(1\wedge \frac{1}{tL(\delta_D(x))} \right)^{1/2}\left(1\wedge \frac{1}{tL(\delta_D(y))} \right)^{1/2} t \nu(\theta_{\eta}(|x-y|,t))\exp\big(-c_4th(\theta_{\eta}(|x-y|,t))\big)\\[3pt]
&\le p_D(t,x,y) \\
&\le c_6 \left(1\wedge \frac{1}{tL(\delta_D(x))} \right)^{1/2}\left(1\wedge \frac{1}{tL(\delta_D(y))} \right)^{1/2} t \nu(\theta_{a_0}(|x-y|,t))\exp\big(-c_5th(\theta_{a_0}(|x-y|,t))\big),
\end{align*}
for all $(t,x,y) \in (0,T] \times D \times D$ where $\theta_a(r,t):=r \vee [\ell^{-1}(a/t)]^{-1}$ and $\ell^{-1}$ is defined as \eqref{ell_inv}.
\end{thm}

\smallskip

If we further assume that $D$ is bounded, then we can obtain the large time estimates for the Dirichlet heat kernel  and the Green function estimates under some mild assumptions.

\smallskip

\begin{thm}\label{t:main2}
Suppose that $Y$ is a pure jump isotropic unimodal L\'evy process satisfying  {\bf(A)} and {\bf (B)}. Let $D$ be a bounded $C^{1,1}$ open set in $\R^d$ with characteristics $(R_0,\Lambda)$ of scale $(r_1, r_2)$. Then, the following estimates hold:

\vspace{2mm}

\noindent(i) If  {\bf (L-1)} holds,  then for every $T>0$, there exist positive constants
$c_1=c_1(d, \psi),$ $c_2=c_2(d,\psi)$ and $c_3=c_3(d, \psi, T, R_0, \Lambda, r_1, r_2)>1$ such that
\begin{align*}
& c_3^{-1} L(\delta_D(x))^{-1/2}L(\delta_D(y))^{-1/2} \Big(\nu(|x-y|)\exp\big(-c_1th(|x-y|)\big) + \exp\big(-\kappa_2 C_4 th(r_1/2)\big)\Big) \\
&  \le p_D(t,x,y) \\
&  \le c_3  L(\delta_D(x))^{-1/2}L(\delta_D(y))^{-1/2} \Big(\nu(|x-y|)\exp\big(-c_2th(|x-y|)\big)+\exp\big(-\frac{\kappa_1 C_5}{2}th(r_2)\big)\Big),
\end{align*}
for all $(t,x,y) \in [T, \infty) \times (D \times D \setminus \diag)$ where $\kappa_1$ and $\kappa_2$ are the positive constants in {\bf (A)} and $C_4$ and $C_5$ are positive constants which only depend on the dimension $d$. 

\vspace{2mm}

\noindent(ii) If {\bf (L-2)} holds, then there exist $T_1 \ge 0$ and $\lambda_1 = \lambda_1(\psi, D)>0$ such that for every fixed $T>T_1$, there exists $c_4=c_4(d, \psi, T, R_0, \Lambda, r_1, r_2)>1$ such that
\begin{align*}
c_4^{-1}e^{- \lambda_1 t}  L(\delta_D(x))^{-1/2}L(\delta_D(y))^{-1/2} \le p_D(t,x,y) \le c_4e^{- \lambda_1 t}  L(\delta_D(x))^{-1/2}L(\delta_D(y))^{-1/2},
\end{align*}
for all $(t,x,y) \in [T,\infty) \times D \times D$. Moreover, we have
\begin{align*}
\frac{\kappa_1C_5}{2}h(r_2) \le \lambda_1 \le \kappa_2C_4 h(r_1/2).
\end{align*}

\noindent(iii) If {\bf (S-2)} holds, then the estimates in (ii) holds with $T_1 = 0$. Moreover, the constant $-\lambda_1<0$ is the largest eigenvalue of the generator of $Y^D$.
\end{thm}

\vspace{3mm}

For a Borel subset $D \subset \R^d$, the Green function $G_D(x,y)$ of $Y$ in $D$ is defined by
\begin{align*}
G_D(x,y):= \int_0^{\infty} p_D(t,x,y) dt \qquad \text{for} \;\; x, y \in D.
\end{align*}

\begin{thm}\label{t:green}
Suppose that $Y$ is a pure jump isotropic unimodal L\'evy process satisfying  {\bf(A)}, {\bf(B)} and {\bf (D)}. Let $D$ be a bounded $C^{1,1}$ open subset in $\R^d$ with characteristics $(R_0,\Lambda)$ of scale $(r_1, r_2)$. Then, the Green function $G_D(x,y)$ of $Y$ in $D$ satisfies the following two-sided estimates: for every $x,y \in D$,
\begin{align}\label{e:Gest}
G_D(x, y) & \asymp  \bigg(1 \wedge \frac{L(|x-y|)}{\sqrt{L(\delta_D(x)) L(\delta_D(y))}} \bigg) \frac{ \ell(|x-y|^{-1}) }{ |x-y|^d L(|x-y|)^2},
\end{align}
where the comparison constants depend only on $d, \psi, R_0, \Lambda$ and $r_2$.
\end{thm}

\begin{remark}\label{r:green_re}
{\rm 
(i) One can obtain \eqref{e:Gest} just by integrating the estimates for $p_D(t,x,y)$ given 
in Theorems \ref{t:main} and \ref{t:main2}
(e.g. \cite[Theorem 7.3]{KM17}). However, to use Theorems \ref{t:main} and \ref{t:main2},  we need conditions more than  {\bf (A)}, {\bf (B)} and {\bf (D)}.  By adopting arguments from \cite{KM14} instead of integrating the Dirichlet heat kernel, we obtained the Green function estimates in more general situations.

\noindent (ii) It is established in \cite[Theorem 1.2]{KM14}  that   for a large class of transient subordinate Brownian motions, the Green function $G_D(x,y)$ in a bounded $C^{1,1}$ open set $D$ enjoys the following sharp two-sided estimates:
\begin{align}\label{e:Gestold}
G_D(x,y) \asymp \bigg(1 \wedge \frac{\vp(|x-y|^{-2})}{\sqrt{\vp(\delta_D(x)^{-2}) \vp(\delta_D(y)^{-2})}}\bigg) \frac{\vp'(|x-y|^{-2})}{|x-y|^{d+2}\vp(|x-y|^{-2})^2}.
\end{align} 
An important novelty of this result is that it was the first explicit Green function estimates even if the lower scaling index of the  L\'evy-Khintchine exponent can be $0$. Note that, in view of Remark \ref{r:remark0} and \cite[Lemma 4.1]{KM12}, when the lower scaling index of the  L\'evy-Khintchine exponent can be $0$,  assumptions {\bf (A-1)}--{\bf (A-5)} in  \cite[Theorem 1.2]{KM14} imply the following:

\vspace{1mm}

\setlength{\leftskip}{4mm}

\noindent(1) The L\'evy-Khintchine exponent $\psi(\xi) = \vp(|\xi|^2)$ for a complete Bernstein function $\vp$.  Thus, {\bf (B)} holds (see Remark \ref{r:remark0}(ii));

\noindent (2) {\bf (A)} holds with $\ell$ such that $\ell(r) \asymp r^2 \vp'(r^2)$ for $r\ge 1$, and constants $\alpha_1 \in  (-d,1)  \cap  [-2, 1)$ and $\alpha_2<1$. Thus,  {\bf (D)} holds.

\vspace{1mm}

\setlength{\leftskip}{0mm}

\noindent Therefore, by Lemma \ref{l:asymLh} and \eqref{e:asym}, we see that  Theorem \ref{t:green} recovers \eqref{e:Gestold}. Here, we note that Theorem \ref{t:green} does not assume the transience of the process unlike \cite[Theorem 1.2]{KM14}.

}
\end{remark}

\section{Heat kernel estimates in $\R^d$}\label{s:hke}

Recall that under {\bf (A)}, we have
\begin{align*}
& K(r)=r^{-2}\int_0^r s\ell(s^{-1})ds \asymp r^{-2}\int_{|y|\le r}|y|^2\,\nu(y)dy, \\
& L(r)=\int_r^{\infty} s^{-1} \ell(s^{-1})ds \asymp \int_{|y|>r} \nu(y) dy, \\
& h(r)=K(r)+L(r) \asymp r^{-2}\int_{\R^d} \left(r^2 \wedge |y|^2 \right) \, \nu(y)dy.
\end{align*}
Clearly, $L(r)$ is decreasing. Moreover, we see that $h'(r)=-2r^{-1}K(r) \le 0$ for all $r>0$ and hence $h(r)$ is also decreasing.
Since the underlying process $Y$ is isotropic unimodal, there are a number of general properties related to these functions. (See, \cite{BGR14a}, \cite{BGR15} and \cite{GK18}.) 

First, since $\nu(r)$ is non-increasing, we have 
\begin{align}\label{Kgenu}
K(r) \ge c \nu(r)r^{-2} \int_0^r s^{d+1}ds= c r^d\nu(r)\quad\text{for all}\;\;r>0.
\end{align}
On the other hand, by Karamata's Tauberian-type theorem, the opposite inequality $K(r)\le c r^d\nu(r)$ holds for $0<r\le 1$ if and only if $\ell(r)$ satisfies $\WUS^{\infty}(\gamma, 1)$ for some $\gamma<2$. Similarly, we have $K(r) \le cr^d \nu(r)$ for $r \ge 1$ if and only if $\ell(r)$ satisfies   $\WUS^0(\gamma', 1)$ for some $\gamma'<2$. (See, \cite[Appendix A]{GK18}.) In particular, {\bf (A)} implies that 
\begin{align}\label{smallKnu}
K(r) \asymp r^d\nu(r) \asymp \ell(r^{-1}) \quad \text{for} \;\; 0<r \le 1
\end{align}
and {\bf (C)} implies that
\begin{align}\label{largeKnu}
K(r) \asymp r^d\nu(r) \asymp \ell(r^{-1}) \quad \text{for} \;\; r \ge 1.
\end{align}

Next, by \cite[(6) and (7)]{BGR14a}, there exist positive constants $C_0$ and $C_1$ which only depend on the dimension $d$ and $\kappa_1$ and $\kappa_2$ in \eqref{A}  such that for all $r>0$,
\begin{align}\label{e:asym}
C_0h(r) \le \psi(r^{-1}) \le C_1 h(r).
\end{align}
Under {\bf (A)}, we can extend this relations by including  $L(r)$ if $r$ is small.

\begin{lemma}\label{l:asymLh}
There exists a constant $c_1>0$ such that
\begin{align*}
L(r) \le h(r) \le c_1L(r) \quad \text{for all} \;\; 0<r \le 1.
\end{align*}
\end{lemma}
\pf
From the definitions of $L$ and $h$, the first inequality is obvious. To prove the second inequality, it suffices to show that there exists $c>0$ such that $L(r) \ge cK(r)$ for $0< r \le 1$. Since {\bf (A)} holds,  by \eqref{e:doubling} and \eqref{smallKnu}, we have $\nu(r) \asymp \nu(2r)$ and $K(r) \asymp r^d \nu(r)$ for $0<r\le1$. Thus, for $0<r \le 1$, we get
\begin{align*}
L(r) \ge c \int_r^{2r} s^{d-1} \nu(s)ds \ge c r^d \nu(r) \ge c K(r). 
\end{align*}
\qed

By Lemma \ref{l:asymLh} and \eqref{e:asym}, we deduce that $L(r) \asymp \psi(r^{-1})$ for small $r$. In view of this relation, to make some computations easier, we define $\Phi : [0, \infty) \to [0, \infty)$ by 
\begin{align*}
\Phi(r):= L(r^{-1}) = \int_{r^{-1}}^{\infty} u^{-1}\ell(u^{-1})du = \int_0^r s^{-1} \ell(s)ds.
\end{align*} 
We used the change of variables $u=s^{-1}$ in the last equality.

\begin{lemma}\label{l:propPhi}
(i) $\Phi(r)$ satisfies $\WS^{\infty}(\alpha_1, \alpha_2 \vee \frac{1}{2}, 1)$. \\
(ii) We have that
\begin{align}\label{e:eq1}
C_0\Phi(r) \le \psi(r) \quad \text{for all} \;\; r \ge 0.
\end{align}
Moreover, there exists a constant $C_2>0$ such that
\begin{align}\label{e:eq2}
C_2\Phi(r) \ge h(r^{-1}) \quad \text{for all} \;\; r\ge1.
\end{align}
\end{lemma}
\pf
(i) Let $\alpha_2'=\alpha_2 \vee \frac{1}{2}$.  By the change of variables and {\bf (A)}, we have that
\begin{equation}\label{e:l2.2}
c_1 \kappa^{\alpha_1} \Phi(r) \le \Phi(\kappa r) =  \int_{1}^r s^{-1} \ell(s) \frac{\ell(\kappa s)}{\ell(s)} ds + \Phi(\kappa) \le c_2\kappa^{\alpha_2'} \Phi(r) + \Phi(\kappa) \;\; \text{ for all } \; \kappa, r>1.
\end{equation}
The first inequality in \eqref{e:l2.2} shows that $\Phi$ satisfies $\WLS^{\infty}(\alpha_1, 1)$. 

Now, we prove that $\Phi$ satisfies $\WUS^{\infty}(\alpha_2', 1)$. Choose any $\kappa>1$ and  $r \ge 2$.  Let $n$ be the smallest integer satisfying $r^{n} \ge \kappa$. By applying the latter inequality in \eqref{e:l2.2} $n$ times, since $\Phi(r)$ is increasing, we obtain
\begin{align*}
&\Phi(\kappa r) \le c_2\kappa^{\alpha_2'} \Phi(r) + \Phi( \frac{\kappa}{r} r) \le c_2\kappa^{\alpha_2'}(1+ r^{-\alpha_2'}) \Phi(r) + \Phi(\frac{\kappa}{r^2}r) \\
&\le \cdots \le c_2\kappa^{\alpha_2'}(1+r^{-\alpha_2'}+\cdots +r^{-(n-1)\alpha_2'})\Phi(r) + \Phi(\frac{\kappa}{r^{n}}r) \le \big((1-2^{-1/2})^{-1}c_2+1\big) \kappa^{\alpha_2'}\Phi(r).
\end{align*}
Besides, for any $\kappa>1$ and $1<r<2$, since $\Phi$ is increasing, we see from the above inequalities that $\Phi(\kappa r)/\Phi(r) \le (\Phi(2)/\Phi(1)) \cdot ( \Phi(2\kappa)/\Phi(2)) \le c_3\kappa^{\alpha_2'}$. Hence, we get the desired result.
 
\noindent (ii) It follows from the definition of $\Phi$, Lemma \ref{l:asymLh} and \eqref{e:asym}.
\qed

Let $C_{\infty}(\R^d)$ be the set of all continuous functions which vanish at infinity. In \cite{HW42}, Hartman and Wintner proved sufficient conditions in terms of the L\'evy exponent $\psi$ under which the transition density $p(t, \cdot)$ of $Y$ is in $C_{\infty}(\R^d)$. Then, in \cite{KS13}, Knopova and Schilling improve that result and they also give some necessary conditions. Using \eqref{e:eq1} and \eqref{e:eq2}, we can formulate these conditions in terms of $\Phi$. Since the underlying process $Y$ is isotropic unimodal, these conditions determine whether $p(t,0) < \infty$ or $p(t,0) = \infty$.

\begin{prop}\label{p:existence}
Let $p(t, \cdot)$ be the transition density of $Y$. Suppose that
\begin{align*}
\liminf_{r \to \infty} \frac{\Phi(r)}{\log(1+r)} = c_1\in [0, \infty], \quad \limsup_{r \to \infty} \frac{\Phi(r)}{\log(1+r)} = c_2\in [0, \infty].
\end{align*}
Then, the following are true.\\
(i) If $c_1=\infty$, then $p(t, 0)< \infty$ for all $t>0$. \\
(ii) If $c_2 = 0$, then $p(t, 0) = \infty$ for all $t>0$. \\
(iii) If $0<c_1 \le c_2< \infty$, then there exist $T_2 \ge T_1>0$ such that $p(t, 0) = \infty$ for $0<t \le T_1$ and $p(t, 0) < \infty$ for $t>T_2$.

\vspace{2mm}

In particular, by l'Hospital's rule, the following are true.\\
(iv) If $\;\liminf_{r \to \infty} \ell(r) = \infty$, then $p(t, 0)<\infty$ for all $t>0$.\\
(v) If $\;\limsup_{r \to \infty} \ell(r) = 0$, then $p(t, 0) = \infty$ for all $t>0$. \\
(vi) If $0<\liminf_{r \to \infty} \ell(r) \le \limsup_{r \to \infty} \ell(r)< \infty$, then there exist $T_2 \ge T_1>0$ such that $p(t, 0) = \infty$ for $0<t \le T_1$ and $p(t, 0) < \infty$ for $t>T_2$.
\end{prop}
\pf
By \eqref{e:eq1} and \eqref{e:eq2}, the first two assertions follow from Part II in \cite{HW42} and the third one follows from \cite[Lemma 2.6]{KS13}.
\qed

Here, we introduce some general estimates which are established in \cite{GRT19}. Note that the following estimates hold no matter $p(t, 0) < \infty$ or $p(t, 0) = \infty$.

\begin{prop}[{\cite[Proposition 5.3]{GRT19}}]\label{p:glb}
There are constants $b_0, c_0 >0$, which only depend on the dimension $d$ and $\kappa_2$ in \eqref{A} such that for all $(t,x) \in (0,\infty) \times \R^d$,
$$
p(t,x) \ge c_0t\nu(|x|)\exp\big(-b_0 t h(|x|)\big).
$$
\end{prop}

\begin{prop}[{\cite[Theorem 5.4]{GRT19}}]\label{p:gub}
There is a constant $c_1>0$, which only depends on the dimension $d$ and  $\kappa_2$ in \eqref{A}  such that for all $t>0$ and $x \in \R^d \setminus \{0\}$,
$$
p(t,x) \le c_1 t|x|^{-d}K(|x|).
$$
\end{prop}

The following lemma will be used several times to obtain heat kernel upper bounds for the whole space. (Cf. \cite[Lemma 4.1 and Corollary 4.4]{GRT19}.)
\begin{lemma}\label{l:simple} For every $\lambda>1$, there exists a constant  $c=c(\lambda)>0$ such that
\begin{equation}\label{e:simple}
\sup_{1 < k \le \lambda}|\psi(k r) - \psi(r)| \le c \ell (r) \quad\;\; \text{for all} \;\; r \ge 1.
\end{equation}
\end{lemma}
\pf Recall the condition {\bf (A)}. We first assume that either $d=1$ or  $\alpha_1>-1$. For $y>0$, set $\nu_1(y)=\nu(y)$ if $d=1$, and 
\begin{equation*}
\nu_1(y):=\int_{\R^{d-1}} \nu\big((y^2+|z|^2)^{1/2}\big) dz \quad \text{if} \;\; d \ge 2.
\end{equation*}
We claim that there exists a constant $c_1>0$ such that
 \begin{equation}\label{e:l2.7}
\nu_1(y) \le  c_1y^{-1}\ell(y^{-1}) \quad \text{for all} \;\; y \in (0,1].
\end{equation} 
If $d=1$, then \eqref{e:l2.7} follows from \eqref{A}. Hence, we assume  $\alpha_1>-1$ and $d \ge 2$. Since $\ell$ is continuous and satisfies $\WS^\infty(\alpha_1, \alpha_2, 1)$, it also satisfies $\WS^\infty(\alpha_1, \alpha_2, 1/2)$. 
Hence, according to \eqref{A} and the change of the variables, we have that, for any $y \in (0,1]$,
\begin{align*}
&\frac{1}{y^{-1}\ell(y^{-1})}\int_0^1 \nu\big((y^2+k^2)^{1/2}\big)k^{d-2}dk \\
&\asymp \int_0^1 \frac{yk^{d-2}}{(y^2+k^2)^{d/2}} \frac{\ell((y^2+k^2)^{-1/2})}{\ell(y^{-1})} dk  = \int_0^{1/y} \frac{k^{d-2}}{(1+k^2)^{d/2}} \frac{\ell(y^{-1}(1+k^2)^{-1/2})}{\ell(y^{-1})}dk  \\
& \le c_2 \int_0^{1/y} \frac{k^{d-2}}{(1+k^2)^{(\alpha_1+d)/2}}dk \le c_2\int_0^1 dk +  c_2\int_1^{\infty} k^{-2-\alpha_1}dk= \frac{c_2(2+\alpha_1)}{1+\alpha_1}.
\end{align*}
Besides, since $\nu$ is non-increasing and a L\'evy measure,  we also have that for any $y \in (0,1]$,
\begin{align*}
&\frac{1}{y^{-1}\ell(y^{-1})}\int_1^\infty \nu\big((y^2+k^2)^{1/2}\big)k^{d-2}dk \le \frac{1}{y^{-1}\ell(y^{-1})}\int_1^\infty \nu(k)k^{d-1}dk \\
&= \frac{c_3}{y^{-1}\ell(y^{-1})}\int_{\xi \in \R^{d}, \, |\xi|>1} \nu(\xi) d\xi= \frac{c_4}{\ell(1)}\frac{\ell(1)}{y^{-1}\ell(y^{-1})} \le c_5 y^{1+\alpha_1}\le c_5.
\end{align*}
Therefore, we obtain \eqref{e:l2.7} with $c_1 = c_2(2+\alpha_1)/(1+\alpha_1)+c_5$.

Observe that for $r>0$,
\begin{align*}
\psi(r) = \int_{\R}\int_{\R^{d-1}} \big(1-\cos (rz_1)\big) \nu\big((z_1^2 + |\wt z|^2)^{1/2}\big) d\wt zdz_1 = 2\int_0^{\infty} \big(1-\cos (ry)\big) \nu_1(y) dy.
\end{align*}
Hence, we see that for any $1<k \le \lambda$ and $r \ge 1$, 
\begin{align*}
&|\psi(k r) - \psi(r)| = 2\,\bigg| \int_0^{\infty}\big(\cos(ry)-\cos(k ry)\big)\nu_1(y)dy\bigg|\\
& =  2r^{-1} \bigg|\int_0^{1} \big(\cos (y) - \cos (k y)\big) \nu_1(y/r) dy+ \int_1^{\infty} \big(\cos (y) - \cos (k y)\big) \nu_1(y/r) dy\bigg| \\
& \le  2r^{-1} \int_0^1 \big|\cos (y) - \cos (k y)\big| \nu_1(y/r) dy \\
& \quad + 2r^{-1} \left|\int_1^{\infty} \cos (y) \nu_1(y/r) dy \right| + 2r^{-1} \left|\int_1^{\infty} \cos (k y) \nu_1(y/r) dy \right|\\
& =: I_1+I_2+I_3.
\end{align*}

By Taylor expansion of the cosine function,  \eqref{e:l2.7} and the assumption that $\ell$ satisfies $\WUS^\infty(\alpha_2,1)$ with $\alpha_2<2$,  we have
\begin{align*}
I_1 \le 2\lambda^2 r^{-1} \int_0^1 y^2  \nu_1(y/r) dy \le 2c_1\lambda^2 \ell(r) \int_0^1 y \frac{\ell(r/y)}{\ell(r)} dy  \le c\lambda^2 \ell(r) \int_0^1 y^{1-\alpha_2} dy  = c\lambda^2 \ell(r).
\end{align*}

Next, to bound $I_2$ and $I_3$, we use a trick from the proof of \cite[Theorem 3.5]{GRT19}. Since $y \mapsto \nu_1(y)$ is non-increasing, there exists a measure $-d\nu_1$ on $(0,\infty)$ such that $\nu_1(y)=\int_y^\infty (-d\nu_1(z))$ for $y>0$. Then by Fubini theorem and \eqref{e:l2.7}, we obtain
\begin{align*}
I_2 &= 2r^{-1}\left| \int_1^\infty \int_{y/r}^{\infty} \cos(y) (-d\nu_1(z)) dy \right| = 2r^{-1}\left| \int_{1/r}^\infty \int_1^{rz} \cos(y)dy  (-d\nu_1(z)) \right| \\
& \le 4r^{-1}\left| \int_{1/r}^\infty (-d\nu_1(z)) \right| = 4r^{-1} \nu_1(1/r) \le 4c_1 \ell(r).
\end{align*}
Similarly, we also have that $I_3 \le 4c_1\ell(r)$. Therefore, we get \eqref{e:simple} in this case.

\smallskip

For the case $\psi(\xi)=\vp(|\xi|^2)$ for a Bernstein function $\vp$, we use  \cite[Lemma 5.13]{GRT19} and  \eqref{smallKnu}, and obtain that for any $1<k\le \lambda$ and $r\ge1$,
\begin{equation*}
|\psi(k r)-\psi(r)|= \int_{r^2}^{(kr)^2}\vp'(u)du \le r^{-d} \int_0^{(\lambda r)^2}u^{d/2}\vp'(u)du \le c_6\lambda^{d}\ell(\lambda r) \le 
c_7 \lambda^{d+\alpha_2}\ell( r).
\end{equation*}
This completes the proof. \qed

Now, we first consider the case when  {\bf (S-2)} holds. Recall that
$
\ll(r):=\sup_{u \in [1,r]} \ell(u)$ 
and $\ell^{-1}$ is the right continuous inverse of $\ll$ (see \eqref{ell_inv}).
Since {\bf (S-2)} holds, we get that $\lim_{r \to \infty} \ll(r) = \infty$ and there exists a constant $C_3 \ge 1$ such that
\begin{equation}\label{e:compS2}
\ell(r) \le \ll(r) \le C_3 \ell(r) \qquad \text{for all} \;\; r > 2.
\end{equation}
Note that in this case, by Proposition \ref{p:existence}, $p(t, 0)<\infty$ for all $t>0$. Here, we give the small time estimates for $p(t,0)$ under {\bf (S-2)}.

\begin{lemma}\label{l:ondiag1}
Assume that {\bf (S-2)} holds. Then, there exists a constant $c_1>0$ such that
\begin{align*}
p(t,x) \le p(t,0) \le c_1 \big[\ell^{-1}(a_1/t)\big]^d \exp\big(- b_1 t h(\ell^{-1}(a_1/t)^{-1})\big),
\end{align*}
for all $0<t \le t_1$ and $x \in \R^d$ where $a_1:=2dC_3/C_0$, $b_1:=C_0/(4C_2C_3)$ and $t_1:=a_1/\ll(3)$.
\end{lemma}
\pf
Let $a_1:=2dC_3/C_0$ and $t_1:=a_1/\ll(3)$. Then, $\ell^{-1}(a_1/t) \ge 3$ for all $t \in (0,t_1]$.
By Fourier inversion theorem, \eqref{e:eq1}, integration by parts and the change of variables $s=\Phi(r)$, we have that for all $t \in (0, t_1]$,
\begin{align*}
p(t,x) &= (2 \pi)^{-d}\int_{\R^d} e^{-i\ip{\xi, x}}e^{-t \psi(\xi)} d\xi \le c \int_0^{\infty} e^{-C_0t \Phi(r)}r^{d-1}dr \\
&\le ct \int_0^{\infty}r^d e^{-C_0t \Phi(r)}  \Phi'(r) dr = ct \int_0^{\infty}\Phi^{-1}(s)^d e^{-C_0ts}  ds  \\
& \le ct + ct \int_{\Phi(1)}^{\Phi(\ell^{-1}(a_1/t))} \Phi^{-1}(s)^d e^{-C_0 ts} ds + ct\int_{\Phi(\ell^{-1}(a_1/t))}^{\infty} \Phi^{-1}(s)^d e^{-C_0 ts} ds\\
& =: ct+ I_1+I_2.
\end{align*}

Observe that for $\Phi(2)<v \le u$, we have
\begin{align*}
u-v &= \Phi(\Phi^{-1}(u))-\Phi(\Phi^{-1}(v)) = \int_{\Phi^{-1}(v)}^{\Phi^{-1}(u)} k^{-1}\ell(k)dk \\
& \ge C_3^{-1}\int_{\Phi^{-1}(v)}^{\Phi^{-1}(u)} k^{-1}\ll(k)dk \ge C_3^{-1}\ll(\Phi^{-1}(v)) \log \frac{\Phi^{-1}(u)}{\Phi^{-1}(v)}. 
\end{align*}
Thus, for all $\Phi(2)<v \le u$, we have that (cf. ~ Section 3.10 in \cite{BGT})
\begin{align}\label{e:l_ond1}
\frac{\Phi^{-1}(u)}{\Phi^{-1}(v)} \le \exp\left(C_3\frac{u-v}{\ll(\Phi^{-1}(v))}\right).
\end{align}
Then, by \eqref{e:l_ond1} and the definition of $a_1$, we get
\begin{align*}
I_2 &= ct \big[\ell^{-1}(a_1/t)\big]^d \int_{\Phi(\ell^{-1}(a_1/t))}^{\infty} \left(\frac{\Phi^{-1}(s)}{\Phi^{-1}(\Phi(\ell^{-1}(a_1/t)))}\right)^d e^{-C_0ts}ds \\
&\le c \big[\ell^{-1}(a_1/t)\big]^d \int_{\Phi(\ell^{-1}(a_1/t))}^{\infty} t\exp\left(-\frac{dC_3t\Phi(\ell^{-1}(a_1/t))}{a_1}+ \frac{dC_3t s}{a_1} -C_0ts\right) ds \\
& \le c \big[\ell^{-1}(a_1/t)\big]^d \exp\big(-\frac{C_0}{2}t\Phi(\ell^{-1}(a_1/t))\big)\int_{\Phi(\ell^{-1}(a_1/t))}^{\infty} \left(-\frac{d}{ds}\exp\big(-\frac{C_0}{2}ts\big)\right) ds \\
& \le c \big[\ell^{-1}(a_1/t)\big]^d \exp\big(-C_0t \Phi(\ell^{-1}(a_1/t))\big).
\end{align*}

On the other hand, define $ g(r):=r^d \exp\big(-\frac{C_0}{2C_3}t\Phi(r)\big)$ for $r \ge 1$. Then, we have
$$
g'(r)=\big(d-\frac{C_0}{2C_3}t\ell(r)\big)r^{d-1}\exp\Big(-\frac{C_0}{2C_3}t \Phi(r)\Big).
$$
It follows that $g$ is strictly increasing on $[1, \ell^{-1}(a_1/t))$. Therefore, we obtain
\begin{align*}
I_1&\le ct \int_{\Phi(1)}^{\Phi(\ell^{-1}(a_1/t))} g(\Phi^{-1}(s)) ds \le 2ct \int_{\Phi(\ell^{-1}(a_1/t))/2}^{\Phi(\ell^{-1}(a_1/t))} g(\Phi^{-1}(s))ds \\
& \le c\big[\ell^{-1}(a_1/t)\big]^d \int_{\Phi(\ell^{-1}(a_1/t))/2}^{\Phi(\ell^{-1}(a_1/t))} \left(-\frac{d}{ds}\exp\Big(-\frac{C_0}{2C_3}ts\Big)\right)ds \\
& \le c\big[\ell^{-1}(a_1/t)\big]^d\exp\Big(-\frac{C_0}{4C_3}t \Phi(\ell^{-1}(a_1/t))\Big).
\end{align*}
We also have that 
\begin{align*}
I_1 \ge ct \int_{\Phi(1)}^{\Phi(3)} \Phi^{-1}(s)^d \exp(-C_0ts) ds \ge  ct.
\end{align*}
Finally, we deduce the result from \eqref{e:eq2}.
\qed

\begin{lemma}\label{l:offdiag}
Assume that {\bf (S-2)} holds. Let $a_1, b_1$ and $t_1$ be the positive constants in Lemma \ref{l:ondiag1}. Then, there exists a constant $c_1>0$ such that
\begin{align*}
p(t,x) \le c_1t|x|^{-d} \ll(|x|^{-1}) \exp\big(-b_1 t h(|x|)\big),
\end{align*}
for all $0<t \le t_1$ and $x \in \R^d$ satisfying $[\ell^{-1}(a_1/t)]^{-1} \le |x| \le 1/2$.
\end{lemma}
\pf
Fix $x \in \R^d$ satisfying $[\ell^{-1}(a_1/t)]^{-1} \le |x| \le 1/2$ and let $r=|x|$. By \cite[(5.4)]{GRT19}, the mean value theorem, \eqref{e:eq1} and Lemma \ref{l:simple}, for $0<t\le t_1$, we have
\begin{align}\label{upperbase}
&r^dp(t, x) \le c \int_{\R^d} \left(e^{-t \psi(|z|/r)} - e^{-t \psi(2|z|/r)} \right) e^{-|z|^2/4} dz\nn \\
&\quad \le ct \int_{\R^d} \sup_{|z| \le y \le 2|z|}e^{-t \psi(y/r)} \big|\psi(2|z|/r)-\psi(|z|/r)\big| e^{-|z|^2/4} dz \nn\\
&\quad \le ctr^d + ct \int_r^1 e^{-C_0t \Phi(u/r)} \ell(u/r)u^{d-1}du +  ct \int_1^{\infty} e^{-C_0t \Phi(u/r)} \ell(u/r)e^{-u^2/4}u^{d-1}du \nn\\
&\quad =: ctr^d +  I_1+I_2.
\end{align}

Since $\ell$ satisfies $\WUS^{\infty}(\alpha_2, 1)$ and  $\Phi$ is increasing, we have
\begin{align*}
I_2 \le ct \ell(1/r)\exp\big(-C_0t\Phi(1/r)\big) \int_1^{\infty} e^{-u^2/4}u^{d-1+\alpha_2} du \le ct \ll(1/r)\exp\big(-C_0t\Phi(1/r)\big).
\end{align*}

On the other hand, define $m(u):=u^{d-1/2} \exp\big(-\frac{C_0}{4C_3}t\Phi(u/r)\big)$ for $r>0$. Then,  for all $u \in (r,1)$, since $1/r \le \ell^{-1}(a_1 /t)$ and $a_1=2dC_3/C_0$, we get
\begin{align*}
m'(u)\exp\big(\frac{C_0}{4C_3}t\Phi(u/r)\big) = \left(d-\frac{1}{2}-\frac{C_0}{4C_3}t \ell(u/r)\right)u^{d-3/2} \ge \left(\frac{d}{2}-\frac{C_0}{4C_3}t\ll(1/r)\right)u^{d-3/2} \ge 0.
\end{align*}
It follows that $m(u)$ is increasing on $[r,1]$.
In particular, we have that 
\begin{align}\label{bigger}
m(1) = \exp\big(-\frac{C_0}{4C_3}t\Phi(1/r)\big) \ge r^{d-1/2}\exp\big(-\frac{C_0}{4C_3}t \Phi(1) \big) \ge cr^d.
\end{align}
Since $C_3 \ge 1$, we obtain
\begin{align*}
I_1 &\le ct \ll(1/r) \int_r^1 u^{d-1/2}\exp\big(-C_0t\Phi(u/r)\big)u^{-1/2}du \le  ct \ll(1/r) \int_r^1 m(u)u^{-1/2}du\\
& \le ct \ll(1/r) m(1) \int_r^1 u^{-1/2}du  \le ct \ll(1/r) \exp\big(-\frac{C_0}{4C_3} t\Phi(1/r)\big).
\end{align*}
Therefore, we deduce the result from \eqref{upperbase}, \eqref{bigger} and \eqref{e:eq2}.
\qed

In view of Lemma \ref{l:ondiag1} and Lemma \ref{l:offdiag}, we define for $a,r,t>0$,
\begin{align}\label{e:deftheta}
\theta_a(r,t):= r \vee [\ell^{-1}(a/t)]^{-1}.
\end{align}
Note that both $r \mapsto \theta_a(r,t)$ and $t \mapsto \theta_a(r,t)$ are increasing, while $a \mapsto \theta_a(r,t)$ is decreasing. 

\begin{prop}\label{p:upper1}
Assume that {\bf (S-2)} holds. For all $T>0$, there exists a constant $c_1>0$ such that for all $(t,x) \in (0,T] \times \R^d$,
\begin{equation}\label{e:upper1}
p(t,x) \le c_1t \frac{K(\theta_{a_1}(|x|,t))}{\big[\theta_{a_1}(|x|,t)\big]^d} \exp\big(-b_1t h(\theta_{a_1}(|x|,t))\big),
\end{equation}
where $a_1$ and $b_1$ are the constants in Lemma \ref{l:ondiag1}.
\end{prop}
\pf
Choose any $x \in \R^d$, and  let $a_1$ and $t_1=a_1/\ell^*(3)$ be the constants in Lemma \ref{l:ondiag1}.

   We first assume that $t \le t_1$. If $|x| < [\ell^{-1}(a_1/t)]^{-1}$, then we have $\theta_{a_1}(|x|,t) = \ell^{-1}(a_1/t)^{-1}$ so that $\ell(\theta_{a_1}(|x|,t)^{-1})\ge a_1C_3^{-1}t^{-1}$ by \eqref{e:compS2}.  Hence,  we obtain \eqref{e:upper1} from \eqref{Kgenu},  \eqref{A}  and Lemma \ref{l:ondiag1}. Else if $[\ell^{-1}(a_1/t)]^{-1} \le |x| \le 1/2$, then \eqref{e:upper1} follows from \eqref{Kgenu},  \eqref{A}, \eqref{e:compS2} and Lemma \ref{l:offdiag}. Otherwise, if $|x| > 1/2$, then since $r\mapsto h(r)$ is decreasing,   we see that
$$
t \frac{K(\theta_{a_1}(|x|,t))}{\big[\theta_{a_1}(|x|,t)\big]^d} \exp\big(-b_1t h(\theta_{a_1}(|x|,t))\big) \asymp \frac{tK(|x|)}{|x|^d}.
$$
Thus, we get \eqref{e:upper1} from Proposition \ref{p:gub}.

Now, suppose that $t\in (t_1, T]$. In this case, we have that $\ell^{-1}(a_1/t) \asymp 1$. Therefore, if $|x| \ge [\ell^{-1}(a_1/t)]^{-1}$, then we get the result from Proposition \ref{p:gub}. Otherwise, if $|x| < [\ell^{-1}(a_1/t)]^{-1}$, then by the semigroup property, Lemma \ref{l:ondiag1} and \eqref{smallKnu}, we have that
\begin{align*}
p(t,x)&= \int_{\R^d} p(t_1/2, x-z) p(t-t_1/2, z) dz \le p(t_1/2,0) \int_{\R^d} p(t-t_1/2,z) dz \\
& \le c \asymp t \frac{K(\theta_{a_1}(|x|,t))}{\big[\theta_{a_1}(|x|,t)\big]^d} \exp\big(-b_1t h(\theta_{a_1}(|x|,t))\big).
\end{align*}
This completes the proof.
\qed

By combining Propositions \ref{p:glb} and  \ref{p:upper1}, we obtain the following two-sided heat kernel estimates under {\bf (S-2)}.

\begin{cor}\label{c:HKEU}
Assume that {\bf (S-2)} holds. For all $T>0$, there exists a constant $c_1>1$ such that for every fixed $\delta>0$, we have that for all $(t,x) \in (0,T] \times \R^d$,
\begin{align}\label{hkes2a}
& c_1^{-1}t \nu(\theta_\delta(|x|,t)) \exp\big(-b_0t h(\theta_\delta(|x|,t))\big) \nn\\
& \qquad \qquad \le p(t,x) \le  c_1t \frac{K(\theta_{a_1}(|x|,t))}{\big[\theta_{a_1}(|x|,t)\big]^d} \exp\big(-b_1t h(\theta_{a_1}(|x|,t))\big),
\end{align}
where $b_0$ is the constant in Proposition \ref{p:glb}, and $a_1$ and $b_1$ are the constants in Lemma \ref{l:ondiag1}.
\end{cor}
\pf The upper bound follows from Proposition \ref{p:upper1}. On the other hand, since $p(t,\cdot)$ is radially non-increasing and  $\theta_{\delta}(|x|,t) \ge |x|$ for all $\delta, t>0$ and $x \in \R^d$, we deduce the lower bound from Proposition \ref{p:glb}.
\qed

\begin{remark}
{\rm
If $\ell$ satisfies $\WLS^{\infty}(\alpha,1)$ for some $\alpha>0$, then  $\ell(r) \asymp \Phi(r)$ for $r \ge 1$. (See, \cite[Theorem 2.6.1]{BGT}.)  Therefore, when $\ell$ satisfies $\WLS^{\infty}(\alpha,1)$ for some $\alpha>0$, we see that the estimate \eqref{hkes2a} can be expressed as follows: For every $(t,x) \in (0,T] \times \R^d$,
\begin{align*}
c_1^{-1}\Phi^{-1}(1/t)^d \wedge t \nu(|x|) \le p(t,x) \le c_1\Phi^{-1}(1/t)^d \wedge t \frac{K(|x|)}{|x|^d}.
\end{align*}
Hence, if {\bf (C)} further holds, then we see from \eqref{smallKnu} and \eqref{largeKnu} that   $p(t, x) \asymp \Phi^{-1}(1/t)^d \wedge t \nu(|x|)$ for $(t,x) \in (0,T] \times \R^d$. In view of  \eqref{e:asym}, \eqref{e:eq1} and \eqref{e:eq2}, this coincides with the main result in \cite{BGR14a}.
}
\end{remark}

In the rest of this section, we assume that {\bf (S-1)} holds. Then, by Proposition \ref{p:existence}, we have that $p(t,0) = \infty$ for all sufficiently small $t$. Recently, some general estimates for such type of heat kernels were established in \cite{GRT19}. Using that results, we obtain the heat kernel estimates in analogous form to \eqref{hkes2a}.

\begin{prop}\label{p:HKEl}
Assume that {\bf (S-1)} holds. Then, there exist constants $t_0, c_1>0$ such that for all $(t,x) \in (0,t_0] \times (\R^d \setminus \{0\})$,
\begin{align}\label{e:HKEl}
p(t,x) \le c_1 t |x|^{-d}K(|x|) \exp\big(-t \psi(|x|^{-1})\big).
\end{align}
\end{prop}
\pf
Let $\omega(r) = K(1) \1_{\{0<r \le 1\}}(r) + K(r^{-1})\1_{\{r>1\}}(r)$ for $r>0$ where $\1_A$ denotes the indicator function on a set $A$.
By \eqref{smallKnu}, Lemma \ref{l:simple}, {\bf (A)} and  {\bf (S-1)}, there exists a constant $c_0>0$ such that $c_0 \omega(r)$ satisfies the assumptions (5.7) and (5.8) in \cite{GRT19}. Therefore, by \cite[Proposition 5.6]{GRT19}, there exist $t_0,c_1>0$ such that for all $t \in (0, t_0]$ and $0<|x| < 1$, the estimate \eqref{e:HKEl} holds. Moreover, for $t \in (0, t_0]$ and $|x| \ge 1$, we have that $e^{-t \psi(|x|^{-1})} \asymp 1$. Then, we get the result from Proposition \ref{p:gub}. 
\qed

\begin{cor}\label{c:HKEB}
Assume that {\bf (S-1)} holds. For all $T>0$, there exist constants $c_1, b_2>0$ such that
\begin{equation}\label{e:cHKEB}
c_1^{-1} t\nu(|x|)\exp\big(-b_0t h(|x|)\big) \le p(t,x) \le c_1 t |x|^{-d}K(|x|) \exp\big(-b_2t h(|x|)\big),
\end{equation}
for all $(t,x) \in (0,T] \times (\R^d \setminus \{0\})$ where $b_0$ is the constant in Proposition \ref{p:glb}.
\end{cor}
\pf
By Propositions \ref{p:glb} and \ref{p:HKEl}, \eqref{e:asym} and induction, it suffices to prove the upper bound in \eqref{e:cHKEB} for $t \in (t_0, 2t_0]$ and $x \in \R^d \setminus \{0\}$, where $t_0$ is the constant in Proposition \ref{p:HKEl}. If $|x| \ge 1$, then $\exp\big(-cth(|x|)\big) \asymp 1$ for each fixed  constant $c>0$ so that the assertion holds by Proposition \ref{p:gub}. Suppose that $|x| < 1$. Without loss of generality, we may assume that
$2b_2 \le b_0$.  Then, by the semigroup property, \eqref{smallKnu}, the induction hypothesis, monotonicity of $p(t, \cdot)$ and Proposition \ref{p:gub}, we get
\begin{align*}
p(t,x) &= \int_{B(x,1)} p(t/2, x-z) p(t/2, z) dz + \int_{\R^d \setminus B(x,1)} p\big(\frac{b_2}{b_0}t, x-z) p\big(\frac{b_0-b_2}{b_0}t, z\big) dz \\
&\le c_1 \int_{B(x,1)} t^2 \nu(|x-z|)\exp\big(-\frac{b_0b_2t}{2b_0}h(|x-z|)\big)  \nu(|z|)\exp\big(-\frac{b_0b_2t}{2b_0}h(|z|)\big) dz  \\
&\quad +  p\big(\frac{b_2}{b_0}t, 1\big)\int_{\R^d} p\big(\frac{b_0-b_2}{b_0}t, z\big) dz \\
& \le c_2 \int_{\R^d} p\big(\frac{b_2}{2b_0}t, x-z\big)p\big(\frac{b_2}{2b_0}t, z\big)dz +  p\big(\frac{b_2}{b_0}t, 1\big) \\
&\le (c_2+1)p\big(\frac{b_2}{b_0}t, x\big) \le c_3t|x|^{-d} K(|x|) \exp\big(-\frac{b_2^2}{b_0}th(|x|)\big).
\end{align*}
\qed

\section{Boundary Harnack principle with explicit decay}\label{s:decay}

In this section, we investigate the boundary behaviour of the process via the renewal function $V$ of $Y$ and the tail of its L\'evy measure. Throughout this section, we assume that {\bf (B)} holds. For an open set $D\subset \R^d$, the first exit time is denoted by $\tau_D:=\inf\{t>0: \, Y_t\notin D\}$. We give the probabilistic definition of a (regular) harmonic function.

\begin{definition}\label{d:har}
{\rm(i) A function 
$u:\bR^d\to \bR$ 
is said to be harmonic in an open set $D\subset \bR^d$ with respect to $Y$ if
for every open set $B$ whose closure is a compact subset of $D$, $\bE_x[|u(Y_{\tau_B})|]<\infty$ and 
$u(x)=\bE_x[u(Y_{\tau_B})]$ for every $x\in B$.

\noindent
(ii) A function 
$u:\bR^d\to \bR$ 
is said to be 
regular harmonic in an open set $D\subset \bR^d$ 
with respect to $Y$ if $\bE_x[|u(Y_{\tau_D})|]<\infty$ and 
$u(x)=\bE_x[u(Y_{\tau_D})]$ for every $x\in D.$
}\end{definition}

Here, we provide the precise definition of the renewal function $V$ of $Y$. Let $Y^d$ be the last coordinate of $Y$, $M_t=\sup_{s\le t}Y^d_s$ and $\scL_t$ be the local time {at $0$} for
{$M_t-Y^d_t$}, the last coordinate of $Y$ reflected at the supremum.
Define the ascending ladder-height process as
$H_t = Y_{\scL^{-1}_t}^{d} = M_{\scL^{-1}_t}$ where $\scL^{-1}$ is the right continuous inverse of $\scL$.
Then, the renewal function $V$ is defined as
\begin{equation*}\label{e:defV}
V(s) = \int_0^{\infty}\P(H_t \le s)dt, \quad
s \in \R.
\end{equation*}
Since the process $Y$ is isotropic unimodal, there are several known properties for the renewal function. (See, \cite[Theorem 1.2]{Si80}, \cite[p.74]{Be96} and \cite[Section 1.2]{BGR14b}.)

\begin{lemma}\label{l:pV}
(i) $V$ is strictly increasing, $V(s) = 0$ if $s<0$ and $\lim_{s \to \infty}V(s)=\infty$.\\
(ii) $V$ is subadditive; that is,
\begin{equation*}
 V(s+r)\le V(s)+V(r) \quad \text{for all} \;\; s,r \in \R.
\end{equation*}
(iii) $V$ is absolutely continuous and harmonic on $(0, \infty)$ for the process $Y_t^d$. Also, $V'$ is a positive harmonic function for $Y_t^d$ on $(0, \infty)$. 
\end{lemma}

According to \cite[Proposition 2.4]{BGR15}, the relation \eqref{e:asym} can be extended to include the renewal function. That is, there exist comparison constants which are only depend on the dimension $d$ and $\kappa_1$ and $\kappa_2$ in \eqref{A}  such that
$
 h(r) \asymp \psi(r^{-1})  \asymp  [V(r)]^{-2}$ for all $r>0$. Then, by Lemmas \ref{l:asymLh} and \ref{l:propPhi}, we have that 
\begin{align}\label{e:asym2}
L(r) \asymp  h(r) \asymp \psi(r^{-1}) \asymp \Phi(r^{-1}) \asymp  [V(r)]^{-2} \qquad \text{for all}\;\; 0<r \le 1.
\end{align}
In particular, by \eqref{e:asym2} and Lemma \ref{l:propPhi}, there are constants $c_1, c_2,c_3,c_4>0$ such that
\begin{align}\label{e:vws}
c_1 \left(\frac{R}{r}\right)^{\alpha_1/2} \le \frac{V(R)}{V(r)} \le c_2 \left(\frac{R}{r}\right)^{(\alpha_2/2) \vee (1/4)} \quad \text{for all} \;\; 0<r \le R \le 1.
\end{align}
and
\begin{align}\label{e:Lws}
c_3 \left(\frac{R}{r}\right)^{\alpha_1} \le \frac{L(r)}{L(R)} \le c_4 \left(\frac{R}{r}\right)^{\alpha_2 \vee (1/2)} \quad \text{for all} \;\; 0<r \le R \le 1.
\end{align}

\begin{prop}\label{p:v1} 
The renewal function $V$ is twice-differentiable on $(0,\infty)$, and there exists $c_1>0$ such that 
\begin{align*}
|V''(r)|\leq c_1\frac{V'(r)}{r\wedge 1}\quad \text{and} \quad
 V'(r)\leq c_1\frac{V(r)}{r\wedge 1}, \quad r>0.
\end{align*}
\end{prop}
\pf
Since {\bf (A)} and {\bf (B)} hold, the scale-invariant Harnack inequality holds for $Y$. (See, \cite[Theorem 1.9]{GK18}.) Then, the results follows from \cite[Theorem 1.1]{KuR16} and Lemma \ref{l:pV}(iii).
\qed

Define $w(x):=V((x_d)^+)$ for $x \in \R^d$ and let $\H:=\{x=(\widetilde{x},x_d) \in \R^d:x_d>0 \}$ the upper half-space. Since the renewal function $V$ is harmonic on $(0, \infty)$ for $Y^d$, by the strong Markov property, $w$ is harmonic in $\H$ with respect to $Y$.

\begin{prop}\label{p:v2}
For all $\lambda>0$, there exists $c_1=c_1(d,\lambda)>0$ such that for any $r>0$,
\begin{align*}
\sup_{\{x\in \R^d\, :\, 0<x_d\le \lambda r\}}\int_{B(x, r)^c} w(y)\nu(|x-y|)dy\le c_1V(r)^{-1}.
\end{align*}
\end{prop}
\pf
See, the proof of \cite[Proposition 3.2]{GKK20}.
\qed

Denote $C_{\infty}^2(\R^d)$ by the set of all twice-differentiable functions in $\R^d$ vanishing at infinity. We define an operator $\sL_Y$ as follows: for $\eps>0$ and $x \in \R^d$,
\begin{align*}
\gener_Y^{\eps} f(x)&:=\int_{B(x, \varepsilon)^c} (f(y)-f(x))\nu(|x-y|)dy, \label{e:gY}\\
\gener_Y f(x)&:=P.V.\int_{\Rd} (f(y)-f(x))\nu(|x-y|)dy =\lim_{\varepsilon\downarrow 0}\gener_Y^{\eps} f(x),\nn\\
{\cal D} (\gener_Y)&:=\left\{f\in C_{\infty}^{2}(\R^d):P.V.\int_{\Rd} (f(y)-f(x))\nu(|x-y|)dy\,\mbox{ exists and is finite.}\right\}.\nn
\end{align*}

\begin{thm}\label{t:v3}
For any $x\in \H$, $\gener_{Y} w(x)$ is well-defined and $\gener_Y w(x)=0$.
\end{thm}
\pf
By Propositions \ref{p:v1} and \ref{p:v2}, using \cite[Lemma 2.3, Theorem 2.11]{Ch09}, the proof is essentially the same as the one given in \cite[Theorem 3.3]{GKK20}. Hence, we omit it.
\qed

\begin{lemma}\label{L:Main}
Let $D$ be a $C^{1,1}$ open set in $\R^d$ with characteristics $(R_0, \Lambda)$. For any $Q \in \partial D$ and $r >0 $, we define
\[
h_r(y)=h_{r,Q}(y) := V(\delta_D(y)) \1_{ D \cap B(Q, r)}(y).\\
\]
Then, there exist $R_1=R_1(R_0,\Lambda,\psi, d) \in (0,(R_0 \wedge 1)/2] $ and $c_1=c_1(R_0,\Lambda,\psi, d)>1$ independent of $Q$ such that for every $r \in (0,R_1)$, $\sL_Y h_r$ is well defined in $D\cap B(Q,
r/4)$ and
\begin{align*}
|\sL_Y h_r(x)|\le
\frac{c_1}{V(r)} \le c_1^2 L(r)^{1/2}\quad \text{ for all } x \in D\cap B(Q, r/4)\, .
\end{align*}

\end{lemma}

\pf
Since the case of $d=1$ is easier, we only give the proof for $d\ge 2$. Fix $Q \in \partial D, r \in(0,(R_0 \wedge 1)/2)$ and $x \in D \cap B(Q,r/4)$. Let $z \in \partial D$ be the point satisfying $\delta_D(x) = |x-z|$ and denote $\Gamma_{z}$ and $CS_{z}$ by the $C^{1,1}$ function and orthonormal coordinate system determined by $z$, respectively. (See, Definition \ref{c11}.) Henceforth, we use the coordinate system $CS_{z}$. Hence, we have $z=0$, $x=(\wt 0, x_d)$ and $D \cap B(z, R_0) = \{y = (\wt y, y_d)\in B(0, R_0)\;\; \text{in} \;\; CS_z : y_d > \Gamma_z(\wt y) \}.$  Since $D$ is a $C^{1,1}$ open set, it satisfies the inner and outer ball conditions. Thus, we may assume that 
\begin{align*}
A_1:=\{ y = (\wt y, y_d) \;\;\text{in}\;\; CS_{z} : |y|<R_0,\; y_d > \phi(\wt y) \} \subset D,
\end{align*}
and
\begin{align*}
A_2:=\{ y = (\wt y, y_d) \;\;\text{in}\;\; CS_{z} : |y|<R_0,\; y_d < - \phi(\wt y) \} \subset D^c,
\end{align*}
where $\phi:\R^{d-1} \to \R$ is defined by $\phi(\wt y):= 1-\sqrt{1-|\wt y|^2}.$

Let $E:= \{ y = (\wt y, y_d) : |\wt y|< r/2, |y_d|< r/2  \}$, $E_1:= \{ y \in E : y_d > 2\phi(\wt y) \}$ and $E_2:= \{ y \in E : y_d < -2\phi(\wt y) \}$. We also let $w_z(y):=V((y_d)^+)$. By Theorem \ref{t:v3}, we get $\sL_Y w_z(x) = 0$. Since $h_r(x)=w_z(x)$ and $h_r(y)=w_z(y)=0$ for $y \in E_2$, we have
\begin{align}\label{e:LY}
\big|\sL_Y h_r(x)\big|&=\big|\sL_Y( h_r-w_z)(x)\big|\nn\\
&=\bigg|\lim_{\epsilon \downarrow 0} \int_{|y-x|>\epsilon} \Big(\big(h_r(y)-w_z(y)\big)-\big(h_r(x)-w_z(x)\big)\Big)\nu(|x-y|)dy\bigg|\nn\\
& = \bigg|\lim_{\epsilon \downarrow 0} \int_{|y-x|>\epsilon} \big(h_r(y)-w_z(y)\big)\nu(|x-y|)dy\bigg|\nn\\
& \le \limsup_{\epsilon \downarrow 0}\bigg(\int_{E_1,|y-x|>\epsilon}+\int_{E\setminus (E_1 \cup E_2),|y-x|>\epsilon}+\int_{E^c} \big|h_r(y)-w_z(y)\big|\nu(|x-y|)dy\bigg)\nn\\
&=: I_1+I_2+I_3. 
\end{align}

First, since $|h_r(y)| \le V(r)$, using Lemma \ref{l:asymLh}, \eqref{e:asym2}, \eqref{e:vws} and Proposition \ref{p:v2}, we have
\begin{align*}
I_3 & \le \int_{B(x,r/2)^c} \big(|h_r(y)|+|w_z(y)|\big)\nu(|x-y|)dy \le cV(r)L(r/2) + cV(r)^{-1} \le cV(r)^{-1}.
\end{align*}
Next, we note that for $y \in E \setminus (E_1 \cup E_2)$,
\begin{align*}
\delta_D(y) \le \delta_{A_2}(y) \le |y_d+\phi(\wt y)|\le 3\phi(\wt y) \le 3 |\wt y|^2 \le 3|\wt y|,
\end{align*}
and hence by subadditivity of $V$, we obtain
\begin{align*}
|h_r(y)|+|w_z(y)| \le V(3|\wt y|)+V(2|\wt y|)\le 5V(|\wt y|).
\end{align*}
Since $1-\sqrt{1-l^2} \le l^2$ for $0\le l<1$, we have for $0<s<1$,
\begin{align*}
m_{d-1}\big(\{y \in E \setminus (E_1 \cup E_2) : |\wt y | =s \}\big) \le  m_{d-1}\big(\{y : |\wt y|=s, -2|\wt y|^2 \le y_d \le 2|\wt y|^2  \}\big)  \le c s^d,
\end{align*}
where $m_{d-1}(dx)$ is the $(d-1)$-dimensional Lebesgue measure. From these observations, using \eqref{e:doubling}, \eqref{Kgenu}, the definitions of $K$ and $h$, \eqref{e:asym2} and \eqref{e:vws},  since $\alpha_2<2$,  we get
\begin{align*}
I_2 & \le c \int_0^r \int_{|\wt y|=s, y \in E \setminus (E_1 \cup E_2) } m_{d-1}(dy) \; V(s) \nu(s)  ds \le c \int_0^r V(s)\nu(s)s^d ds\\
&  \le c \int_0^r V(s)h(s)ds \le  c \int_0^r V(s)^{-1} ds \le  c V(r)^{-1} \int_0^r \left(\frac{r}{s}\right)^{(\alpha_2/2) \vee (1/4)}ds \le c V(r)^{-1}.
\end{align*}

Lastly, to estimate $I_1$ we first claim that
\begin{align}\label{e:distasym}
\delta_D(y) \asymp y_d, \;\; |\delta_D - y_d| \le 2|\wt y|^2 \quad \text{for all} \;\; y \in E_1.
\end{align}
Indeed, for any $y \in E_1$, if $0<y_d \le \delta_D(y)$, then we have
\begin{align*}
\delta_D(y) - y_d \le \delta_{A_2}(y)-y_d \le \phi(\wt y) \le (y_d/2) \wedge (2|\wt y|^2).  
\end{align*}
Otherwise, if $\delta_D(y) < y_d$, then we have
\begin{align*}
y_d - \delta_D(y) \le y_d - \delta_{A_1}(y) = y_d - 1 + \sqrt{|\wt y|^2+(1-y_d)^2} = \frac{|\wt y|^2}{(1-y_d) + \sqrt{|\wt y|^2 + (1-y_d)^2}}.
\end{align*}
 Hence, since $|\wt y|, |y_d|<r/2<1/4$, we get $y_d-\delta_D(y) \le (2/3)|\wt y|^2 \le (4/3)\phi(\wt y)<(2/3)y_d$. Therefore, \eqref{e:distasym} holds.

Recall that by Lemma \ref{l:pV}(iii), $V'$ is a harmonic function for $Y^d_t$ on $(0,\infty)$. Since the scale-invariant Harnack inequality holds for $Y$ (see, \cite[Theorem 1.9]{GK18}), by \eqref{e:distasym}, we deduce that for every $y \in E_1$, 
\begin{equation*}
|h_r(y)-w_z(y)| \le \Big( \sup_{\delta_D(y) \wedge y_d \le l \le \delta_D(y) \vee y_d} V'(l) \Big) |\delta_D(y)-y_d|\le c V'(y_d)|\wt y|^2,
\end{equation*}
for some constant $c>0$ independent of choice of $Q, r$ and $x$.  Hence, we obtain
\begin{align*}
I_1 & \le c \int_0^r \int_{2\phi(k)}^{k} V'(y_d) \nu\big(\sqrt{k^2+(y_d-x_d)^2}\big)k^d dy_d  dk\\
& \quad + c \int_0^r \int_{k}^{r} V'(y_d) \nu\big(\sqrt{k^2+(y_d-x_d)^2}\big)k^d dy_d  dk=:I_{1,1}+I_{1,2}.
\end{align*}
By the monotonicity of $\nu$, \eqref{Kgenu}, the definition of $h$, \eqref{e:asym2} and \eqref{e:vws}, since $\alpha_2<2$, we have
\begin{equation*}
I_{1,1} \le c \int_0^r \int_0^k V'(y_d)dy_d \nu(k)k^d dk\le c \int_0^r V(k)h(k) dk \le cV(r)^{-1} \int_0^r \frac{V(r)}{V(k)}dk \le cV(r)^{-1}.
\end{equation*}
Besides, set $\rho:=4^{-1}((2-\alpha_2) \wedge 1) \in (0, 1/4]$. 
By Proposition \ref{p:v1},  \eqref{Kgenu}, the definition of $h$, \eqref{e:asym2} and \eqref{e:vws}, since $\sqrt{a^2+b^2} \ge (a+b)/\sqrt2$ for all $a,b \ge 0$, we have that
\begin{align}
 I_{1,2} & \le c \int_0^r \int_{k}^{r} \frac{V(y_d)k^{-\rho}}{y_d^{1-\rho}} \nu\big(\sqrt{k^2+(y_d-x_d)^2}\big)k^d dy_d  dk\nn\\
 &\le  c \int_0^r \int_{k}^{r} \Big(\sup_{u \in [y_d,1]} \frac{V(u)}{u^{1-\rho}}\Big) \nu\big(\sqrt{k^2+(y_d-x_d)^2}\big)\big(\sqrt{k^2+(y_d-x_d)^2}\big)^{d-\rho} dy_d  dk  \nn\\
  &\le  c \int_0^r \int_{k}^{r} \Big( \sup_{u \in [y_d,1]} \frac{V(u)}{u^{1-\rho}}\Big)\frac{(k+|y_d-x_d|)^{-\rho}}{V(k+|y_d-x_d|)^2}  dy_d  dk. \label{e:arrange}
\end{align}
By  \eqref{e:vws}, since $1-\rho>(\alpha_2/2)\vee (1/4)$, we see that for all $u \in [y_d,1]$,
\begin{equation}\label{e:3.7}
\frac{V(u)}{u^{1-\rho}} \le  \frac{cV(y_d)}{y_d^{1-\rho}} \Big(\frac{u}{y_d}\Big)^{-(1-\rho)+(\alpha_2/2)\vee (1/4)}\le \frac{cV(y_d)}{y_d^{1-\rho}}.
\end{equation} 
Since $s\mapsto \sup_{u \in [s,1]} (V(u)u^{\rho-1})$ and $s\mapsto s^{-\rho}V(s)^{-2}$ are non-increasing, by \eqref{e:arrange}, \cite[Lemma 4.4]{KSV14}, \eqref{e:3.7} and \eqref{e:vws}, we obtain  that 
\begin{align*}
I_{1,2} &\le c \int_0^{3r/2} \int_0^u \frac{V(s)}{s^{1-\rho}} ds \frac{u^{-\rho}}{V(u)^2}du \le c \int_0^{3r/2} \int_0^u \frac{ds}{s^{1-\rho}}  \frac{u^{-\rho}}{V(u)}du\\
&\le cV(r)^{-1} \int_0^{3r/2}  \frac{V(r)}{V(u)}du \le cV(r)^{-1}.
\end{align*}

Finally, we get from \eqref{e:LY} that $|\sL_Y h_r(x)| \le cV(r)^{-1}$. By \eqref{e:asym2}, this finishes the proof. \qed

For $l>0$, we define $D_{int}(l):= \{y \in D : \delta_D(y)>l \}$.

\begin{lemma}\label{l:explictdecay}
Let $D$ be a $C^{1,1}$ open set in $\R^d$ with characteristics $(R_0, \Lambda)$ and $R_1$ be the constant in Lemma \ref{L:Main}. Then, there exist constants $R_2=R_2(R_0,\Lambda, \psi,d) \in (0, R_1/16]$ and $c_1=c_1(R_0,\Lambda,\psi, d)>1$ such that for every $r \in (0, R_2]$ and $x \in D$ with $\delta_D(x)<r/2$,
\begin{align}\label{e:survivaltime}
\frac{c_1^{-1}}{L(\delta_D(x))^{1/2}L(r)^{1/2}} \le \E_x[\tau_{D\cap B(z,r)}]\le  \frac{c_1}{L(\delta_D(x))^{1/2}L(r)^{1/2}}.
\end{align}
and
\begin{align}\label{e:exitdist}
\P_x\left(Y_{\tau_{D \cap B(z,r)}} \in D_{int}(r/4) \right) \ge c_1^{-1} \left(\frac{L(r)}{L(\delta_D(x))} \right)^{1/2},
\end{align}
where $z \in \partial D$ is the point satisfying $\delta_D(x) = |x-z|$.
\end{lemma}
\pf
Let $R_1$ be the constant in Lemma \ref{L:Main}. Fix $r \in (0, R_2]$ and $x \in D$ with $\delta_D(x)<r/2$ where the constant $R_2 \in (0, R_1/16]$ will be selected later. Let $z \in \partial D$ be the point satisfying $\delta_D(x) = |x-z|$. As in Lemma \ref{L:Main}, we denote by $\Gamma_{z}:\R^{d-1} \to \R$ and $CS_{z}$ for a $C^{1,1}$ function and coordinate system with respect to $z$, respectively and hereinafter we use the coordinate system $CS_z$.

Denote by $U(s):= D \cap B(0,s)$ for $s>0$. Then, we define 
$$
u(y) = V(\delta_D(y)) \1_{U(R_1)}(y).
$$
Using Dynkin's formula and approximation argument, (see, \cite[Proposition 4.7]{KM17},) by Lemma \ref{L:Main}, there exists a positive constant $a$ independent of choice of $R_2$ and $x$ such that
\begin{align}\label{e:Dynkin}
\E_x\left[u(Y_{\tau_W})\right]-aL(R_1)^{1/2}\E_x[\tau_{W}] \le u(x) \le \E_x\left[u(Y_{\tau_W})\right]+aL(R_1)^{1/2}\E_x[\tau_{W}],
\end{align} 
for every open subset $W \subset U(R_1/8)$. 

\smallskip

Let
$$\sC_1: = \{(\wt y, y_d) : 2\Lambda |\wt y| < y_d, \; 0<|y|<R_0 \} \;\; \text{and} \;\; \sC_2: = \{(\wt y, y_d) : 4\Lambda |\wt y| < y_d, \; 0<|y|<R_0 \}.$$
We claim that $\sC_2 \subset \sC_1 \subset D$. Indeed, the first inclusion is obvious. Moreover, for all $y \in \sC_1$, 
$y_d-\Gamma_{z}(\wt y) \ge y_d - \Lambda|\wt y| \ge y_d/2>0.$ Hence, the second inclusion holds.

Observe that for $0<s \le R_1$ and $y \in \sC_2 \cap \partial U(s)$, we have
\begin{align}\label{cone}
s \ge \delta_D(y) \ge \delta_{\sC_1}(y) \ge c_0 s,
\end{align}
for some constant $c_0$ which only depends on $\Lambda$.
By the L\'evy system, \eqref{cone}, integration by parts, Lemma \ref{l:asymLh}, \eqref{e:asym2}, \eqref{e:Lws} and monotonicity of $V$,  for $0<4s<R_1$,
\begin{align}\label{e:Levy}
\E_x\big[u(Y_{\tau_{U(s)}})\big] & \ge \E_x\Big[u(Y_{\tau_{U(s)}}) : Y_{\tau_{U(s)}} \in \sC_2 \cap \big(U(R_1) \setminus U(2s)\big)\Big] \nn\\
& = \E_x \bigg[\int_0^{\tau_{U(s)}} \int_{\sC_2 \cap (U(R_1) \setminus U(2s))} V(\delta_D(y)) \nu(|Y_k-y|) dy dk  \bigg] \nn\\
& \ge c \E_x[\tau_{U(s)}] \int_{2s}^{R_1} V(k) \nu(k) k^{d-1} dk = c_1 \E_x[ \tau_{U(s)}] \int_{2s}^{R_1} (-L'(k))V(k) dk \nn\\
& = c_1 \E_x[\tau_{U(s)}]\bigg(L(2s)V(2s) - L(R_1)V(R_1) + \int_{2s}^{R_1} L(k)V'(k)dk\bigg) \nn\\
& \ge c_1\E_x[\tau_{U(s)}] \Big(c_2 L(s)^{1/2} - c_3L(R_1)^{1/2} \Big),
\end{align}
for some constants $c_1, c_2, c_3>0$ independent of  $s$.
Moreover, by the similar argument, we also have that
\begin{align}\label{need1}
\P_x\Big(Y_{\tau_{U(r)}} \in D_{int}(r/4)\Big) & \ge c\E_x \bigg[\int_0^{\tau_{U(s)}} \int_{\sC_2 \, \cap \,  (U(R_1) \setminus U(2r))} \nu(|Y_k-y|) dy dk  \bigg] \nn\\
&\ge c \E_x[\tau_{U(r)}] \int_{2r}^{R_1} \nu(k)k^{d-1}dk = c_1 \E_x[\tau_{U(r)}] \int_{2r}^{R_1} (-L'(k))dk \nn\\[3pt]
&  \ge c_1\E_x[\tau_{U(r)}] \big(c_4 L(r) - c_5L(R_1) \big)
\end{align}
and
\begin{align}\label{need2}
& \E_x\big[u(Y_{U(r)}) : Y_{U(r)} \in D_{int}(2r)\big]  \le \E_x \bigg[\int_0^{\tau_{U(s)}} \int_{U(R_1) \setminus U(2r)} V(\delta_D(y))\nu(|Y_k-y|) dy dk  \bigg] \nn\\
&   \le c \E_x[\tau_{U(r)}] \int_{2r}^{R_1} V(k)\nu(k/2)k^{d-1} dk  \le c \E_x[\tau_{U(r)}] \int_{r}^{R_1} L(k)^{-1/2}(-L'(k)) dk  \nn\\[3pt]
& \le c \E_x[\tau_{U(r)}] L(r)^{1/2} \le  c \E_x[\tau_{U(r)}] V(r)^{-1}.
\end{align} 
We used \eqref{e:doubling} in the third inequality.

For selected constants $a, c_1, c_2, c_3, c_4$ and $c_5$ in \eqref{e:Dynkin}, \eqref{e:Levy} and \eqref{need1}, we set 
$$
R_2=V^{-1}\left(\frac{c_1c_2}{2(a+c_1c_3)} V(R_1)\right) \wedge V^{-1}\left(\frac{c_4}{2c_5} V(R_1)\right) \wedge \frac{R_1}{4}.
$$
 Then, by \eqref{e:Dynkin} and \eqref{e:Levy}, we get
\begin{align*}
L(\delta_D(x))^{-1/2} &\asymp V(\delta_D(x)) = u(x) \ge \E_x[u(Y_{\tau_{U(r)}})]-aL(R_1)^{1/2}\E_x[\tau_{U(r)}]\\
&\ge \left(c_1c_2L(r)^{1/2} - (c_1c_3+a)L(R_1)^{1/2}\right) \E_x[\tau_{U(r)}]\ge 2^{-1}c_1c_2L(r)^{1/2}\E_x[\tau_{U(r)}].
\end{align*}
This proves the upper bound of (\ref{e:survivaltime}).

On the other hand, by \eqref{need1}, we get
\begin{align}\label{need11}
\P_x\Big(Y_{\tau_{U(r)}} \in D_{int}(r/4)\Big) & \ge 2^{-1}c_1c_4 \E_x[\tau_{U(r)}]L(r).
\end{align}
By \cite[Lemma 2.1]{BGR15} and \eqref{e:asym2}, there exists $c_6>0$ such that
\begin{align}\label{need3}
\P_x(Y_{\tau_{U(r)}} \in D) \le c_6 \E_x[\tau_{U(r)}] V(r)^{-2}.
\end{align}
Then, by \eqref{e:Dynkin}, \eqref{need2}, \eqref{need3}, \eqref{e:vws} and the monotonicity of $V$, we get
\begin{align*}
V(\delta_D(x))  & \le \E_x\big[u(Y_{U(r)}) : Y_{U(r)} \in D_{int}(2r)\big]+\E_x\big[u(Y_{U(r)}) : Y_{U(r)} \in D \setminus D_{int}(2r)\big] + c \E_x[\tau_{U(r)}]\\[2pt]
&\le c \E_x[\tau_{U(r)}]V(r)^{-1} + cV(2 r \sqrt{1+\Lambda^2}) \P_x\big(Y_{U(r)} \in D \setminus D_{int}(2r)\big) + c \E_x[\tau_{U(r)}] \\[2pt]
& \le c \E_x[\tau_{U(r)}]\big(V(r)^{-1} + V(2 r \sqrt{1+\Lambda^2}) V(r)^{-2} + 1 \big) \le c \E_x[\tau_{U(r)}] V(r)^{-1}.
\end{align*}
This proves the lower bound of \eqref{e:survivaltime} in view of \eqref{e:asym2}. Finally, we get \eqref{e:exitdist} from \eqref{need11}.
\qed

\bigskip

\section{Estimates of survival probability}

In this section, we obtain two-sided estimates for the survival probability $\P_x(\tau_D>t)$ which play a crucial role in factorization of the Dirichlet heat kernel. We first state the general two-sided estimates for the survival probability in balls which are recently established in \cite{GRT19}. 

\begin{prop}[{\cite[Proposition 5.2]{GRT19}}]
There exist positive constants $c_1, c_2, C_4$ and $C_5$ which only depend on the dimension $d$ such that for all $t,r>0$,
\begin{align}\label{exitball}
\begin{split}
c_1\exp\big(- \kappa_2C_4 t h(r)\big)& \le \bP_x(\tau_{B(x,r)}>t)\\
&\le \sup_{z \in B(x,r)} \P_z(\tau_{B(x,r)}>t) \le c_2 \exp\big(- \kappa_1C_5 t h(r)\big),
\end{split}
\end{align}
where $\kappa_1$ and $\kappa_2$ are constants in {\bf (A)}.
As a consequence, for all $r>0$,
\begin{align}\label{meanexit}
\E_x[\tau_{B(x,r)}] = \int_0^{\infty} \P_x(\tau_{B(x,r)}>s)ds  \asymp h(r)^{-1}.
\end{align}
\end{prop}
Note that  the last inequality in \eqref{exitball}, and \eqref{meanexit} were obtained for a large class of Feller processes in $\R^d$. See \cite[Corollaries 5.3, 5.8 and Theorem 5.9]{BSW13}.

\medskip

In the rest of this section, we assume that {\bf (B)} holds. Fix $T>0$ and $D$ a $C^{1,1}$ open set in $\R^d$ with characteristics $(R_0, \Lambda)$. Let $R_2$ be the constant in Lemma \ref{l:explictdecay}. For $t \in (0,T]$, we set 
$$
r_t=r_t(T, R_0, \Lambda, \psi, d) := \frac{L^{-1}(1/t)}{L^{-1}(1/T)}R_2.
$$

For $x \in D$ with $\delta_D(x)<r_t/2$, we define an open neighborhood $U(x,t)$ of $x$ and an open ball $W(x,t) \subset D \setminus U(x,t)$ as follows:

\smallskip

\textit {Find $z_x \in  \partial D$ satisfying $\delta_D(x) = |x- z_x|$ and let $v_x :=z_x+ 2r_t(x-z_x)/|x-z_x|$. Then, we have $\delta_D(v_x) \ge r_t/\sqrt{1+\Lambda^2}.$ We define 
}
\begin{align}\label{e:dUW}
U(x,t):= D \cap B(z_x, r_t) \quad \text{\textit{and}} \quad W(x,t):= B\big(v_x, \frac{r_t}{2\sqrt{1+\Lambda^2}}\big) \subset D.
\end{align}

Note that by the construction, we have that for all $u \in U(x,t)$ and $w \in W(x,t)$,
$$
|u-w| \ge |z_x-v_x| - |u-z_x| - |v_x -w| \ge 2r_t - r_t - r_t/2 \ge r_t/2
$$
and
$$
|u-w| \le |z_x-v_x| + |u-z_x| + |v_x -w| \le 4r_t.
$$
It follows that
\begin{align}\label{UWdist}
|u-w| \asymp r_t \quad \text{for all} \;\; u \in U(x,t), \; w \in W(x,t).
\end{align}

\begin{prop}\label{p:survival}
Let $D$ be a $C^{1,1}$ open set in $\R^d$ with characteristics $(R_0, \Lambda)$. Let $r_t$ and $U(x,t)$ be defined as in just before the Proposition. For all $T>0$ and $M \ge 1$, we have that for every $t \in (0, T]$ and $x \in D$ with $\delta_D(x)< r_t/2$,
\begin{align*}
\P_x(\tau_D>t) &\asymp \P_x(\tau_D>Mt) \asymp \P_x(Y_{\tau_{U(x,t)}} \in D)\asymp t^{-1} \E_x[\tau_{U(x,t)}] \asymp \big(tL(\delta_D(x))\big)^{-1/2},
\end{align*}
where the comparison constants depend only on $T, M, \psi, R_0, \Lambda$ and $d$.
\end{prop}
\pf
Recall that $z_x \in \partial D$ is the point satisfying $\delta_D(x) = |x-z_x|$. Let 
$$
o_x= z_x+ \frac{r_t(x-z_x)}{2|x-z_x|} \in D.
$$ 
Indeed, we have $r_t/(2\sqrt{1+\Lambda^2}) \le \delta_D(o_x) \le r_t/2.$ By {\bf (A)} and {\bf (B)}, we see that assumptions in \cite[Theorem 1.9]{GK18} hold and hence by that theorem, the (scale-invariant) boundary Harnack principle holds. Therefore, we get 
\begin{align}\label{BHP}
\P_x(Y_{\tau_{U(x,t)}} \in D) \le c\frac{\P_{x}(Y_{\tau_{U(x,t)}} \in W(x,t))}{\P_{o_x}(Y_{\tau_{U(x,t)}} \in W(x,t))}\P_{o_x}(Y_{\tau_{U(x,t)}} \in D)  \le c\frac{\P_{x}(Y_{\tau_{U(x,t)}} \in W(x,t))}{\P_{o_x}(Y_{\tau_{U(x,t)}} \in W(x,t))},
\end{align}
where $W(x,t)$ is the subset of $D$ defined as in just before the Proposition.
By the L\'evy system, \eqref{e:doubling}, \eqref{UWdist} and Lemma \ref{l:explictdecay}, we have
\begin{align*}
\P_{x}(Y_{\tau_{U(x,t)}} \in W(x,t)) &= \E_x\left[\int_0^{\tau_{U(x,t)}} \int_{W(x,t)} \nu(|Y_s-w|)dw ds\right] \\
&\asymp \E_x[\tau_{U(x,t)}] \nu(r_t)r_t^d \asymp L(r_t)^{-1/2}L(\delta_D(x))^{-1/2} \nu(r_t)r_t^d.
\end{align*}
Similarly, we also have that  $\P_{o_x}(Y_{\tau_{U(x,t)}} \in W(x,t)) \asymp \E_{o_x}[\tau_{U(x,t)}] \nu(r_t)r_t^d \asymp L(r_t)^{-1}\nu(r_t)r_t^d.$

Then, by the strong Markov property, Chebyshev's inequality, \eqref{BHP} and Lemma \ref{l:explictdecay}, since $L(r_t) \asymp t^{-1}$, we obtain
\begin{align*}
\P_x(\tau_D>t) & \le \P_x(\tau_{U(x,t)}>t) + \P_x(Y_{\tau_{U(x,t)}} \in D) \le t^{-1} \E_x[\tau_{U(x,t)}] + c L(r_t)^{1/2}L(\delta_D(x))^{-1/2}\\
& \le t^{-1}L(r_t)^{-1/2}L(\delta_D(x))^{-1/2}+c L(r_t)^{1/2}L(\delta_D(x))^{-1/2} \le ct^{-1/2}L(\delta_D(x))^{-1/2}.
\end{align*}

On the other hand, for any $a>0$, by the strong Markov property,  \eqref{exitball}, \eqref{e:asym2}, Lemma \ref{l:explictdecay} and Chebyshev's inequality,
\begin{align*}
&\P_x(\tau_D>at)\\
 &\ge \P_x\big(\tau_{U(x,t)}<at, Y_{\tau_{U(x,t)}} \in D_{int}(r_t/4), |Y_{\tau_{U(x,t)}}-Y_{\tau_{U(x,t)}+s}| \le r_t/4 \text{ for all } 0<s<at\big)\\
&\ge \P_x\big(\tau_{U(x,t)}<at, Y_{\tau_{U(x,t)}} \in D_{int}(r_t/4)\big) \P_0(\tau_{B(0, r_t/4)}>at)\\
&\ge c_1\big(\P_x(Y_{\tau_{U(x,t)}} \in D_{int}(r_t/4))-\P_x(\tau_{U(x,t)}\ge at)\big)\\
&\ge c_1\big(c_2t^{-1}\E_x[\tau_{U(x,t)}]-a^{-1}t^{-1}\E_x[\tau_{U(x,t)}]\big).
\end{align*}
Take $a=(2c_2^{-1}) \vee M$. By Lemma \ref{l:explictdecay} and the fourth line in the above inequalities, we obtain
\begin{align*}
\P_x(\tau_D>Mt)\ge \P_x(\tau_D>at) &\ge \frac{c_1}{2}\P_x(Y_{\tau_{U(x,t)}} \in D_{int}(r_t/4)) \ge c t^{-1/2}L(\delta_D(x))^{-1/2}.
\end{align*}
This completes the proof.
\qed

\begin{cor}\label{c:survival}
Let $D$ be a $C^{1,1}$ open set in $\R^d$ with characteristics $(R_0, \Lambda)$. For all $T>0$, there exists a constant $c_1=c_1(d, T, \psi, R_0, \Lambda)>1$ such that for every $t \in (0, T]$ and $x \in D$,
\begin{align*}
c_1^{-1} \left(1 \wedge \frac{1}{tL(\delta_D(x))}\right)^{1/2} \le \P_x(\tau_D>t) \le c_1 \left(1 \wedge \frac{1}{tL(\delta_D(x))}\right)^{1/2}.
\end{align*}
\end{cor}
\pf
We use the same notations as those in Proposition \ref{p:survival}. If $\delta_D(x)<r_t/2$, then the result follows from Proposition \ref{p:survival}. If $\delta_D(x) \ge r_t/2$, then $tL(\delta_D(x))\le tL(r_t/2)\le c$. Also, by \eqref{exitball} and \eqref{e:asym2}, we get $1 \ge \P_x(\tau_D>t) \ge \P_x(\tau_{B(x,r_t/2)}>t) \ge c$.
\qed

\begin{cor}\label{c:bddsurvival}
Let $D$ be a bounded $C^{1,1}$ open subset in $\R^d$ with characteristics $(R_0, \Lambda)$ of scale $(r_1, r_2)$. Then, there exists $c_1=c_1(R_0, \Lambda, \psi, d)>0$ such that for all $t>0$ and $x \in D$,
\begin{align*}
&c_1^{-1} \left(1 \wedge \frac{1}{(t \wedge 2)L(\delta_D(x))}\right)^{1/2} \exp\big(-\kappa_2 C_4 t h(r_1/2)\big) \\
&\qquad \qquad \le \P_x(\tau_D>t) \le c_1  \left(1 \wedge \frac{1}{(t \wedge 2)L(\delta_D(x))}\right)^{1/2} \exp\big(-\kappa_1 C_5 t h(r_2)\big),
\end{align*}
where $\kappa_1, \kappa_2$ are constants in {\bf (A)} and $C_4, C_5$ are constants in \eqref{exitball}.
\end{cor}
\pf
Fix $(t,x) \in (0, \infty) \times D$. If $t \le 2$, then the assertion follows from Corollary \ref{c:survival}. Hence, we assume that $t >2$. Let $x_1, x_2 \in \R^d$ be the points satisfying $B(x_1, r_1) \subset D \subset B(x_2, r_2)$. By the semigroup property, \eqref{exitball} and Corollary \ref{c:survival}, we get
\begin{align*}
\P_x( \tau_D>t) &= \int_D p_D(t,x,y) dy \le \int_D \int_D p_D(1,x,z) p_{B(x_2, r_2)}(t-1,z,y) dz dy \\
& \le  \P_x(\tau_D>1)\sup_{z \in D} \P_z(\tau_{B(x_2, r_2)}>t-1) \le \frac{c }{L(\delta_D(x))^{1/2}} \exp\big(-\kappa_1C_5t h(r_2)\big).
\end{align*}

To prove the lower bound, we first assume that $\delta_D(x)<R_2/2$ where $R_2$ is the constant in Lemma \ref{l:explictdecay}. Without loss of generality, we may assume that $R_2 \le r_1/2$. Let $z \in \partial D$ be the point satisfying $\delta_D(x)=|x-z|$ and $\theta$ be the shift operator defined as $Y_t \circ \theta_s= Y_{s+t}$. Then, by the strong Markov property, \eqref{e:exitdist}, the L\'evy system and \eqref{exitball}, we have
\begin{align*}
\P_x(\tau_D>t) & \ge \E_x\Big[Y_{\tau_{D \cap B(z, R_2)}} \in D_{int}(R_2/4),  \; Y_{\tau_{B(Y_0,R_2/4)}} \circ \theta_{\tau_{D \cap B(z, R_2)}} \in B(x_1, r_1/2), \\
&\qquad\quad \tau_D\circ \theta_{\tau_{B(Y_0,R_2/4)}} \circ \theta_{\tau_{D \cap B(z, R_2)}}  >t \Big] \\
 \ge&\; \frac{cL(R_2)^{1/2}}{L(\delta_D(x))^{1/2}} \inf_{w \in D_{int}(R_2/4)} \P_w\big(Y_{\tau_{B(w,R_2/4)}} \in B(x_1, r_1/2)\big)\inf_{y \in B(x_1, r_1/2)} \P_y(\tau_{B(x_1, r_1)}>t) \\
\ge& \;  \frac{c}{L(\delta_D(x))^{1/2}}\exp\big(-\kappa_2 C_4t h(r_1/2)\big).
\end{align*}
Indeed, on $\{Y_{\tau_{D \cap B(z, R_2)}} \in D_{int}(R_2/4)\}$, we have $B(Y_{\tau_{D \cap B(z, R_2)}}, R_2/4) \subset D$. Also, since $R_2 \le r_1/2$, we can always find $A \subset B(x_1, r_1/2) \setminus B(Y_{\tau_{D \cap B(z, R_2)}}, R_2/2)$ such that $|A| \ge c_1>0$ for some constant $c_1>0$. Then, by the L\'evy system and \eqref{meanexit}, we obtain
\begin{align*}
 \P_x\big(Y_{\tau_{B(Y_0,R_2/4)}} \circ \theta_{\tau_{D \cap B(z, R_2)}} \in B(x_1, r_1/2)\big)  \ge \E_0 \left[\int_0^{\tau_{B(0,R_2/4)}} \int_{A} \nu(|Y_s- y|) dy ds\right] \ge c>0.
\end{align*} 
Similarly, if $\delta_D(x) \ge R_2/2$, then we have
\begin{align*}
\P_x(\tau_D>t) & \ge \E_x[Y_{\tau_{B(x,R_2/4)}} \in B(x_1, r_1/2), \tau_D\circ \theta_{\tau_{B(x,R_2/4)}}>t] \\
& \ge c \inf_{y \in B(x_1, r_1/2)} \P_y(\tau_{B(x_1, r_1)}>t) \ge c \exp\big(-\kappa_2 C_4t h(r_1/2)\big).
\end{align*}
\qed

\section{Small time Dirichlet heat kernel estimates in $C^{1,1}$ open set}\label{s:DHKE}

In this section, we provide the proof of Theorem \ref{t:main}.
Let $T>0$ be a fixed constant and $D$ be a fixed $C^{1,1}$ open set in $\R^d$ with characteristics $(R_0, \Lambda)$. We assume that {\bf (B)} holds. If $D$ is unbounded, then we further assume that {\bf (C)} holds. Then, by \eqref{e:doubling} and \eqref{e:doubling2}, we have
\begin{align}\label{doubling}
\nu(|x-y|) \asymp \nu(2|x-y|) \quad \text{for all} \;\; x,y \in D.
\end{align}
By \eqref{smallKnu}, \eqref{largeKnu}, Corollary \ref{c:HKEU} and Corollary \ref{c:HKEB}, we have the following heat kernel estimates for small $t$. Let $b_0$ be the constant in Proposition \ref{p:glb}.

\smallskip

\setlength{\leftskip}{4mm}
\noindent (1) If {\bf (S-1)} holds, then there exist constants $c_1>1$ and $b_2>0$ such that
\begin{align}\label{e:HKEB}
c_1^{-1} t \nu(|x|)\exp\big(-b_0th(|x|)\big) \le p(t,x) \le c_1t \nu(|x|)\exp\big(-b_2th(|x|)\big),
\end{align}
for all $(t,x) \in (0,T] \times (D\setminus \{0\})$.

\smallskip

\noindent (2) If {\bf (S-2)} holds, then there exist a constant $c_2>1$ such that 
\begin{align}\label{e:HKEU}
&c_2^{-1} t\nu(\theta_{\eta}(|x|,t))\exp\big(-b_0t h(\theta_{\eta}(|x|,t))\big) \nn\\
&\quad \qquad\qquad \le p(t,x) \le c_2t\nu(\theta_{a_1}(|x|,t))\exp\big(-b_1th(\theta_{a_1}(|x|,t))\big),
\end{align}
for all $(t,x) \in (0,T] \times D$ and $\eta>0$ where $a_1$ and $b_1$ are the constants in Lemma \ref{l:ondiag1}, and $\theta_{a}(r,t)=r \vee [\ell^{-1}(a/t)]^{-1}$ is defined as \eqref{e:deftheta}.

\setlength{\leftskip}{0mm}

\medskip

Before proving  Theorem \ref{t:main}, we obtain a lower bound of $p_D(t,x,y)$ without   {\bf (S-1)} and {\bf (S-2)}. This result will be used later to obtain Green function estimates.

\begin{prop}\label{p:pDl}
For every $T>0$, there exist positive constants $c_1=c_1(d, \psi, T, R_0, \Lambda)$ and $c_2=c_2(d, \psi, T, R_0, \Lambda)$ such that for all $(t,x,y) \in (0,T] \times (D \times D \setminus \diag)$,
\begin{align*}
p_D(t,x,y) \ge c_2 \left(1\wedge \frac{1}{tL(\delta_D(x))} \right)^{1/2}\left(1\wedge \frac{1}{tL(\delta_D(y))} \right)^{1/2} t \nu(|x-y|)\exp\big(-c_1th(|x-y|)\big).
\end{align*}
\end{prop}
\pf
Let $R_2$ be the constant in Lemma \ref{l:explictdecay}. Fix $(t,x,y) \in (0,T] \times (D \times D \setminus \diag)$ and set 
\begin{align}
\label{e:rtlt}
r_t = \frac{L^{-1}(1/t)}{L^{-1}(1/T)}R_2 \qquad \text{and} \qquad l_t(x,y) = r_t \wedge \frac{|x-y|}{4}.
\end{align}
Note that by \eqref{e:asym2}, \eqref{e:vws} and \eqref{e:Lws}, we have $V(r_t) \asymp t^{1/2}$ and $L(r_t) \asymp h(r_t) \asymp t^{-1}.$
 
 Let $z_x, z_y \in \partial D$ be the points satisfying $\delta_D(x) = |x-z_x|$ and $\delta_D(y) = |y-z_y|$.  By \eqref{e:vws}, there exists a constant $m>1$ such that 
 \begin{equation}\label{e:scaleV}
 mV(\delta k) \ge \delta V(k) \quad \text{for all} \;\; 0<\delta \le 1, \; 0<k \le 1.
 \end{equation} 

\smallskip

\textit{Case 1.} Suppose that $|x-y| \le R_2$. Define open neighborhoods of $x$ and $y$ as follows:
\begin{center}
$\sO(x)=
\begin{cases}
B\big(x,V^{-1}[\frac{1}{8m}V(|x-y|)]\big), & \mbox{if } \;\; 8mV(\delta_D(x)) \ge V(|x-y|); \\
D \cap B(z_x,\frac{1}{3}|x-y|), & \mbox{if }  \;\; 8mV(\delta_D(x)) < V(|x-y|),
\end{cases}
$
\end{center}
and
\begin{center}
$\sO(y)=
\begin{cases}
B\big(y,V^{-1}[\frac{1}{8m}V(|x-y|)]\big), & \mbox{if } \;\; 8mV(\delta_D(y)) \ge V(|x-y|); \\
D \cap B\big(z_y,\frac{1}{3}|x-y|\big), & \mbox{if } \;\; 8mV(\delta_D(y)) < V(|x-y|).
\end{cases}
$
\end{center}
Then, we see that $x \in \sO(x) \subset D, \, y \in \sO(y) \subset D$ and 
\begin{align*}
|u-w| \asymp |x-y| \qquad \text{for all} \;\; u \in \sO(x), \; w \in \sO(y).
\end{align*}
Thus, by the strong Markov property and \eqref{doubling}, we have (cf. ~ \cite[Lemma 3.3]{CKS12b},)
\begin{align}\label{e:general}
p_D(t,x,y) &\ge t \P_x(\tau_{\sO(x)}>t) \P_y(\tau_{\sO(y)}>t) \inf_{u \in \sO(x), w \in \sO(y)}\nu(|u-w|) \nn\\
&\ge ct\nu(|x-y|) \P_x(\tau_{\sO(x)}>t) \P_y(\tau_{\sO(y)}>t).
\end{align}

To calculate the survival probability $\P_x(\tau_{\sO(x)}>t)$, we first assume that $8mV(\delta_D(x)) \ge V(|x-y|)$. In this case, we see from \eqref{exitball} and \eqref{e:asym2} that
\begin{align}\label{intlow}
\P_x(\tau_{\sO(x)}>t) &\ge c \exp\big(-c_1tV(|x-y|)^{-2}\big) \ge c \exp\big(-c_2 th(|x-y|)\big).
\end{align}

Next, assume that $8mV(\delta_D(x)) < V(|x-y|)$. Note that by the monotonicity of $V$ and \eqref{e:scaleV}, we get $|x-y| > 8 \delta_D(x)$ in this case. Let $\rho:=V^{-1}(\varepsilon V(l_t(x,y)))$ where $\varepsilon \in (0,(8m)^{-1})$ will be chosen later. Then, we see from \eqref{e:asym2} and \eqref{e:vws} that
\begin{equation}\label{e:rhol}
V(\rho) \asymp V(l_t(x,y)) \asymp t^{1/2} \wedge V(|x-y|) \quad \text{and} \quad h(\rho) \asymp h(l_t(x,y)) \asymp t^{-1} \vee h(|x-y|).
\end{equation}
Note that we can not expect that $\rho \asymp l_t(x,y)$ in general.

If $8\delta_D(x)\ge \rho$, then by \eqref{exitball} and \eqref{e:rhol}, we have
\begin{align}\label{sO}
\P_x(\tau_{\sO(x)}>t) &\ge \P_x(\tau_{B(x,\rho/8)}>t)\ge c\exp\big(-\kappa_2 C_4th(\rho/8)\big) \ge c\exp\big(-c_3th(|x-y|)\big).
\end{align}
Indeed, by Lemma \ref{l:propPhi}(i) and \eqref{e:asym2}, we see that $h(\rho/8) \asymp h(4\rho)$. Thus, if $l_t(x,y) = |x-y|/4$, then we get \eqref{sO}. Otherwise, if $l_t(x,y) = r_t$, then $\P_x(\tau_{\sO(x)}>t) \asymp 1 \asymp \exp\big(-c_3th(|x-y|)\big)$ and hence \eqref{sO} holds.

If $8\delta_D(x)<\rho$, then we can find a piece of annulus $\sA(x) \subset \{w \in \sO(x) : \rho<|w-z_x|<|x-y|/4\}$ such that $\dist(\sA(x), \partial \sO(x))>\rho/8$. Recall that $\theta$ is shift operator. Then, by the strong Markov property, the L\'evy system,  \eqref{exitball}, \eqref{e:survivaltime},  \eqref{e:asym2} and \eqref{e:vws}, we have
\begin{align*}
\P_x(\tau_{\sO(x)}>t) &\ge \P_x\big(Y_{\tau_{B(z_x,\rho/2) \cap D}} \in \sA(x), \; \tau_{\sO(x)}\circ \theta_{\tau_{B(z_x,\rho/2) \cap D}} >t \big)\\[4pt]
&\ge \P_x\big(Y_{\tau_{B(z_x,\rho/2) \cap D}} \in \sA(x)\big) \inf_{z \in \sA(x)} \P_z(\tau_{\sO(x)}>t)\\
&\ge c \E_x\left[\int_0^{\tau_{B(z_x,\rho/2) \cap D}}\int_{\sA(x)} \nu(|Y_s-w|) dw ds \right]\P_0(\tau_{B(0,\rho/8)}>t)\\
&\ge c \E_x\big[\tau_{B(z_x,\rho/2) \cap D}\big] \int_{\rho}^{|x-y|/4} (-L'(k)) dk \exp\big(-c_4t h(|x-y|)\big) \\
&\ge c  \big(L(\rho)-L(|x-y|/4)\big)L(\delta_D(x))^{-1/2}L(\rho/2)^{-1/2}\exp\big(-c_4t h(|x-y|)\big)\\[4pt]
&\ge c\big(c_5^{-1}V(\rho)^{-2}-c_5V(|x-y|)^{-2}\big)  L(\delta_D(x))^{-1/2} V(\rho) \exp\big(-c_4t h(|x-y|)\big),
\end{align*}
where $c_5>1$ is a constant independent of choice of $\varepsilon$. Now, we choose $\varepsilon=(2c_5)^{-1} \wedge (16m)^{-1}$. Then, we get from \eqref{e:rhol} that
\begin{align*}
\P_x(\tau_{\sO(x)}>t) &\ge  c V(\rho)^{-1}L(\delta_D(x))^{-1/2}\exp\big(-c_4t h(|x-y|)\big) \\
&\ge c t^{-1/2}L(\delta_D(x))^{-1/2}\exp\big(-c_4t h(|x-y|)\big).
\end{align*}

Finally, by combining the above inequality with \eqref{intlow} and \eqref{sO}, we deduce that 
$$
\P_x(\tau_{\sO(x)}>t) \ge c \left(1 \wedge \frac{1}{tL(\delta_D(x))}\right)^{1/2}\exp\big(-c_5t h(|x-y|)\big).
$$
By the same way, we get $\P_y(\tau_{\sO(y)}>t) \ge c \big(1 \wedge \frac{1}{tL(\delta_D(y))}\big)^{1/2}\exp\big(-c_5t h(|x-y|)\big)$. Then, \eqref{e:general} yields the desired lower bound.

\smallskip

\textit{Case 2.} Suppose that $|x-y|>R_2$. In this case, we let $D_x:=D \cap B(x,R_2/4)$ and $D_y:=D \cap B(y, R_2/4)$. By the same argument as \eqref{e:general}, \eqref{doubling} and Corollary \ref{c:survival}, we get
\begin{align*}
p_D(t,x,y) &\ge t\P_x(\tau_{D_x}>t) \P_y(\tau_{D_y}>t) \inf_{u \in D_x, w \in D_y} \nu(|u-w|) \\
&\ge c\left(1 \wedge \frac{1}{tL(\delta_D(x))}\right)^{1/2}\left(1 \wedge \frac{1}{tL(\delta_D(y))}\right)^{1/2}t\nu(|x-y|).
\end{align*}
This completes the proof.  \qed

\vspace{2mm}

Now, we are ready to prove Theorem \ref{t:main}.

\vspace{2mm}
\noindent{\bf Proof of  Theorem \ref{t:main}.} 
 Fix $(t,x,y) \in (0,T] \times (D \times D \setminus \diag)$ and continue using the notation $r_t$ and $l_t(x,y)$ in \eqref{e:rtlt}.
\smallskip

\noindent (i) Since we have proved the lower bound in Proposition \ref{p:pDl}, it suffices to show that there exist $c_1>0, b_3 \in (0,b_0]$ such that for all $(t,x,y) \in (0,T] \times (D \times D \setminus \diag)$,
\begin{align}\label{upper0}
p_D(t,x,y) \le c_1\left(1 \wedge \frac{1}{tL(\delta_D(x))}\right)^{1/2} t\nu(|x-y|)\exp\big(-b_3th(|x-y|)\big),
\end{align}
where $b_0$ is the constant in Proposition \ref{p:glb}.
Indeed, if \eqref{upper0} holds, then by the semigroup property and \eqref{e:HKEB}, we get
\begin{align*}
p_D(t,x,y)&=\int_D p_D(t/2,x,z)p_D(t/2,y,z) dz\\
& \le c\left(1\wedge \frac{1}{tL(\delta_D(x))} \right)^{1/2}\left(1\wedge \frac{1}{tL(\delta_D(y))} \right)^{1/2} \int_D p\big(\frac{b_3}{2b_0}t,x,z\big)p\big(\frac{b_3}{2b_0}t,y,z\big) dz \\
& \le c\left(1\wedge \frac{1}{tL(\delta_D(x))} \right)^{1/2}\left(1\wedge \frac{1}{tL(\delta_D(y))} \right)^{1/2} t\nu(|x-y|)\exp\big(-\frac{b_2b_3}{b_0}th(|x-y|)\big),
\end{align*}
which yields the result.

\smallskip

Now, we prove \eqref{upper0}. If $\delta_D(x) \ge r_t/2$, then \eqref{upper0} is a consequence of \eqref{e:HKEB} and the trivial bound that $p_D(t,x,y) \le p(t, x-y)$. Hence, we assume that $\delta_D(x) < r_t/2$. By \eqref{e:Lws}, there exists a constant $M>1$ such that 
\begin{equation}\label{e:Lscale}
ML(16k) \ge L(k) \quad \text{for all} \;\;  k \le 1/16.
\end{equation}
Observe that by the semigroup property, monotonicity of $p(t, \cdot)$ and Proposition \ref{p:survival}, we have
\begin{align*}
p_D(t,x,y)&=\left( \int_{\{z\in D : |y-z|>|x-y|/2 \}}+\int_{\{z\in D : |y-z| \le |x-y|/2 \}}\right) p_D(t/2,x,z)p_D(t/2,z,y) dz\\
&\le\left( \int_{\{z\in D : |y-z|>|x-y|/2 \}}+\int_{\{z\in D : |x-z| > |x-y|/2 \}}\right) p_D(t/2,x,z)p_D(t/2,z,y) dz\\[3pt]
&\le p(t/2,|x-y|/2)\big(\P_x(\tau_D>t/2)+\P_y(\tau_D>t/2)\big)\\[3pt]
&\le cp(t/2,|x-y|/2)\big(t^{-1/2}L(\delta_D(x))^{-1/2}+t^{-1/2}L(\delta_D(y))^{-1/2}\big).
\end{align*}
Thus, if $ML(\delta_D(y)) \ge L(\delta_D(x))$, then we get \eqref{upper0}. Therefore, we assume that $ML(\delta_D(y))<L(\delta_D(x))$. Since $L$ is strictly decreasing, it follows from \eqref{e:Lscale} that $\delta_D(y)>16\delta_D(x)$ and hence $|x-y| \ge |y-z_x|-|z_x-x|\ge \delta_D(y)-\delta_D(x)>15\delta_D(x)$ where $z_x \in \partial D$ is the point satisfying $\delta_D(x)=|x-z_x|$. Define
$$
W_1:=D \cap B(z_x, l_t(x,y)), \quad W_3:=\big\{w\in D : |w-y| \le |x-y|/2 \big\}
$$ and $W_2:=D \setminus (W_1 \cup W_3)=\{w\in D \setminus W_1 : |w-y| > |x-y|/2 \}$. Then, for $u \in W_1$ and $w \in W_3$, we obtain
\begin{align}\label{dist13}
|u-w| \ge |x-y|-|z_x-x|-|u-z_x|-|y-w|\ge \big(1-\frac{1}{15}-\frac{1}{4}-\frac{1}{2}\big)|x-y|> \frac{|x-y|}{6}.
\end{align}

Observe that by the strong Markov property,
\begin{align}\label{upper_decom}
p_D(t,x,y) &= \E_x\big[\,p_D(t-\tau_{W_1}, Y_{\tau_{W_1}},y) : \tau_{W_1}<t\big]\nn\\
&=\E_x[\;p_D(t-\tau_{W_1}, Y_{\tau_{W_1}},y) : \tau_{W_1}<t, Y_{\tau_{W_1}} \in W_3]\nn\\
&\quad+\E_x[\,p_D(t-\tau_{W_1}, Y_{\tau_{W_1}},y) : \tau_{W_1}\in (0,2t/3], Y_{\tau_{W_1}} \in W_2]\nn\\
&\quad+\E_x[\,p_D(t-\tau_{W_1}, Y_{\tau_{W_1}},y) : \tau_{W_1}\in (2t/3,t), Y_{\tau_{W_1}} \in W_2]\nn\\
&=: I_1+I_2+I_3.
\end{align}

First, by the L\'evy system and \eqref{dist13}, we get
\begin{align}\label{i1}
I_1 &=\int_0^t \int_{W_3} \int_{W_1} p_{W_1}(s,x,u)\nu(|w-u|)p_D(t-s,w,y) du dw ds \nn\\
&\le \nu(|x-y|/6) \int_0^t \P_x(\tau_{W_1}>s) \int_{W_3}p(t-s,y-w)dw ds.
\end{align}
By \eqref{e:HKEB} and Lemma \ref{l:asymLh}, for all $s \in (0,T]$ and $l \in (0, 2r_t]$, we have
\begin{align}\label{e:intball}
\int_{B(y,l)}p(s,y-w) dw &\le c\int_0^l -s L'(k) \exp\big(-c_2s L(k)\big)dk \le c\exp\big(-c_3 sh(l)\big).
\end{align}
It follows that for all $s \in (0, t]$, 
\begin{align}\label{e:W3}
\int_{W_3} p(s, y-w)dw & \le c
\begin{cases}
\exp\big(-c_3 sh(|x-y|)\big), & \mbox{if } \;\; |x-y| \le 2r_t ; \\
1, & \mbox{if } \;\; |x-y| > 2r_t
\end{cases}
\nn\\
&\le c\exp\big(-c_3 sh(|x-y|)\big).
\end{align}
Indeed, if $|x-y| >2r_t$, then we have $sh(|x-y|) \le sh(2r_t) \asymp st^{-1} \le 1$.
Moreover, by the semigroup property, Proposition \ref{p:survival}, \eqref{e:intball} and monotonicity of $h$, we get
\begin{align}\label{e:W1}
&\P_x(\tau_{W_1}>2t/3)=\int_{W_1} \int_{W_1} p_{W_1}(t/3,x,v) p_{W_1}(t/3,v,u) dvdu \nn \\
&\le \P_x(\tau_D>t/3) \int_{B(0,2l_t(x,y))} p(t/3,u)du\le ct^{-1/2}L(\delta_D(x))^{-1/2} \exp\big(-3^{-1}c_3 th(2l_t(x,y)) \big) \nn \\
& \le c t^{-1/2}L(\delta_D(x))^{-1/2} \exp\big(-3^{-1}c_3 th(|x-y|) \big).
\end{align}
Then, using \eqref{i1}, \eqref{doubling}, \eqref{e:intball}, \eqref{e:W1} and Proposition \ref{p:survival}, we obtain
\begin{align}\label{e:i1}
&I_1 \le c \nu(|x-y|) \int_0^{t} \P_x(\tau_{W_1}>s) \int_{W_3}p(t-s,y-w)dw ds \nn\\
&\quad \le c \nu(|x-y|)\exp\big(-c_3t h(|x-y|)/3\big)\int_0^{2t/3}\P_x(\tau_D>s)ds \nn\\
&\quad \quad +c \nu(|x-y|)\P_x(\tau_{W_1}>2t/3)\int_0^{t/3} \exp\big(-c_3 s h(|x-y|)\big)ds \nn\\
&\quad \le cL(\delta_D(x))^{-1/2} \nu(|x-y|)\exp\big(-c_3t h(|x-y|)/3\big) \Big(\int_0^{2t/3} s^{-1/2}ds+t^{-1/2} \int_0^{t/3}  ds\Big) \nn\\[3pt]
&\quad =  ct^{-1/2}L(\delta_D(x))^{-1/2} t\nu(|x-y|)\exp\big(-3^{-1}c_3t h(|x-y|)\big).
\end{align}

Second, by monotonicity of $p(t, \cdot)$, \eqref{e:HKEB}, \eqref{doubling} and Proposition \ref{p:survival}, we get
\begin{align}\label{e:i2}
I_2 & \le c\P_x(Y_{\tau_{W_1}} \in W_2)  \sup_{s \in [t/3, t), l \ge |x-y|/2} p(s, l) = c \P_x(Y_{\tau_{W_1}} \in W_2)  \sup_{s \in [t/3, t)} p(s, |x-y|/2) \nn\\[3pt]
&\le c\P_x(Y_{\tau_{W_1}} \in W_2)  \nu(|x-y|) \Big(\sup_{s \in [t/3,t)} s\exp\big(-b_2 s h(|x-y|)\big) \Big) \nn\\
& \le c
\begin{cases}
L(r_t)^{1/2}L(\delta_D(x))^{-1/2} t \nu(|x-y|) \exp\big( -3^{-1}b_2th(|x-y|) \big), & \mbox{if } \;\; |x-y| \ge 4r_t ; \\[2pt]
L(|x-y|)^{1/2}L(\delta_D(x))^{-1/2} t \nu(|x-y|) \exp\big( -3^{-1}b_2th(|x-y|) \big), & \mbox{if } \;\; |x-y| < 4r_t
\end{cases}
\nn\\[4pt]
& \le ct^{-1/2}L(\delta_D(x))^{-1/2} t \nu(|x-y|) \exp\big( -4^{-1}b_2th(|x-y|) \big).
\end{align}
In the last inequality, we used \eqref{e:asym2}, $L(r_t) \asymp t^{-1}$ and the fact thats $e^x \ge 2e\sqrt{x}$ for all $x>0$ and $h(r) \ge L(r)$ for all $r>0$.

Lastly, we note that $t \mapsto te^{-at}$ is increasing on $(0, 1/a)$ and decreasing on $(1/a, \infty)$. Thus, using similar calculation as the one given in \eqref{e:W1}, by monotonicity of $p(t, \cdot)$, \eqref{e:HKEB}, \eqref{doubling}, Proposition \ref{p:survival} and \eqref{e:asym2}, we have
\begin{align*}
I_3&\le c\P_x(\tau_{W_1}>2t/3) \nu(|x-y|) \Big(\sup_{s \in (0,t/3)} s\exp\big( -b_2sh(|x-y|) \big) \Big) \\
& \le c
\begin{cases}
\P_x(\tau_{W_1}>2t/3)\nu(|x-y|)h(|x-y|)^{-1}  , & \mbox{if }\;\; b_2th(|x-y|) \ge 3 ; \\[2pt]
\P_x(\tau_D>2t/3)\nu(|x-y|) t\exp\big( -3^{-1}b_2th(|x-y|) \big) , & \mbox{if }\;\; b_2th(|x-y|) < 3
\end{cases}
\\
& \le c
\begin{cases}
t^{-1/2}L(\delta_D(x))^{-1/2} t \nu(|x-y|)\exp\big(-3^{-1}c_3 th(|x-y|) \big), & \mbox{if }\;\; b_2th(|x-y|) \ge 3 ; \\[2pt]
t^{-1/2}L(\delta_D(x))^{-1/2} t \nu(|x-y|)\exp\big(-2^{-1}b_2 th(|x-y|) \big), & \mbox{if }\;\; b_2th(|x-y|) < 3.
\end{cases}
\end{align*}
Combining the above inequality with \eqref{e:i1}, \eqref{e:i2} and \eqref{upper_decom}, we get \eqref{upper0}.

\vspace{4mm}

\noindent (ii) We use the same notations as in the proof of (i) and follow the proof of (i).

\smallskip

 {\bf(Upper bound)} By the semigroup property and \eqref{e:HKEU}, it suffices to show that there exist positive constants $c_1$ and $b_4$ such that 
\begin{align}\label{upper1}
p_D(t,x,y) \le c_1\Big(1\wedge \frac{1}{tL(\delta_D(x))} \Big)^{1/2} t \nu(\theta_{3a_1}(|x-y|, t)) \exp\big(-b_4th(\theta_{3a_1}(|x-y|, t))\big).
\end{align}
By the similar argument to the one given in the proof of (i), we may assume that $\delta_D(x)<r_t/2$. Moreover, observe that for every $u, v \in D$, by the triangle inequality, $\max\{|x-u|, |u-v|, |v-y|\} \ge |x-y|/3$. Thus, by the semigroup property, monotonicity of $p(t, \cdot)$ and Proposition \ref{p:survival}, we have that
\begin{align*}
p_D(t,x,y) & =\int_{u \in D, |x-u| \ge |x-y|/3} \int_{D}  p_D(t/3, x, u) p_D(t/3, u, v) p_D(t/3, v, y) dvdu \\
& \quad +  \int_{u \in D, |x-u| < |x-y|/3} \int_{v \in D, |u-v| \ge |x-y|/3 } p_D(t/3, x, u) p_D(t/3, u, v) p_D(t/3, v, y) dudv  \\
& \quad +  \int_{u \in D, |x-u| < |x-y|/3} \int_{v \in D, |u-v| < |x-y|/3 } p_D(t/3, x, u) p_D(t/3, u, v) p_D(t/3, v, y) dudv  \\
& \le p(t/3, |x-y|/3) \int_D p_D(t/3, v, y) \int_D p_D(t/3, u, v) du  dv \\
& \quad + p(t/3, |x-y|/3) \int_D p_D(t/3, x, u) du \int_D p_D(t/3, v, y) dv \\
& \quad + p(t/3, |x-y|/3) \int_D p_D(t/3, x, u) \int_D p_D(t/3, u, v) dv du \\
& \le 2p(t/3, |x-y|/3) \big(\P_x(\tau_D>t/3) + \P_y(\tau_D>t/3)\big)\\[3pt]
& \le cp(t/3, |x-y|/3) \big(t^{-1/2} L(\delta_D(x))^{-1/2} + t^{-1/2} L(\delta_D(y))^{-1/2}\big).
\end{align*}
Therefore, if $ML(\delta_D(y)) \le L(\delta_D(x))$ for some $M>0$, then we get \eqref{upper1} from \eqref{e:HKEU}. Hence, we may assume that $\delta_D(y)>16 \delta_D(x)$ by the same argument as the one given in the proof of (i).

To prove \eqref{upper1},  we first assume that $|x-y| \le [\ell^{-1}(3a_1/t)]^{-1}$. In this case, we have that $\theta_{a_1}(|x-y|, t/3) = \theta_{3a_1}(|x-y|, t) = [\ell^{-1}(3a_1/t)]^{-1}$. Then, by the semigroup property, \eqref{e:HKEU} and Proposition \ref{p:survival}, we get
\begin{align*}
p_D(t,x,y) &= \int_D p_D(2t/3,x,z)p_D(t/3,z,y)dz \le c \P_x(\tau_D>2t/3) p(t/3,0) \\
&\le c t^{-1/2}L(\delta_D(x))^{-1/2} t \nu(\theta_{3a_1}(|x-y|, t)) \exp\big(-3^{-1}b_1th(\theta_{3a_1}(|x-y|, t))\big).
\end{align*}

Now, suppose that $|x-y| >[\ell^{-1}(3a_1/t)]^{-1}$. In this case, we use \eqref{upper_decom} and find upper bounds for $I_1$, $I_2$ and $I_3$. Observe that for all $s \in (0,T]$ and $l \in (0,2r_t]$, by \eqref{e:HKEU} and the similar calculation to the one given in \eqref{e:intball},
\begin{align}\label{e:intball2}
\int_{B(y, l)} p(s, y-w) dw 
&\le c
\begin{cases}
l^d[\ell^{-1}(a_1/s)]^{d}  \exp\big( -b_1sh([\ell^{-1}(a_1/s)]^{-1}) \big), & \mbox{if } \;\; l\le [\ell^{-1}(a_1/s)]^{-1} ; \\[2pt]
\exp\big( -c_2sh(l) \big), & \mbox{if } \;\; l> [\ell^{-1}(a_1/s)]^{-1}
\end{cases}
\nn\\
&\le c \exp\big(-c_3 s h(\theta_{a_1}(l,s))\big).
\end{align}
Then, by using \eqref{e:intball2} instead of \eqref{e:intball}, we have that for all $0<s \le T$,
\begin{align*}
\P_x(\tau_{W_1}>s) & = \int_{W_1}\int_{W_1} p_{W_1}(s/3,x,u)p_{W_1}(2s/3,u,v)dudv \\
& \le cs^{-1/2}L(\delta_D(x))^{-1/2}\exp\big(-c_4s h(\theta_{a_1}(|x-y|,2s/3))\big).
\end{align*}
Hence, by the similar arguments to the ones given in \eqref{e:W3} and \eqref{e:i1}, we obtain
\begin{align*}
 I_1 \le c t^{-1/2}L(\delta_D(x))^{-1/2} t \nu(|x-y|) \exp\big(-c_5 t h(|x-y|)\big).
\end{align*}

Next, by \eqref{e:HKEU}, \eqref{doubling}, monotonicity of $h$, we have
\begin{align*}
\sup_{s \in [t/3, t)} p(s, |x-y|/2) \le ct \sup_{s \in [t/3, t)} \Big[\nu(\theta_{a_1}(|x-y|, s)) \exp\big(-3^{-1}b_1 t h(\theta_{a_1}(|x-y|,s))\big)\Big].
\end{align*}
Let $f(r):=r^{-d}\exp\big(-c_7th(r)\big)$ where the constant $c_7 \in (0, b_1/3)$ will be chosen later. Then, by \eqref{smallKnu}, there exists a constant $c_6>0$ such that for $r \in (0, [\ell^{-1}(a_1/t)]^{-1})$,
\begin{align*}
r^{d+1}\exp\big(c_7th(r)\big)f'(r) =-d + 2c_7tK(r) \le -d + c_6c_7t \ell(r^{-1}). 
\end{align*}
Set $c_7=d/(3a_1c_6) \wedge b_1/3$. Then, we see that $f$ is decreasing on $\big([\ell^{-1}(3a_1/t)]^{-1}, [\ell^{-1}(a_1/t)]^{-1} \big)$. Using this fact, since $\ell$ is almost increasing, we deduce that
\begin{align*}
&\sup_{s \in [t/3, t)} \Big[\nu(\theta_{a_1}(|x-y|, s)) \exp\big(-3^{-1}b_1 t h(\theta_{a_1}(|x-y|,s))\big)\Big] \\
&\quad \le c\nu(\theta_{a_1}(|x-y|, t/3)) \exp\big(-c_7 t h(\theta_{a_1}(|x-y|,t/3))\big)= c\nu(|x-y|) \exp\big(-c_7 t h(|x-y|)\big).
\end{align*}
It follows that by the same argument as in the one given in \eqref{e:i2}, 
\begin{align*}
I_2 &\le ct^{-1/2}L(\delta_D(x))^{-1/2} t \nu(|x-y|) \exp \big(-c_8th(|x-y|)\big).
\end{align*}

Lastly, we note that since $|x-y| > [\ell^{-1}(3a_1/t)]^{-1}$,
\begin{align*}
&\sup_{s \in (0,t/3)} \Big[s \nu(\theta_{a_1}(|x-y|,s)) \exp\big(-b_1 s h(\theta_{a_1}(|x-y|,s))\big)\Big]\\
&= \sup_{s \in (0,t/3)} \Big[s \nu(|x-y|) \exp\big(-b_1 s h(|x-y|)\big)\Big].
\end{align*}
Therefore, by the same proof as in the one given in (i), we obtain 
$$
I_3 \le ct^{-1/2}L(\delta_D(x))^{-1/2} t \nu(|x-y|) \exp \big(-c_9th(|x-y|)\big).
$$
This finishes the proof for the upper bound.
\smallskip

 {\bf (Lower bound)}
Fix $\eta>0$. By Proposition \ref{p:pDl}, it remains to prove the lower bound when $|x-y|<[\ell^{-1}(\eta/t)]^{-1} \wedge R_2$, where $R_2$ is the constant in Lemma \ref{l:explictdecay}. Let $\zeta_t:=[\ell^{-1}(\eta/t)]^{-1} \wedge R_2$ and define open neighborhoods of $x$ and $y$ as follows. Recall that $z_x, z_y \in \partial D$ are the points satisfying $\delta_D(x) = |x-z_x|$ and $\delta_D(y) = |y-z_y|$. We define
\begin{center}
$\sU(x)=
\begin{cases}
B\big(x,V^{-1}(\frac{1}{8m}V(\zeta_t))\big), & \mbox{if} \;\; 8mV(\delta_D(x)) \ge V(\zeta_t); \\[2pt]
B(z_x,\frac{1}{3}\zeta_t) \cap D, & \mbox{if} \;\; 8mV(\delta_D(x)) < V(\zeta_t),
\end{cases}
$
\end{center}
and
\begin{center}
$\sU(y)=
\begin{cases}
B\big(y,V^{-1}(\frac{1}{8m}V(\zeta_t))\big), & \mbox{if } \;\; 8m V(\delta_D(y)) \ge V(\zeta_t); \\[2pt]
B(z_y,\frac{1}{3}\zeta_t) \cap D, & \mbox{if } \;\; 8mV(\delta_D(y)) < V(\zeta_t),
\end{cases}
$
\end{center}
where $m$ is the constants in \eqref{e:scaleV}.
Then, we can see that $x \in \sU(x) \subset D$ and $y \in \sU(y) \subset D$.

\smallskip

We claim that there exist a constant $c_4>0$ independent of the choice of $\eta$, and a constant $c_3>0$ such that
\begin{align}\label{cl:bdd}
&\P_x(\tau_{\sU(x)}>t) \ge c_3 \Big(1 \wedge \frac{1}{tL(\delta_D(x))}\Big)^{1/2} \exp \big(-c_4 t h(\zeta_t)\big).
\end{align}

Indeed, if $8m V(\delta_D(x)) \ge V(\zeta_t)$, then by \eqref{exitball} and \eqref{e:asym2}, we have
\begin{align}\label{case1}
\P_x(\tau_{\sU(x)}>t) \ge c \exp\big(-c_1 t h(\zeta_t)\big).
\end{align}

Suppose that $8m V(\delta_D(x)) < V(\zeta_t)$. If $\zeta_t = R_2$, then by Corollary \ref{c:survival}, we get
\begin{align}\label{case2}
\P_x(\tau_{\sU(x)}>t)  \ge c \Big(1 \wedge \frac{1}{tL(\delta_D(x))}\Big)^{1/2}.
\end{align}
Otherwise, if $\zeta_t=[\ell^{-1}(\eta/t)]^{-1} < R_2$, then by the similar argument to the one given in the proof of  Proposition \ref{p:pDl},
\begin{align}\label{case3}
\P_x(\tau_{\sU(x)}>t) &\ge cL(\zeta_t)^{1/2}L(\delta_D(x))^{-1/2} \exp\big(-c_2th(\zeta_t)\big) \nn\\
&\ge ct^{-1/2}L(\delta_D(x))^{-1/2} \exp\big(-c_2th(\zeta_t)\big).
\end{align}
In the second inequality, we used  $L(\zeta_t) \ge cK(\zeta_t)\asymp\ell(\zeta_t^{-1}) \asymp t^{-1}$ which follows from the proof of Lemma \ref{l:asymLh} and \eqref{smallKnu}.  By combining \eqref{case1}, \eqref{case2} and \eqref{case3}, we obtain \eqref{cl:bdd}.

 Let $w_x:=z_x+  4\zeta_t (x-z_x)/|x-z_x| \in D$ and define
$$
\sW_{0}:=B\big(w_x,\frac{\zeta_t}{2\sqrt{1+\Lambda^2}}\big) \quad \text{and} \quad  \sW:= B\big(w_x,\frac{\zeta_t}{\sqrt{1+\Lambda^2}}\big) \subset D.
$$
Then, for all $u \in \sU(x)$ and $v \in \sW$, we have $|u-v| \asymp \zeta_t$. Moreover, since $|x-y| < \zeta_t$, we also have $|u'-v| \asymp \zeta_t$ for all $u' \in \sU(y)$ and $v \in \sW$.
Thus, for every $v \in \sW_{0}$, by the similar argument to \eqref{e:general}, \eqref{exitball} and \eqref{cl:bdd}, we have
\begin{align*}
p_D(t/2,x,v) &\ge ct \nu(\zeta_t) \P_x(\tau_{\sU(x)}>t/2) \P_v\big(\tau_{B(v,\zeta_t/(2\sqrt{1+\Lambda^2}))}>t/2\big) \\
& \ge c  \Big(1 \wedge \frac{1}{tL(\delta_D(x))}\Big)^{1/2} t \nu(\zeta_t) \exp\big(-c_5 t h(\zeta_t)\big).
\end{align*}
Similarly, we also have that
$$
p_D(t/2,v,y) \ge c   \Big(1 \wedge \frac{1}{tL(\delta_D(y))}\Big)^{1/2} t \nu(\zeta_t) \exp\big(-c_5 t h(\zeta_t)\big).
$$
It follows that by the semigroup property and {\bf (A)},
\begin{align*}
&p_D(t,x,y) \ge \int_{\sW}p_D(t/2,x,v)p_D(t/2,v,y) dv \\
&\quad \ge c \Big(1 \wedge \frac{1}{tL(\delta_D(x))}\Big)^{1/2} \Big(1 \wedge \frac{1}{tL(\delta_D(y))}\Big)^{1/2} t^2 |\sW| \, \nu(\zeta_t)^2\exp\big(-2c_5 t h(\zeta_t)\big) \\
& \quad \ge c\Big(1 \wedge \frac{1}{tL(\delta_D(x))}\Big)^{1/2} \Big(1 \wedge \frac{1}{tL(\delta_D(y))}\Big)^{1/2} t^2 \ell(\zeta_t^{-1}) \nu(\zeta_t) \exp\big(-2c_5 t h(\zeta_t)\big).
\end{align*}
If $\zeta_t = [\ell^{-1}(\eta/t)]^{-1}$, then since $\ell$ is almost increasing, we get $\ell(\zeta_t^{-1}) \asymp t^{-1}$. Hence, we are done. Otherwise, if $\zeta_t = R_2$, then we have $t \asymp 1$ and hence $ t^2 \ell(\zeta_t^{-1}) \nu(\zeta_t) \exp\big(-2c_5 t h(\zeta_t)\big) \asymp t \nu([\ell^{-1}(\eta/t)]^{-1}) \exp\big(-cth([\ell^{-1}(\eta/t)]^{-1})\big) \asymp 1$.
This completes the proof. \qed

\section{Large time estimates}\label{s:DHKE2}

In this section, we give the proof of Theorem \ref{t:main2}. Let $D$ be a fixed bounded $C^{1,1}$ open subset in $\R^d$ of scale $(r_1, r_2)$ and $x_1, x_2 \in \R^d$ be the fixed points satisfying $B(x_1, r_1) \subset D \subset B(x_2, r_2)$. We mention that under condition {\bf (L-1)}, the transition semigroup $\{P_t^D, t \ge 0\}$ of $Y_t^D$ may not be compact operators in $L^2(D;dx)$, though $D$ is bounded. (See, Proposition \ref{p:existence}.) Hence, in that case, we need some lemmas to obtain the large time estimates instead of the general spectral theory.

\begin{lemma}\label{l:large1}
There exists a constant $c_1>0$ which only depend on the dimension $d$ such that for all $(t,x,y) \in (0, \infty) \times(D \times D \setminus \diag)$,
\begin{equation*}
 p_D(t,x,y) \le c_1p(t/2, |x-y|/2) \exp\big(-2^{-1}\kappa_1C_5th(r_2)\big).
\end{equation*}
\end{lemma}
\pf
By the semigroup property, we have
\begin{align*}
p_D(t,x,y)&=\left( \int_{\{z\in D : |y-z|>|x-y|/2 \}}+\int_{\{z\in D : |y-z| \le |x-y|/2 \}}\right) p_D(t/2,x,z)p_D(t/2,z,y) dz\\
&\le\left( \int_{\{z\in D : |y-z|>|x-y|/2 \}}+\int_{\{z\in D : |x-z| > |x-y|/2 \}}\right) p_D(t/2,x,z)p_D(t/2,z,y) dz\\[3pt]
&\le p(t/2,|x-y|/2)\big(\P_x(\tau_{B(x_2, r_2)}>t/2)+\P_y(\tau_{B(x_2, r_2)}>t/2)\big).
\end{align*}
Hence, we obtain the result from \eqref{exitball}.
\qed

Define for $r \ge 1$,
\begin{align*}
\tll(r):=\sup_{s \in [1, r]} \frac{1}{\ell(s)} \qquad \text{and} \qquad \tPhi(r):= \int_1^r \frac{1}{k \tll(k)} dk.
\end{align*}
Note that if {\bf (L-1)} holds, we have that
\begin{align}\label{tll}
\tll(r)^{-1} \asymp \ell(r) \qquad \text{for all} \;\; r \ge 2.
\end{align}
Moreover, by the same argument as the one given in the proof of Lemma \ref{l:propPhi}, there exist positive constants $C_6$ and $C_7$ which only depend on the dimension $d$ and $\kappa_1$ and $\kappa_2$ in \eqref{A}  such that
\begin{align}\label{e:asym3}
C_6\tPhi(r) \le \psi(r) \quad \text{and} \quad h(r^{-1}) \le C_7\tPhi(r) \quad\;\; \text{for all} \;\; r\ge 2.
\end{align}
We also note that $\tPhi$ satisfies $\WS^{\infty}(0, 2, 1)$.
Here, we get the large time on-diagonal estimates for $p(t,x)$ under condition {\bf (L-1)}.

\begin{lemma}\label{l:ondiag3}
Assume that {\bf (L-1)}  holds. Then, there exists a positive constant $b_5=b_5(d,\psi, r_2)$ such that for any $T>0$, there exist  $c_1, c_2>0$ such that 
for all $t \in [T, \infty)$ and $|x| \le 2r_2$,
\begin{equation}\label{e:l6.2}
p(t,x) \le c_1 + c_2\nu(|x|)\exp\big(-b_5t h(|x|)\big).
\end{equation}
\end{lemma}
\pf 
Fix $x \in \R^d$ satisfying $|x| \le 2r_2$ and let $r=|x|$. By \cite[(5.4)]{GRT19}, the mean value theorem, Lemma \ref{l:simple}, \eqref{tll} and \eqref{e:asym3},   we have that, for all $t>0$ (cf. \eqref{upperbase}),
\begin{align}\label{fourier}
r^dp(t, x) &\le c \int_{\R^d} \left(e^{-t \psi(|z|/r)} - e^{-t \psi(2|z|/r)} \right) e^{-|z|^2/4} dz \nn \\
&\le cr^d + ct \int_{|z|>2r} \sup_{|z| \le y \le 2|z|}e^{-t \psi(y/r)} \big|\psi(2|z|/r)-\psi(|z|/r)\big| e^{-|z|^2/4} dz \nn \\
& \le cr^d + ct \int_{2r}^{4r_2} e^{-C_6t \tPhi(u/r)} \frac{u^{d-1}}{\tll(u/r)}du +  ct \int_{4r_2}^{\infty} e^{-C_6t \tPhi(u/r)} \frac{u^{d-1}}{\tll(u/r)}e^{-u^2/4}du \nn \\
& =: cr^d + I_1+I_2.
\end{align}

First, by \eqref{tll}, and the monotonicity and the scaling properties of $\tll$, $\ell$ and $\tPhi$, we have
\begin{align}\label{e:l6.2-2}
I_2 &\le ct\big[\tll(4r_2/r)\big]^{-1} \exp\big(-C_6 t \tPhi(4r_2/r)\big) \int_{4r_2}^{\infty} u^{d-1}e^{-u^2/4} du  \nn\\
&\le ct \ell(4r_2/r) \exp\big(-c_1 t \tPhi(1/r)\big) \le c\ell(1/r)\exp\big(-2^{-1}c_1 t \tPhi(1/r)\big).
\end{align}
In the last inequality above, we used the facts that $e^x\ge x$ for $x>0$ and $\tPhi(1/r) \ge \tPhi(1/(2r_2))$.
 
On the other hand, set $q_{\gamma, k}(u)= q_{\gamma, k}(u, t):=  u^{\gamma}\exp\big(-kt\tPhi(u)\big)$ for $u \ge 2$ and $\gamma, k>0$. We observe that for any $\gamma,k>0$,
\begin{align*}
\frac{d}{du}q_{\gamma,k}(u) =  \big(\gamma-kt \tll(u)^{-1}\big)q_{\gamma-1,k}(u).
\end{align*}
Since $\tll$ is increasing, it follows that there exists $u_0 \in [2, \infty)$   such that $q$ is decreasing on $[2, u_0]$ and increasing on $[u_0, \infty)$. Thus, for any $[a,b] \subset [2,\infty)$ and $\gamma,k>0$, we have that 
\begin{equation}\label{e:qmax}
\sup_{u \in [a,b]}q_{\gamma,k}(u) = q_{\gamma,k}(a) \vee q_{\gamma,k}(b).
\end{equation}
Choose any constant $\eps\in (0,1)$ such that $(1-\eps)d+\alpha_1>0$.  This is possible since $\alpha_1>-d$. Then we set $\rho:=((1-\eps)d+\alpha_1)/(1+\eps) \in (0, d+\alpha_1)$ so that
\begin{equation}\label{e:defrho}
\eps = (d+\alpha_1-\rho)/(d+\rho).
\end{equation}

 First, suppose that $q_{d+\rho,\,C_6}(2) \ge q_{d+\rho,\, C_6}(4r_2/r)$. Then by the change of variables and \eqref{e:qmax}, since $\tll$ is increasing,  we have
\begin{align*}
I_1&=c tr^d \int_{2}^{4r_2/r} \frac{u^{d}}{u\tll(u)} e^{-C_6t \tPhi(u)} du = c tr^d  \int_{2}^{4r_2/r} \frac{q_{d+\rho, \, C_6}(u)}{u^{1+\rho}\tll(u)} du\nn\\
& \le c tr^d \frac{q_{d+\rho,\, C_6}(2)}{\tll(2)} \int_2^{4r_2/r} \frac{du}{u^{1+\rho}} \le ctr^d \exp \big(- C_6 t \tPhi(2)\big) \le \frac{c}{C_6\tPhi(2)}r^d.
\end{align*}
Hence, we obtain \eqref{e:l6.2} from \eqref{fourier}, \eqref{e:l6.2-2}, \eqref{A} and  \eqref{e:asym3}  in this case.

Now, we suppose that $q_{d+\rho,\,C_6}(2) < q_{d+\rho,\, C_6}(4r_2/r)$. Note that by \eqref{tll} and the scaling property of $\ell$, it holds that for any $2<u \le s$, 
\begin{align*}
\frac{\tll(s)}{\tll(u)} \le c \frac{\ell(u)}{\ell(s)} \le c \bigg(\frac{u}{s}\bigg)^{\alpha_1}, \qquad \text{i.e.,} \quad \frac{1}{u^{\alpha_1}\tll(u)} \le \frac{c}{s^{\alpha_1}\tll(s)}.
\end{align*}
We also note that for any $\gamma, k>0$ and $u \ge 2$, the map $q_{\gamma,k}(u) \le q_{\gamma, \eps k}(u)$ and $q_{\gamma,  k}(u)^s = q_{s\gamma, \, sk}(u)$ for all $s>0$. Thus, by  the change of variables,  \eqref{tll}, \eqref{e:qmax}, \eqref{e:defrho} and the scaling properties of $\ell$ and $\tPhi$, since we have assumed $q_{d+\rho,\,C_6}(2) < q_{d+\rho,\, C_6}(4r_2/r)$, we get that
\begin{align*}
I_1&=c tr^d  \int_{2}^{4r_2/r} \frac{q_{d+\alpha_1-\rho,\, C_6}(u)}{u^{1-\rho}u^{\alpha_1}\tll(u)} du \le\frac{c tr^d}{(4r_2/r)^{\alpha_1} \tll(4r_2/r)} \int_{2}^{4r_2/r} \frac{q_{d+\alpha_1-\rho,\, \eps C_6}(u)}{u^{1-\rho}} du \\
& \le ctr^{d+\alpha_1} \ell(1/r) \big( q_{d+\alpha_1-\rho,\, \eps C_6}(2) \vee q_{d+\alpha_1-\rho,\, \eps C_6}(4r_2/r) \big) \int_2^{4r_2/r} \frac{du}{u^{1-\rho}}\\[3pt]
& \le ctr^{d+\alpha_1-\rho} \ell(1/r) \big( q_{d+\rho,\, C_6}(2) \vee q_{d+\rho,\, C_6}(4r_2/r) \big)^{(d+\alpha_1-\rho)/(d+\rho)} \\[3pt]
& = ctr^{d+\alpha_1-\rho} \ell(1/r)  q_{d+\alpha_1-\rho,\, \eps C_6}(4r_2/r)= ct  \ell(1/r) \exp \big( -\eps C_6 t \tPhi(4r_2/r)\big)\\[3pt]
&\le c \ell(1/r) \exp \big( -2^{-1}\eps C_6 t \tPhi(4r_2/r)\big) \le c \ell(1/r) \exp \big( -c_2 t \tPhi(1/r)\big).
\end{align*}
Then we get \eqref{e:l6.2} by using \eqref{fourier}, \eqref{e:l6.2-2}, \eqref{A} and  \eqref{e:asym3} again. \qed

Now, we give the proof of Theorem \ref{t:main2}.

\smallskip

\noindent{\bf Proof of  Theorem \ref{t:main2}.} 
Let $a(x,y):=L(\delta_D(x))^{-1/2}L(\delta_D(y))^{-1/2}$.

\smallskip

\noindent (i) Choose $(t,x,y) \in [T, \infty) \times (D \times D \setminus \diag)$ and let $x_1 \in D$ be the point satisfying $B(x_1, r_1) \subset D$.  By the semigroup property, Theorem \ref{t:main}(i), \eqref{e:HKEB} and \eqref{exitball}, we have 
\begin{align}\label{e:large1}
& p_D(t,x,y) \ge \int_{B(x_1, r_1/4) \times B(x_1, 3r_1/4)} p_D(T/4, x, u) p_D(t-T/2, u, v) p_D(T/4, v, y) dudv \nn\\
& \ge c a(x,y) \int_{B(x_1, r_1/4) \times B(x_1, 3r_1/4)} p(cT, 2r_2)^2 p_D(t-T/2, u, v) dudv \nn\\
&  \ge c a(x,y) \inf_{u \in B(x_1, r_1/4)} \P_u(\tau_{B(x_1,3r_1/4)}>t-T/2) \ge c a(x,y) \exp\big(-\kappa_2 C_4 t h(r_1/2)\big).
\end{align}

On the other hand, since $D$ is a bounded set, one can follow the proof of Proposition \ref{p:pDl}, after changing the definition of $l_t(x,y)$ therein from $r_t \wedge (|x-y|/4)$ to $|x-y|/4$, and see that
\begin{align}\label{e:large2}
p_D(t,x,y) &\ge c \left(1 \wedge \frac{L(|x-y|)}{L(\delta_D(x))}\right)^{1/2}\left(1 \wedge \frac{L(|x-y|)}{L(\delta_D(x))}\right)^{1/2} t\nu(|x-y|) \exp \big(-c_1 t h(|x-y|)\big)\nn\\
& \ge cT\left(1 \wedge \frac{L(2r_2)}{L(\delta_D(x))}\right)^{1/2}\left(1 \wedge \frac{L(2r_2)}{L(\delta_D(x))}\right)^{1/2} \nu(|x-y|) \exp \big(-c_1 t h(|x-y|)\big) \nn\\[3pt]
&\ge ca(x,y)\nu(|x-y|)\exp \big(-c_1 t h(|x-y|)\big).
\end{align}
In the last inequality above, we used the fact that $L(2r_2)/L(\delta_D(z)) \le L(2r_2)/L(r_2) \le c$ for all $z \in D$, which comes from the monotonicity of $L$ and \eqref{e:Lws}. By combining \eqref{e:large1} with \eqref{e:large2},   we get the desired lower bound.  

Now, we prove the upper bound. By the semigroup property, Theorem \ref{t:main}(i), Corollary \ref{c:HKEB}, Lemma \ref{l:large1} and Lemma \ref{l:ondiag3}, we get
\begin{align*}
p_D(t,x,y) &= \int_{D \times D} p_D(T/4, x, u) p_D(t-T/2, u, v) p_D(T/4, v, y) dudv \\
& \le c a(x,y) \exp\big(-2^{-1}\kappa_1 C_5 th(r_2)\big)  \\
& \quad \times\int_{D \times D} p(cT/4, |x-u|/2) p(t/2-T/4, |u-v|/2) p(cT/4, |v-y|/2) dudv \\
& \le c a(x,y) p(t/2-cT, |x-y|/2) \exp\big(-2^{-1}\kappa_1 C_5 th(r_2)\big) \\
& \le c a(x,y) \Big[c+c\nu(|x-y|)\exp \big(-2^{-1}b_5 t h(|x-y|)\big) \Big]\exp\big(-2^{-1}\kappa_1 C_5 th(r_2)\big),
\end{align*}
which yields the upper bound.

\smallskip

\noindent (ii) \& (iii) Since the proof of (iii) is similar and easier, we only provide the proof of (ii). By Proposition \ref{p:existence}, there exist $T_0>0$ such that the transition semigroup $\{P_t^D : t \ge T_0 \}$ of $Y_t^D$ consists of compact operators. Let $0<\mu_1<1$ be the largest eigenvalue of the operator $P_{T_0}^D$ and $\phi_1 \in L^2(D;dx)$ be the corresponding eigenfunction with unit $L^2$-norm. For $n \ge 1$, we denote by $\{\mu_{n,k}; k \ge 1\} \subset (0,1)$ the discrete spectrum of $P_{nT_0}^D$, arranged in decreasing order and repeated according to their multiplicity and $\{\phi_{n,k}; k \ge 1\}$ be the corresponding eigenfunctions with unit $L^2$-norm. Then, by the semigroup property, we have $\mu_{n,1} = \mu_1^n$ and $\phi_{n,1} = \phi_1$ for all $n \ge 1$. From the eigenfunction expansion of $p_D(nT_0, x, \cdot)$ and Parseval's identity, we have for  $n \ge 1$,
\begin{align}\label{spectral1}
\int_{D \times D} p_D(nT_0, x, y) dx dy & = \sum_{k=1}^{\infty} \mu_{n,k} \left(\int_D \phi_{n,k}(y) dy \right)^2 \le \sup_{k} \mu_{n,k} \int_D 1^2 dy = \mu_1^n |D|.
\end{align}
On the other hand, for all $s>0$ and $x \in D$, since $p(T_0, 0) \le c_0< \infty$, we have
\begin{align*}
\phi_1(x) &\le \int_{D \times D} p_D(s, x, z) p_D(T_0, z, y) \phi_1(y) dz dy \le c_0 \P_x(\tau_D>s) \int_D \phi_1(y) dy \\
& \le c_0 \P_x(\tau_D>s) \lVert \phi_1 \rVert_{L^2(D)}\left(\int_D 1^2 dy\right)^{1/2} = c_0 |D|^{1/2} \P_x(\tau_D>s).
\end{align*}
Thus, we obtain for all $0< s \le T_0$ and $n \ge 1$,
\begin{align}\label{spectral2}
&\int_{D \times D} \P_x(\tau_D>s)p_D(nT_0, x, y)\P_y(\tau_D>s) dx dy \nn\\
&\quad \ge \mu_1^n \left(\int_D \P_y(\tau_D>s) \phi_1(y) dy \right)^2  \ge \mu_1^n \left( \int_D c_0^{-1}|D|^{-1/2} \phi_1(y)^2 dy \right)^2 \ge c_0^{-2} |D|^{-1}\mu_1^n.
\end{align}

For $t \ge 4T_0$ and $x,y \in D$, we let $n:= \lfloor (t-3T_0)/T_0 \rfloor \ge 1$ and $s:= t - (n+2)T_0 \in [T_0, 2T_0)$. Recall  $a(x,y)=L(\delta_D(x))^{-1/2}L(\delta_D(y))^{-1/2}$.   By \eqref{spectral1} and Corollary \ref{c:survival}, we have 
\begin{align*}
p_D(t,x,y) &= \int_{D \times D \times D \times D} p_D(s/2, x, z_1) p_D(T_0, z_1, z_2) p_D(nT_0, z_2, z_3)\\
& \qquad \qquad \qquad  \times p_D(T_0, z_3, z_4) p_D(s/2, z_4, y) dz_1 dz_2 dz_3 dz_4\\
&\le c_0^2 \int_D p_D(s/2, x, z_1) dz_1 \int_{D \times D} p_D(nT_0, z_2, z_3) dz_2dz_3 \int_D p_D(s/2,z_4,y)dz_4 \\[3pt]
&\le c_0^2 |D| \P_x(\tau_D>s/2) \P_y(\tau_D>s/2)  \mu_1^n  \le ca(x,y) e^{-\lambda_1 t},
\end{align*}
where $\lambda_1 := T_0^{-1} \log (\mu_1^{-1})$. Moreover, by Theorem \ref{t:main}, Corollary \ref{c:survival} and \eqref{spectral2}, we get
\begin{align*}
& p_D(t,x,y)  = \int_{D \times D} p_D(s/2, x, z_1) p_D((n+2)T_0, z_1, z_2) p_D(s/2, z_2, y) dz_1dz_2 \\
&  \ge ca(x,y) \int_{D \times D} \P_{z_1}(\tau_D>s/2) p_D((n+2)T_0, z_1, z_2) \P_{z_2}(\tau_D>s/2) dz_1dz_2  \ge ca(x,y) e^{-\lambda_1 t}.
\end{align*}
This completes the proof.
\qed

\section{Green function estimates}

In this section, we provide the proof of Theorem \ref{t:green}. Throughout this section, we assume that {\bf (D)} holds.
\begin{lemma}\label{l:green1}
For all Borel set $D$ and $x,y \in D$, we have
\begin{align}\label{e:neG1}
\left(1 \wedge \frac{V(\delta_D(x))}{V(|x-y|)}\right)\left(1 \wedge \frac{V(\delta_D(y))}{V(|x-y|)}\right) \asymp \left(1 \wedge \frac{V(\delta_D(x))V(\delta_D(y))}{V(|x-y|)^2}\right).
\end{align}
In particular, if $D$ is bounded then for all $x,y \in D$
\begin{align*}
\left(1 \wedge \frac{L(|x-y|)}{L(\delta_D(x))}\right)^{1/2}\left(1 \wedge \frac{L(|x-y|)}{L(\delta_D(y))}\right)^{1/2} \asymp \bigg(1 \wedge \frac{L(|x-y|)}{\sqrt{L(\delta_D(x))L(\delta_D(y))}}\bigg).
\end{align*}
\end{lemma}
\pf
Since $(1 \wedge p)(1 \wedge q) \le 1 \wedge (pq)$ for every $p, q > 0$, one inequality  in \eqref{e:neG1} is trivial. On the other hand, since $1 \wedge (p/q) \asymp p/(p+q)$ for every $p, q>0$, it suffices to prove that
\begin{align*}
\big(V(\delta_D(x))+V(|x-y|)\big)\big(V(\delta_D(y))+V(|x-y|)\big) \le 2V(\delta_D(x))V(\delta_D(y))+4V(|x-y|)^2.
\end{align*}
By symmetry, we may assume that $\delta_D(x) \le \delta_D(y)$. According to the subadditivity of $V$,
\begin{align*}
&\big(V(\delta_D(x))+V(|x-y|)\big)\big(V(\delta_D(y))+V(|x-y|)\big) \\
& \quad \le V(\delta_D(x))V(\delta_D(y)) + V(|x-y|)^2 + 2V(|x-y|)V(\delta_D(y)) \\
& \quad \le V(\delta_D(x))V(\delta_D(y)) + V(|x-y|)^2 + 2V(|x-y|)(V(\delta_D(x))+V(|x-y|)) \\
& \quad \le V(\delta_D(x))V(\delta_D(y)) + V(|x-y|)^2 + V(|x-y|)^2 +  V(\delta_D(x))^2 + 2V(|x-y|)^2 \\
& \quad \le 2V(\delta_D(x))V(\delta_D(y))+4V(|x-y|)^2.
\end{align*}
This proves \eqref{e:neG1}. The second claim follows from \eqref{e:asym2}.
\qed

\begin{lemma}\label{l:green2}
It holds that
\begin{align*}
\liminf_{r \to 0} \frac{\nu(r)}{L(r)} = \liminf_{r \to 0} \frac{\nu(r)}{L(r)^2} = \infty.
\end{align*} 
\end{lemma}
\pf
Since the L\'evy measure $\nu$ is infinite, we have $\lim_{r \to 0} L(r) = \infty$. Thus, it suffices to show that the second equality holds. By l'Hospital's rule, \cite[Lemmas 3.1 and 3.2]{CKKW18}, \eqref{A} and \eqref{e:Lws},  since $\alpha_2<d$,  we have
\begin{align*}
\liminf_{r \to 0} \frac{\nu(r)}{L(r)^2} \ge c\liminf_{r \to 0} \frac{r^{-1}\nu(r)}{2r^{-1}L(r)\ell(r^{-1})} \ge \frac{c}{L(1)} \liminf_{r \to 0} \frac{L(1)r^{-d}}{L(r)} \ge c \liminf_{r \to 0}  r^{-d+(\alpha_2 \vee 2^{-1})} = \infty.
\end{align*}
Indeed, since  $r \mapsto \nu(r^{-1})$ satisfies $\WS^{\infty}(d+\alpha_1,  d+\alpha_2 , 1)$ and $d+\alpha_1>0$, according to \cite[Lemma 3.1 and 3.2]{CKKW18}, there exists a function $\wt{\nu}(r)$ such that for all $0<r<1$, $\wt{\nu}(r) \asymp \nu(r)$ and $-\wt{\nu}'(r) \asymp r^{-1}\wt{\nu}(r) \asymp r^{-1}\nu(r)$. Hence, the first inequality in the display  holds.
\qed

Recall that for a Borel subset $D \subset \R^d$, the Green function $G_D(x,y)$ is defined by
\begin{align*}
G_D(x,y):= \int_0^{\infty} p_D(t,x,y)dt.
\end{align*}
Since the  process $Y$ can be recurrent, we can not expect to obtain upper estimates for $G_{\R^d}(x,y)$ in general. However, when $D$ is bounded, we can establish a prior upper estimates for $G_D(x,y)$ regardless of transience of $Y$. By $\diam(D)$ we denote the diameter of $D$.

\begin{lemma}\label{intgreen}
Let $D \subset \R^d$ be a bounded Borel set. Then, there exists a constant $c_1=c_1(d, \psi, \diam(D))>0$ such that for all $x,y \in D$,
\begin{align*}
G_D(x,y) \le \frac{c_1\ell(|x-y|^{-1})}{|x-y|^dL(|x-y|)^2} \asymp \frac{\nu(|x-y|)}{L(|x-y|)^2}.
\end{align*}
\end{lemma}
\pf
Fix $x,y \in D$ and let $r=|x-y|$. If $x=y$, by Lemma \ref{l:green2}, there is nothing to prove. Hence, we assume that $r>0$.

By Lemma \ref{l:large1}, \eqref{upperbase} and Fubini's Theorem, we have
\begin{align*}
&r^d G_D(x,y) \le c\int_0^{\infty} r^d p(t/2, r/2) \exp \big(- 2^{-1}\kappa_1 C_5th(\diam(D))\big) dt  \\
&\le c \int_0^{\infty} ctr^d\exp \big(- 2^{-1}\kappa_1 C_5th(\diam(D))\big) dt + c \int_0^{\infty}t\int_r^{\infty} e^{-C_0t\Phi(u/r)} \ell(u/r)u^{d-1}e^{-u^2/4}dudt \\
&\le cr^d + c\int_r^{\infty} \ell(u/r)u^{d-1}e^{-u^2/4} \int_0^{\infty} t e^{-C_0t\Phi(u/r)}dtdu \\
&\le cr^d + c \int_r^{1} \frac{\ell(u/r)u^{d-1}}{C_0^2 \Phi(u/r)^2}du + c \int_1^{\infty} \frac{\ell(u/r)u^{d-1}e^{-u^2/4}}{C_0^2 \Phi(u/r)^2}du\\
&=: cr^d + I_1 + I_2.
\end{align*}

First, by the change of the variables, we have
\begin{align*}
I_1 = c r^d\int_1^{1/r} \frac{\ell(u)u^{d-1}}{\Phi(u)^2}du = cr^d \int_1^{1/r} \frac{\ell(u)u^{d-1}}{L(u^{-1})^2} du.
\end{align*}
Observe that by applying l'Hospital's rule three times and the same argument as the one, coming from \cite[Lemmas 3.1 and 3.2]{CKKW18},  given in the proof of Lemma \ref{l:green2}, we obtain
\begin{align*}
&\limsup_{s \to 0} \frac{s^d\int_1^{1/s} \ell(u) u^{d-1}L(u^{-1})^{-2}du}{\ell(s^{-1})L(s)^{-2}} = \limsup_{s \to 0} \frac{L(s)^2\int_1^{1/s} \ell(u) u^{d-1}L(u^{-1})^{-2}du}{s^{-d}\ell(s^{-1})} \\
&\quad \le c+ c\limsup_{s \to 0} \frac{s^{-1}\ell(s^{-1})L(s)\int_1^{1/s} \ell(u)u^{d-1}L(u^{-1})^{-2} du}{s^{-d-1}\ell(s^{-1})}\\
&\quad = c+c\limsup_{s \to 0} \frac{L(s)\int_1^{1/s} \ell(u)u^{d-1}L(u^{-1})^{-2}du}{s^{-d}}\\
&\quad \le c+c\limsup_{s \to 0} \frac{\ell(s^{-1})L(s)^{-1}s^{-d-1}}{s^{-d-1}}+c\limsup_{s \to 0} \frac{s^{-1}\ell(s^{-1})\int_1^{1/s} \ell(u)u^{d-1}L(u^{-1})^{-2} du}{s^{-d-1}} \\
&\quad \le c+c\limsup_{s \to 0} \frac{\int_1^{1/s} \ell(u)u^{d-1}L(u^{-1})^{-2} du}{s^{-d}/\ell(s^{-1})} \\
&\quad \le c+c\limsup_{s \to 0} \frac{\ell(s^{-1})L(s)^{-2}s^{-d-1}}{s^{-d-1}/\ell(s^{-1})} = c+c\limsup_{s \to 0} \frac{\ell(s^{-1})^2}{L(s)^2} \le c.
\end{align*}
The fourth inequality above is valid, since  we can assume that $-(r^{-d}/\ell(r^{-1}))' \asymp r^{-d-1}/\ell(r^{-1})$ for $0<r<1$ by the argument given in the proof of Lemma \ref{l:green2} because  $r \mapsto r^d/\ell(r)$ satisfies $\WS^{\infty}(d-\alpha_2,  d-\alpha_1,  1)$ and $\alpha_2<d$.  In the third and fifth inequalities above, we used the fact that $\ell(r^{-1}) \le cL(r)$ for $0<r<1$, which follows from \eqref{smallKnu} and the proof of Lemma \ref{l:asymLh}.
Thus, since $D$ is bounded, we get that $I_1 \le c\ell(r^{-1})L(r)^{-2}$.

On the other hand, by the scaling property of $\ell$ and the monotonicity of $\Phi$, we obtain
\begin{align*}
I_2 \le c \frac{\ell(r^{-1})}{\Phi(r^{-1})^2} \int_1^{\infty} \frac{\ell(u/r)}{\ell(1/r)}u^{d-1} e^{-u^2/4}du \le   c \frac{\ell(r^{-1})}{\Phi(r^{-1})^2} \int_1^{\infty}  u^{d-1+\alpha_2} e^{-u^2/4}du =  c \frac{\ell(r^{-1})}{L(r)^2}.
\end{align*}
This completes the proof. \qed

\noindent{\bf Proof of  Theorem \ref{t:green}.}
Fix $x, y \in D$ and set $a(x,y):=L(\delta_D(x))^{-1/2}L(\delta_D(y))^{-1/2}$. By \eqref{A} and Lemma \ref{l:asymLh}, it suffices to prove that
$$
G_D(x,y) = \int_0^{\infty} p_D(t,x,y) dt \asymp \Big(1 \wedge \big[a(x,y)L(|x-y|)\big] \Big) \frac{\nu(|x-y|)}{h(|x-y|)^2}.
$$

{\bf(Lower  bound)} 
By Proposition \ref{p:pDl}, we have that
\begin{align}\label{green_lower}
&G_D(x,y) \ge \int_0^1 p_D(t, x, y) dt \nn\\
& \quad  \ge c\nu(|x-y|) \int_0^1 \left(1 \wedge \frac{1}{tL(\delta_D(x))} \right)^{1/2} \left(1 \wedge \frac{1}{tL(\delta_D(y))} \right)^{1/2}  t\exp\big(-c_1th(|x-y|)\big)dt\nn\\
& \quad = c \frac{\nu(|x-y|)}{h(|x-y|)^2} \int_0^{h(|x-y|)} \left(1 \wedge \frac{h(|x-y|)}{sL(\delta_D(x))} \right)^{1/2} \left(1 \wedge \frac{h(|x-y|)}{sL(\delta_D(y))} \right)^{1/2} se^{-c_1 s}ds\nn\\
& \quad \ge c \frac{\nu(|x-y|)}{h(|x-y|)^2} \left(1 \wedge \frac{L(|x-y|)}{L(\delta_D(x))} \right)^{1/2} \left(1 \wedge \frac{L(|x-y|)}{L(\delta_D(y))} \right)^{1/2}\int_0^{h(2r_2)\wedge 1}  se^{-c_1 s}ds\nn\\
& \quad \ge c\bigg(1 \wedge \big[a(x,y)L(|x-y|)\big] \bigg) \frac{\nu(|x-y|)}{h(|x-y|)^2}.
\end{align} 
In the above, we used the change of the variables $s=th(|x-y|)$ in the third line, the fact that $h(r) \ge L(r)$ for all $r>0$ in the fourth line, and Lemma \ref{l:green1} and \eqref{e:asym2} in the fifth line.

\smallskip

{\bf(Upper bound)} 
Using boundary Harnack principle and Lemma \ref{intgreen}, one can prove  the upper bound following
 the proofs of \cite[Theorem 1.2 and Theorem 6.4]{KM14} and \cite[Theorem 4.6]{KSV12b} line by line.
 Thus, we provide the main steps of the proof only. 
 
By the boundary Harnack principle (see, \cite[Theorem 1.9]{GK18}), Lemma \ref{intgreen} and \eqref{green_lower}, we can follow the proof of \cite[Theorem 6.4]{KM14} to obtain
\begin{align}\label{e:BHPdecomp}
G_D(x,y)\le c\frac{ g_D(x)g_D(y)}{g_D(A)^2}\frac{\nu(|x-y|)}{h(|x-y|)^2},
\end{align} 
where $g_D(z):=G_D(z, z_0)\wedge c_1$ for some fixed constant $c_1>0$, $z_0\in D$ is a fixed point in $D$ and $A\in \sB(x,y)$, where $\sB(x,y)$ is given by \cite[(6.7)]{KM14}. Moreover, we can also follow the proof of \cite[Theorem 4.6]{KSV12b} to show that for all $z\in D$,
\begin{align}\label{e:gDesti}
g_D(z)\asymp L(\delta_D(z))^{-1/2}.
\end{align} 
Indeed, let $R_3:=\delta_D(z_0) \wedge R_2$ where $R_2$ is the constant in Lemma \ref{l:explictdecay}. 
If $\delta_D(z) \ge R_3/8$, then we get $L(\delta_D(z))^{-1/2} \ge L(R_3/8)^{-1/2} \ge c_1^{-1} L(R_3/8)^{-1/2}g_D(z)$. Moreover, by \eqref{green_lower} and Lemma \ref{l:green2}, we also get $g_D(z) \ge c \ge cL(r_2)^{-1/2} \ge cL(\delta_D(z))^{-1/2}$. Hence, \eqref{e:gDesti} holds in this case.

Next, we assume that $\delta_D(z) <R_3/8$. Then, we get that $|z-z_0| \ge \delta_D(z_0)-\delta_D(z) \ge 7R_3/8$. Therefore, by Lemma \ref{l:green2}, $g_D(z) \asymp G_D(z,z_0)$. Choose $w_z \in \partial D$ satisfying $\delta_D(z)=|z-w_z|$. Let $z^*:=w_z + R_3(z-w_z)/(4|z-w_z|) \in D$ and define $U(z,1)$ as \eqref{e:dUW}. Then, by the boundary Harnack principle, \eqref{e:exitdist}, Lemma \ref{intgreen}, \eqref{green_lower} and Proposition \ref{p:survival}, we get
$$
g_D(z) \asymp G_D(z,z_0) \asymp G_D(z^*,z_0) \frac{\P_z\big(Y_{\tau_{U(z,1)}} \in D\big)}{\P_{z^*}\big(Y_{\tau_{U(z,1)}} \in D\big)} \asymp \P_z\big(Y_{\tau_{U(z,1)}} \in D\big) \asymp L(\delta_D(z))^{-1/2}.
$$
 Hence, we obtain \eqref{e:gDesti}.

We see from the definition of $\sB(x,y)$ that $\delta_D(A) \ge c|x-y|$ for some constant $c>0$. Thus,
by combining \eqref{e:BHPdecomp} and \eqref{e:gDesti}, we get from \eqref{e:Lws} that 
$$
G_D(x,y) \le ca(x,y)L(\delta_D(A))\frac{\nu(|x-y|)}{h(|x-y|)^2} \le ca(x,y)L(|x-y|)\frac{\nu(|x-y|)}{h(|x-y|)^2}.
$$
This together with Lemma \ref{intgreen} completes the proof.  \qed

\section{Example}

In this section, we give an example that is covered by our results.

\begin{example}\label{ex:logp}
{\rm Let $Y=(Y_t : t \ge 0)$ be a pure jump isotropic unimodal L\'evy process with L\'evy measure $\nu$ satisfying {\bf (A)} and {\bf (B)}, and $D$ be a $C^{1,1}$ open set in $\R^d$ with characteristics $(R_0,\Lambda)$. Suppose that there exists $p \in [-1, \infty)$ such that
\begin{align}\label{logp}
\nu(r) \asymp r^{-d} |\log r|^p \qquad \text{for} \;\; 0<r \le 1/2.
\end{align}
Typical examples of isotropic unimodal L\'evy processes satisfying \eqref{logp} are geometric stable processes $(p=0)$ and iterated geometric stable processes $(p=-1)$. The condition $p \ge -1$ is necessary to make the L\'evy measure $\nu$ be infinite. We let
\begin{align*}
\lo(r):= \log (e + r) \quad \text{for} \;\; r>0.
\end{align*}
Then, for every fixed $R>0$, we have that for $0<r \le R$,
\begin{align*}
\ell(r) \asymp \lo(r^{-1})^p \quad \text{and} \quad L(r) \asymp h(r) \asymp \begin{cases}
\lo(r^{-1})^{p+1}, & \mbox{if } \;\; p>-1 ; \\
\lo \circ \lo (r^{-1}), & \mbox{if } \;\; p=-1.
\end{cases}
\end{align*}

We first obtain the small time estimates for the Dirichlet heat kernel. Define for $p>-1$,
\begin{equation*}
\mathfrak{B}_p(x,y):=\left(1 \wedge \frac{\lo(\delta_D(x)^{-1})^{-(p+1)/2}}{\sqrt{t}}\right)\left(1 \wedge \frac{\lo(\delta_D(y)^{-1})^{-(p+1)/2}}{\sqrt{t}}\right).
\end{equation*}
and
\begin{equation*}
\mathfrak{B}_{-1}(x,y):=\left(1 \wedge \frac{[\lo \circ \lo(\delta_D(x)^{-1})]^{-1/2}}{\sqrt{t}}\right)\left(1 \wedge \frac{[\lo \circ \lo(\delta_D(y)^{-1})]^{-1/2}}{\sqrt{t}}\right).
\end{equation*}

\noindent{(Case 1) $p>0.$}

In this case, {\bf (S-2)} holds. Note that we do not need the condition {\bf (C)} when we estimate $p_D(t,x,y)$ only for $|x-y| \le 1$. Thus, according to Theorem \ref{t:main}(ii), for every $T>0$, there are constants $c_0>1$, $c_1,...,c_6>0$ such that for all $(t,x,y) \in (0,T] \times D \times D$ satisfying $|x-y| \le 1$,
\begin{align*}
&c_0^{-1}\mathfrak{B}_p(x,y) F_p(t,|x-y|,c_1,c_2,c_3)\le p_D(t,x,y) \le c_0\mathfrak{B}_p(x,y) F_p(t,|x-y|,c_4,c_5,c_6),
\end{align*} 
where
\begin{align*}
F_p(t,r,a_1,a_2,a_3):=
\begin{cases}
\exp\big(a_1t^{-1/p}\big), & \mbox{if } \, r \le \exp\big(-a_3t^{-1/p}\big); \\[3pt]
t \lo(r^{-1})^p  \exp \Big( \big(-d+a_2t \lo(r^{-1})^p\big)\log r \Big), & \mbox{if } \, r >\exp\big(-a_3t^{-1/p}\big).
\end{cases}
\end{align*}

\smallskip

\noindent{(Case 2) $-1<p \le 0.$}
  
Since {\bf (S-1)} holds, by Theorem \ref{t:main}(i), for every $T>0$, there are constants $c_0>1$, $c_1,c_2>0$ such that for all $(t,x,y) \in (0,T] \times (D \times D \setminus \diag)$ satisfying $|x-y| \le 1$,
\begin{align*}
&c_0^{-1}\mathfrak{B}_p(x,y)
t\lo(|x-y|^{-1})^p\exp\Big(\big(-d+c_1t\lo(|x-y|^{-1})^p\big) \log|x-y|\Big)  \\[3pt]
&\le p_D(t,x,y)\le c_0\mathfrak{B}_p(x,y)t\lo(|x-y|^{-1})^p\exp\Big(\big(-d+c_2t\lo(|x-y|^{-1})^p\big) \log|x-y|\Big).
\end{align*}

\smallskip

\noindent{(Case 3) $p = -1.$}
  
Since {\bf (S-1)} holds, by Theorem \ref{t:main}(i), for every $T>0$, there are constants $c_0>1$, $c_1,c_2>0$ such that for all $(t,x,y) \in (0,T] \times (D \times D \setminus \diag)$ satisfying $|x-y| \le 1$,
\begin{align*}
&c_0^{-1}\mathfrak{B}_{-1}(x,y) t|x-y|^{-d} \lo(r^{-1})^{-1-c_1t}\\[3pt]
&\le p_D(t,x,y)  \le c_0\mathfrak{B}_{-1}(x,y) t|x-y|^{-d} \lo(r^{-1})^{-1-c_2t}.
\end{align*}

\vspace{6mm}

Now, we further assume that $D$ is bounded and of scale $(r_1, r_2)$. Then, we get the following large time estimates.

\vspace{3mm}

\noindent{(Case 1) $p \ge 0.$}

Since either {\bf (S-2)} or {\bf (L-2)} holds, by Theorem \ref{t:main2}(ii, iii), there exist $T_1 \ge 0$ (if $p>0$, then $T_1 = 0$) and $\lambda_1>0$ such that for every  $T>T_1$, we have that for all $(t,x,y) \in [T,\infty) \times D \times D$,
\begin{align*}
p_D(t,x,y) \asymp e^{-\lambda_1 t}\lo(\delta_D(x)^{-1})^{-(p+1)/2}\lo(\delta_D(y)^{-1})^{-(p+1)/2}.
\end{align*}

\smallskip

\noindent{(Case 2) $-1< p < 0$.}

Since  {\bf (L-1)} holds, by Theorem \ref{t:main2}(i), for every $T>0$, there are constants $c_0,c_1,c_2, \lambda_2, \lambda_3>0$ such that for all $(t,x,y) \in [T,\infty) \times (D \times D \setminus \diag)$,
\begin{align*}
&c_0^{-1} \lo(\delta_D(x)^{-1})^{-(p+1)/2}\lo(\delta_D(y)^{-1})^{-(p+1)/2}\left(|x-y|^{-d+c_1t\lo(|x-y|^{-1})^p} \lo(|x-y|^{-1})^p  + e^{- \lambda_2 t}\right) \\
& \le p_D(t,x,y) \\
& \le c_0 \lo(\delta_D(x)^{-1})^{-(p+1)/2}\lo(\delta_D(y)^{-1})^{-(p+1)/2}\left(|x-y|^{-d+c_2t\lo(|x-y|^{-1})^p} \lo(|x-y|^{-1})^p  + e^{- \lambda_3 t}\right).
\end{align*}

\smallskip

\noindent{(Case 3) $p=-1$.}

Since {\bf (L-1)}  holds, by Theorem \ref{t:main2}(i), for every $T>0$, there are constants $c_0,c_1,c_2, \lambda_2, \lambda_3>0$ such that for all $(t,x,y) \in [T,\infty) \times (D \times D \setminus \diag)$,
\begin{align*}
&c_0^{-1} \big[\lo \circ \lo(\delta_D(x)^{-1})\big]^{-1/2}\big[\lo \circ \lo(\delta_D(y)^{-1})\big]^{-1/2}\Big( |x-y|^{-d} \lo(r^{-1})^{-1-c_1t}  + e^{- \lambda_2 t}\Big) \\
& \le p_D(t,x,y) \\
& \le c_0  \big[\lo \circ \lo(\delta_D(x)^{-1})\big]^{-1/2}\big[\lo \circ \lo(\delta_D(y)^{-1})\big]^{-1/2}\Big( |x-y|^{-d} \lo(r^{-1})^{-1-c_2t} + e^{- \lambda_3 t}\Big).
\end{align*}

\vspace{2mm}

Finally, we obtain the Green function estimates by Theorem \ref{t:green}. Let $D$ be a bounded $C^{1,1}$ open set in $\R^d$. Note that the L\'evy measure $\nu$ satisfies {\bf (D)}. Hence, we have that, for all $x,y \in D$, if $p>-1$, then
\vspace{2mm}
\begin{align*}
G_D(x,y) \asymp \left(1 \wedge \frac{\lo(\delta_D(x)^{-1})^{-(p+1)/2}\lo(\delta_D(y)^{-1})^{-(p+1)/2}}{\lo(|x-y|^{-1})^{-(p+1)}}\right)  \frac{\lo(|x-y|^{-1})^{-(p+2)}}{|x-y|^d}
\end{align*}
\vspace{2mm}
and if $p=-1$, then
\begin{align*}
G_D(x,y) \asymp &\left(1 \wedge \frac{[\lo \circ \lo(\delta_D(x)^{-1})]^{-1/2}[\lo \circ \lo(\delta_D(y)^{-1})]^{-1/2}}{[\lo \circ \lo(|x-y|^{-1})]^{-1}}\right) \\
& \;\; \times  \frac{\lo(|x-y|^{-1})^{-1}[\lo \circ \lo(|x-y|^{-1})]^{-2}}{|x-y|^d}.
\end{align*}
}
\end{example}


\vspace{.1in}


\begin{thebibliography}{99}
\bibitem{BK17}
J.~Bae and P.~Kim.
\newblock On estimates of transition density for subordinate Brownian motions with Gaussian components in $C^{1,1}$-open sets.
\newblock {\em Potential Anal.}. 52(4):661--687, 2020.
 
\bibitem{Be96}
J.~Bertoin.
\newblock {\em L\'{e}vy processes}, 
\newblock Cambridge University Press, Cambridge, 1996.



\bibitem{BGT}
N.~H. Bingham, C.~M. Goldie, and J.~L. Teugels.
\newblock {\em Regular variation}, 
\newblock Cambridge University Press, Cambridge, 1989.

\bibitem{BBKRSV}
K. Bogdan, T. Byczkowski, T. Kulczycki, M. Ryznar, R. Song and Z. ~Vondra\v{c}ek.
{\em Potential analysis of stable processes and its extensions}. Lecture Notes in Mathematics, 1980. Springer-Verlag, Berlin, 2009. 

\bibitem{BGPR}
K.~Bogdan, T.~Grzywny, K.~Pietruska-Pałuba and A. Rutkowski. 
\newblock Extension and trace for nonlocal operators. 
\newblock {\em J. Math. Pures Appl.}, 137: 33--69, 2020.


\bibitem{BGR10}
K.~Bogdan, T.~Grzywny, and M.~Ryznar.
\newblock Heat kernel estimates for the fractional {L}aplacian with {D}irichlet
  conditions.
\newblock {\em Ann. Probab.}, 38(5):1901--1923, 2010.

\bibitem{BGR14a}
K.~Bogdan, T.~Grzywny, and M.~Ryznar.
\newblock Density and tails of unimodal convolution semigroups.
\newblock {\em J. Funct. Anal.}, 266(6):3543--3571, 2014.

\bibitem{BGR14b}
K.~Bogdan, T.~Grzywny, and M.~Ryznar.
\newblock Dirichlet heat kernel for unimodal {L}\'{e}vy processes.
\newblock {\em Stochastic Process. Appl.}, 124(11):3612--3650, 2014.

\bibitem{BGR15}
K.~Bogdan, T.~Grzywny, and M.~Ryznar.
\newblock Barriers, exit time and survival probability for unimodal {L}\'{e}vy
  processes.
\newblock {\em Probab. Theory Related Fields}, 162(1-2):155--198, 2015.

\bibitem{BSW13}
B. Böttcher, R. L. Schilling and J. Wang.
\newblock {\em L\'{e}vy matters. III.}
\newblock Lecture Notes in Mathematics, 2099. Lévy Matters. Springer, Cham, 2013. 

\bibitem{Ch09}
Z.-Q. Chen.
\newblock On notions of harmonicity.
\newblock {\em Proc. Amer. Math. Soc.}, 137(10):3497--3510, 2009.

\bibitem{CK16}
Z.-Q. Chen and P.~Kim.
\newblock Global {D}irichlet heat kernel estimates for symmetric {L}{\'e}vy
  processes in half-space.
\newblock {\em Acta Appl. Math.}, 146:113--143, 2016.


\bibitem{CKKW18} 
Z.-Q. Chen, P.~Kim, T.~Kumagai and J.~Wang.
\newblock Heat kernel estimates for time fractional equations.
\newblock {\em Forum Math.}, 30(5):1163--1192, 2018.

\bibitem{CKS10}
Z.-Q. Chen, P.~Kim, and R.~Song.
\newblock Heat kernel estimates for the {D}irichlet fractional {L}aplacian.
\newblock {\em J. Eur. Math. Soc.}, 12(5):1307--1329, 2010.


\bibitem{CKS11}
Z.-Q. Chen, P.~Kim, and R.~Song.
\newblock Heat kernel estimates for {$\Delta+\Delta^{\alpha/2}$} in {$C^{1,1}$}
  open sets.
\newblock {\em J. Lond. Math. Soc. (2)}, 84(1):58--80, 2011.


\bibitem{CKS12a}
Z.-Q. Chen, P.~Kim, and R.~Song.
\newblock Sharp heat kernel estimates for relativistic stable processes in open
  sets.
\newblock {\em Ann. Probab.}, 40(1):213--244, 2012.



\bibitem{CKS12b}
Z.-Q. Chen, P.~Kim, and R.~Song.
\newblock Global heat kernel estimates for relativistic stable processes in
  half-space-like open sets.
\newblock {\em Potential Anal.}, 36(2):235--261, 2012.


\bibitem{CKS14}
Z.-Q. Chen, P.~Kim, and R.~Song.
\newblock Dirichlet heat kernel estimates for rotationally symmetric {L}\'{e}vy
  processes.
\newblock {\em Proc. Lond. Math. Soc. (3)}, 109(1):90--120, 2014.


\bibitem{CKS16}
Z.-Q. {Chen}, P.~{Kim}, and R.~{Song}.
\newblock {Dirichlet heat kernel estimates for subordinate Brownian motions
  with Gaussian components.}
\newblock {\em {J. Reine Angew. Math.}}, 711:111--138, 2016.


\bibitem{CT11}
Z.-Q. Chen and J.~Tokle.
\newblock Global heat kernel estimates for fractional Laplacians in unbounded
  open sets.
\newblock {\em Probab. Theory and Related Fields}, 149(3):373--395,   2011.

\bibitem{CKSV20}
S. Cho, P. Kim, R. Song and Z. Vondra\v{c}ek.
\newblock Factorization and estimates of Dirichlet heat kernels for non-local operators with critical killings
\newblock {\em 	J. Math. Pures Appl.}, 143:208--256, 2020.


\bibitem{GKK20}
T.~Grzywny, K.-Y. Kim, and P.~Kim.
\newblock Estimates of {D}irichlet heat kernel for symmetric {M}arkov
  processes.
\newblock {\em Stochastic Process. Appl.}, 130(1):431--470, 2020.



\bibitem{GK18}
T.~Grzywny and M.~Kwa\'{s}nicki.
\newblock Potential kernels, probabilities of hitting a ball, harmonic
  functions and the boundary {H}arnack inequality for unimodal {L}\'{e}vy
  processes.
\newblock {\em Stochastic Process. Appl.}, 128(1):1--38, 2018.

\bibitem{GR}
T.~Grzywny and M.~Ryznar, \newblock{Potential theory of one-dimensional geometric
	stable processes}.
\newblock {\em  Colloq. Math.}, 129(1):7--40, 2012.


\bibitem{GRT19}
T.~Grzywny, M.~Ryznar, and B.~Trojan.
\newblock Asymptotic behaviour and estimates of slowly varying convolution
  semigroups.
\newblock {\em Int. Math. Res. Not. IMRN}, (23):7193--7258, 2019.

\bibitem{HW42}
P.~Hartman and A.~Wintner.
\newblock On the infinitesimal generators of integral convolutions.
\newblock {\em Amer. J. Math.}, 64:273--298, 1942.

\bibitem{KK14}
K.-Y. Kim and P.~Kim.
\newblock Two-sided estimates for the transition densities of symmetric
  {M}arkov processes dominated by stable-like processes in {$C^{1,\eta}$} open
  sets.
\newblock {\em Stochastic Process. Appl.}, 124(9):3055--3083, 2014.


\bibitem{KM12}
P.~Kim and A.~Mimica.
\newblock Harnack inequalities for subordinate {B}rownian motions.
\newblock {\em Electron. J. Probab.}, 17(37):1--23, 2012.


\bibitem{KM14}
P.~Kim and A.~Mimica.
\newblock Green function estimates for subordinate Brownian motions: Stable and
  beyond.
\newblock {\em Trans. Amer. Math. Soc.},
  366(8):4383--4422, 2014.


\bibitem{KM17}
P.~Kim and A.~Mimica.
\newblock Estimates of {D}irichlet heat kernels for subordinate {B}rownian
  motions.
\newblock {\em Electron. J. Probab.}, 23(64):1--45, 2018.


\bibitem{KSV12}
P.~Kim, R.~Song, and Z.~Vondra\v{c}ek.
\newblock Potential theory of subordinate {B}rownian motions revisited.
\newblock {\em Stochastic analysis and applications to finance}, 243–290, Interdiscip. Math. Sci., 13, World Sci. Publ., Hackensack, NJ, 2012.


\bibitem{KSV12b}
P.~Kim, R.~Song, and Z.~Vondra\v{c}ek.
\newblock Two-sided Green function estimates for killed subordinate Brownian motions.
\newblock {\em Proc. Lond. Math. Soc. (3)},  104(5):927--958, 2012. 

\bibitem{KSV14}
P.~Kim, R.~Song, and  Z.~Vondra\v{c}ek.
\newblock Global uniform boundary Harnack principle with explicit decay rate and its application.
\newblock {\em 
Stochastic Process. Appl.}, 124(1):235--267, 2014.

\bibitem{KS13}
V.~Knopova and R.~L. Schilling.
\newblock A note on the existence of transition probability densities of
  {L}\'evy processes.
\newblock {\em Forum Math.}, 25(1):125--149, 2013.

\bibitem{KuR16}
T.~Kulczycki and M.~Ryznar.
\newblock Gradient estimates of harmonic functions and transition densities for
  {L}{{\'e}}vy processes.
\newblock {\em Trans. Amer. Math. Soc.}, 368(1):281--318, 2016.

\bibitem{M16}
A.~Mimica.
\newblock Heat kernel estimates for subordinate {B}rownian motions.
\newblock {\em Proc. Lond. Math. Soc. (3)}, 113(5):627--648, 2016.

\bibitem{Pr81}
W.~E. Pruitt.
\newblock The growth of random walks and L\'evy processes.
\newblock {\em Ann. Probab.}, 9(6):948--956, 1981.

\bibitem{Si80}
M.~L. Silverstein.
\newblock Classification of coharmonic and coinvariant functions for a
  {L}\'{e}vy process.
\newblock {\em Ann. Probab.}, 8(3):539--575, 1980.

\bibitem{W83}
T.~Watanabe.
\newblock The isoperimetric inequality for isotropic unimodal {L}{\'e}vy
  processes.
\newblock {\em Z. Wahrsch. Verw. Gebiete}, 63(4):487--499, 1983.

\end{thebibliography}
\end{document}